\documentclass{amsart}
\usepackage{amsmath,amssymb}

\numberwithin{equation}{section}

\newtheorem{theorem}[equation]{Theorem}
\newtheorem{proposition}[equation]{Proposition}
\newtheorem{lemma}[equation]{Lemma}
\newtheorem{corollary}[equation]{Corollary}

\theoremstyle{definition}
\newtheorem{definition}[equation]{Definition}
\newtheorem{remark}[equation]{Remark}
\newtheorem{example}[equation]{Example}

\def\C{\mathbb C}
\def\Dom{\mathcal{D}}
\def\R{\mathbb R}
\def\ceil#1{\lceil#1\rceil}

\def\embed{\hookrightarrow}
\def\eps{\varepsilon}
\def\im{i}
\def\seq#1{\{#1\}_{n\in\mathbb N}}
\def\set#1{\{#1\}}
\def\st{\,|\,}
\def\vt{\vartheta}
\def\Wedge{\raise2ex\hbox{$\mathchar"0356$}}
\def\minus{\backslash}
\def\Dot#1{\overset{\,{}_\circ}{#1}}
\DeclareMathOperator{\Ind}{ind}
\DeclareMathOperator{\LinSpan}{span}
\DeclareMathOperator{\diff}{Diff}
\DeclareMathOperator{\rg}{rg}
\DeclareMathOperator{\spec}{spec}

\begin{document}

\title{Adjoints of elliptic cone operators}
\author{Juan B. Gil}
\author{Gerardo A. Mendoza}
\address{Department of Mathematics\\ Temple University\\
 Philadelphia, PA 19122}

\begin{abstract}
We study the adjointness problem for the closed extensions of a general
$b$-elliptic operator $A\in x^{-\nu}\diff^m_b(M;E)$, $\nu>0$, initially defined as an
unbounded operator $A:C_c^\infty(M;E)\subset x^\mu L^2_b(M;E)\to x^\mu L^2_b(M;E)$,
$\mu \in \R$. The case where $A$ is a symmetric semibounded operator is of particular
interest, and we give a complete description of the domain of the Friedrichs
extension of such an operator.
\end{abstract}

\maketitle


\section{Introduction}

Let $M_0$ be a smooth paracompact manifold, $\mathfrak m$ a smooth positive
measure on $M_0$. Suppose $A:C_c^\infty(M_0)\to C_c^\infty(M_0)$ is
a scalar linear partial differential operator with smooth coefficients.
Among all possible domains $\Dom\subset L^2(M_0,\mathfrak m)$ for $A$ as an
unbounded operator on $L^2(M_0)$ there are two that stand out:
\begin{equation*}
 \Dom_{\max}(A) = \set{u\in L^2(M_0,\mathfrak m)\st Au\in L^2(M_0,\mathfrak m)}
\end{equation*}
where $Au$ is computed in the distributional sense, and
\begin{equation*}
 \Dom_{\min}(A)= \text{ completion of } C_c^\infty(M_0) \text{ with
 respect to the norm } \|u\|+\|Au\|,
\end{equation*}
which can be regarded as a subspace of $L^2(M_0,\mathfrak m)$. Both domains are
dense in $L^2(M_0,\mathfrak m)$, since $C_c^\infty(M_0)\subset\Dom_{\min}(A)\subset
\Dom_{\max}(A)$, and with each domain $A$ is a closed operator. Clearly
$\Dom_{\min}(A)$ is the smallest domain containing $C_c^\infty(M_0)$ with respect to
which $A$ is closed, and $\Dom_{\max}(A)$ contains any domain on which the action of
the operator coincides with the action of $A$ in the distributional sense.
Also clearly, $A$ with domain $\Dom_{\max}(A)$ is an extension of
$A$ with domain $\Dom_{\min}(A)$.

If $M_0$ is compact without boundary and $A$ is elliptic then $\Dom_{\max}(A) =
\Dom_{\min}(A)$; the interesting situations occur when $A$ is non-elliptic or $M_0$
is noncompact. In this paper we shall analyze the latter problem, assuming that
$M_0$ is the interior of a smooth compact manifold $M$ with boundary and  that $A \in
x^{-\nu}\diff^m_b(M;E)$,
$\nu>0$, is a $b$-elliptic `cone' operator acting on sections of a smooth
vector bundle $E\to M$; here $x:M\to \R$ is a smooth defining function for
$\partial M$, positive in $M_0$.

The elements of $\diff^m_b(M;E)$ are the totally characteristic differential
operators introduced and analyzed systematically by Melrose \cite{RBM1}. These are
linear operators with smooth coefficients which near the boundary can be written in
local coordinates $(x,y_1,\dots,y_n)$ as $P=\sum_{k+|\alpha|\leq m}
a_{k\alpha}(xD_x)^kD_y^\alpha$. Such an operator is $b$-elliptic if it is elliptic in
the interior in the usual sense and in addition $\sum_{k+|\alpha|= m}
a_{k\alpha}(0,y)\xi^k\eta^\alpha$ is invertible for $(\xi,\eta)\ne 0$; this is
expressed more concisely by saying that the principal symbol of $P$, as an object on
the compressed cotangent bundle (see Melrose, op.cit.), is an isomorphism. It follows
from the definition of $b$-ellipiticity that the family of differential operators on
$\partial M$ given locally by
$\hat P_0(\sigma)=\sum_{j+|\alpha|\leq m} a_{\alpha,j}(0,y)\sigma^jD_y^\alpha$,
called the indicial operator or conormal symbol of $P$, is elliptic of
order $m$ for any $\sigma\in \C$. For the general theory of these operators and the
associated pseudodifferential calculus the reader is referred to the paper of Melrose
cited above, his book
\cite{RBM2}, as well as Schulze \cite{Sz91,Sz94}. An operator $A\in
x^{-\nu}\diff^m_b(M;E)$ is $b$-elliptic if $P=x^{\nu}A$ is $b$-elliptic. By
definition the conormal symbol of $A$ is that of $P$. For more details see Section
\ref{SetUp}.

The measure used to define the $L^2$ spaces will be
of the form $x^\mu\mathfrak m$ for some real $\mu$, where $\mathfrak m$ is a
$b$-density, that is, $x\mathfrak m$ is a smooth positive density. The spaces
$L^2(M;E;x^\mu\mathfrak m)$ are defined in the usual way with the aid of a smooth but
otherwise arbitrary hermitian metric on $E$. These spaces are related among themselves
by canonical isometries with which $L^2(M;E;x^\mu\mathfrak m) = x^{-\mu/2}L^2_b(M,E)$,
where
$L^2_b(M,E)$ is the space defined by the measure $\mathfrak m$ itself.

This said, the general problem we are concerned with is the description of the
adjoints of the closed extensions of a general $b$-elliptic operator $A\in
x^{-\nu}\diff^m_b(M;E)$ initially defined as an unbounded operator 
 \begin{equation}\label{InitialOp}
 A:C_c^\infty(M;E)\subset x^\mu L^2_b(M;E)\to x^\mu L^2_b(M;E).
 \end{equation}
 The case where $A$ is a symmetric semibounded operator is of particular interest, and we
give, in Theorem \ref{FriedExtension}, a complete description of the domain of the
Friedrichs extension of such an operator.

Differential operators in $x^{-\nu}\diff^m_b(M;E)$ arise in the study of manifolds
with conical singularities. The study of such manifolds from the
geometric point of view began with Cheeger~\cite{Ch79}, and by now there is an
extensive literature on the subject. In the specific context of our work, probably
the most relevant references, aside from those cited above, are the book by
Lesch \cite{Le97}, and the papers by Br\"uning and Seeley \cite{BrSe85,BrSe88}, and
Mooers \cite{Moo99}. See also Coriasco, Schrohe, and Seiler \cite{CSS01}. 

As already mentioned, the domains $\Dom_{\min}(A)$ and $\Dom_{\max}(A)$ need not be
the same. The object determining the closed extensions of $A$
is its conormal symbol $\hat P_0(\sigma):C^\infty(\partial M)\to C^\infty(\partial M)$.
Because of the $b$-ellipticity, this operator is invertible for all $\sigma\in \C$
except a discrete set $\spec_b(A)$, the boundary spectrum of $A$ (or $P$), a set which
(again due to the $b$-ellipticity) intersects any strip $a<\Im \sigma < b$ in a
finite set.  It was noted by Lesch
\cite{Le97} that $\Dom_{\min}(A) = \Dom_{\max}(A)$ if and only if
$\spec_b(A)\cap\set{-\mu-\nu <\Im\sigma<-\mu}=\varnothing$. Also proved by
Lesch [op.cit., Proposition 1.3.16] was the fact that
$\Dom_{\max}(A)/\Dom_{\min}(A)$ is finite dimensional. This provides a
simple description of the closed extensions of $A$: they are given by the
operator $A$ acting in the distributional sense on subspaces $\Dom\subset
x^{\mu}L^2_b(M;E)$ with $\Dom_{\min}(A)\subset \Dom\subset \Dom_{\max}(A)$. These
results and others, also due to Lesch (op.cit.) are reproved in Section
\ref{SectGround} for the sake of completeness, with a slightly different approach
emphasizing the use of the pseudodifferential calculus, for totally characteristic
operators \cite{MM}, or for the cone algebra \cite{Sz91}. In that section we also
prove a relative index theorem for the closed extensions of $A$.

 From the point of view of closed extensions there is nothing more to understand than
that they are in one to one correspondence with the subspaces of
$\mathcal E(A) = \Dom_{\max}(A)/\Dom_{\min}(A)$. The problem of finding the domain of
the adjoint is more delicate, and forces us to pass to the Mellin transforms of
representatives of elements of $\mathcal E$. Our approach entails rewriting the
pairing $[u,v]_A=(Au,v)-(u,A^\star v)$, defined for $u\in \Dom_{\max}(A)$ and
$v\in \Dom_{\max}(A^\star)$, where $A^\star$ is the formal adjoint of
$A$, in terms of the Mellin transforms of $u$ and $v$, and proving certain specific
nondegeneracy properties of the pairing. It is generally true (and well known) that
in a general abstract setting, $[\cdot,\cdot]_A$ induces a nonsingular pairing
of $\mathcal E(A)$ and $\mathcal E(A^\star)$. In the case at hand, there is an
essentially well defined notion of an element $u\in \Dom_{\max}(A)$ with pole `only'
at $\sigma_0\in \spec_b(A)\cap \set{-\mu-\nu<\Im\sigma<-\mu}$. If $\sigma_0\in 
\spec_b(A)$ then $\overline \sigma_0 - \im(\nu+2\mu)\in \spec_b(A^\star)$, and we
show, for example, that the restriction of the pairing to elements $u\in
\Dom_{\max}(A)$ with pole `only' at $\sigma_0$ and elements $v$ in
$\Dom_{\max}(A^\star)$ with pole `only' at $\overline \sigma_0 - \im(\nu+2\mu)$ is
nonsingular (modulo the respective minimal domains). 

Our analysis of the pairing begins in Section~\ref{MeromorphicSolutions} with a
careful description of the Mellin transforms of elements in
$\Dom_{\max}(A)$. The main result in that section is Proposition~\ref{mero20}, a
result along the lines of part of the work of Gohberg and Sigal \cite{GohSi}. In
Section~\ref{CanonicalPairing} we prove in a somewhat abstract setting that the
pairing alluded to above for solutions with poles at conjugate points is nonsingular
(Theorem~\ref{mero5}).  In
Section~\ref{DomainPairing} we link the pairing
$[\cdot,\cdot]_A$ with the pairing of Section~\ref{CanonicalPairing}. The main
results are Theorems~\ref{BasicPairing} and \ref{FlatPairing}. The first of these
gives a formula for the pairing which in particular shows that the pairing is null
when the poles in question are not conjugate, and the second, which is based on the
formula in Theorem~\ref{BasicPairing} and Theorem~\ref{mero5}, shows that the
pairing of elements associated to conjugate poles is nonsingular. The formulas are
explicit enough that in simple cases it is easy to determine the domain of the
adjoint of a given extension of $A$. 

We undertake the study of the domain of the Friedrichs extension of $b$-elliptic
semibounded operators in Section~\ref{FriedrichsExtension}. The main result there
is Theorem~\ref{FriedExtension}, which is a complete description of the domain of the
Friedrichs extension. Loosely speaking, the domain consists of the sum of those
elements $u\in \Dom_{\max}(A)$ with poles `only' in $-\mu-\nu<\Im\sigma<-\mu -\nu/2$
and those with pole `only' at $\Im\sigma=-\mu -\nu/2$ and `half' of the order of the
pole.  Finally, in Section~\ref{Examples} we collect a number of examples that
illustrate the use of Theorems ~\ref{BasicPairing},
\ref{FlatPairing}, and \ref{FriedExtension}.

Any closed extension of a $b$-elliptic operator $A\in x^{-\nu}\diff^m_b(M;E)$ as
an operator $\Dom\subset x^{\mu}L^2_b(M;E)\to x^{\mu}L^2_b(M;E)$ has the space
$x^{\mu+\nu} H^m_b(M;E)$ in its domain. The space $ H^m_b(M;E)$ is the subspace of
$L^2_b(M;E)$ whose elements $u$ are such that for any smooth vector fields
$V_1,\dotsc,V_k$, $k\leq m$, on $M$ tangent to the boundary,
$V_1\dots V_k u \in L^2_b(M;E)$. Other than this, the set of domains of the closed
extensions of different $b$-elliptic operators in $x^{-\nu}\diff^m_b(M;E)$ are
generally not equal. In Section \ref{EqOfDom} we provide simple sufficient conditions
for the domains of two such operators to be the same. Not unexpectedly, these
conditions are on the equality of Taylor expansion at the boundary of the operators
involved, up to an order depending on $\nu$. We prove in particular that under the
appropriate condition the domains of the Friedrichs extensions of two different
symmetric semibounded operators are the same. This is used in Section
\ref{FriedrichsExtension} as an intermediate step to determine the domain of the
Friedrichs extension of such an operator, and a refinement of the condition is
obtained as a consequence.

Operators of the kind we investigate arise naturally as geometric operators. In such
applications $\mu$ is determined by the actual situation. It is convenient for us,
however, to work with the normalization
 \begin{equation*}
 x^{-\mu-\nu/2}A
x^{\mu+\nu/2}:C_c^\infty(M;E)\subset x^{-\nu/2} L^2_b(M;E)\to x^{-\nu/2}
L^2_b(M;E).
  \end{equation*}
 rather than \eqref{InitialOp}. Since the mappings 
$x^s:x^{\mu} L^2_b(M;E)\to x^{\mu+s} L^2_b(M;E)$
 are surjective isometries, this represents no loss. In particular, adjoints and
symmetry properties of operators are preserved. Also to be noted is that these
transformations represent translations on the Mellin transform side, so it is a
simple exercise to recast information presented in terms of Mellin transforms of the
modified operator as information on the original operator.

 \section{Geometric preliminaries}\label{SetUp}
 
 Throughout the paper $M$ is a compact manifold with (nonempty) boundary
with a fixed positive $b$-density $\mathfrak m$, that is, a smooth
density
$\mathfrak m$ such that for some (hence any) defining
function $x$, $x\mathfrak m$ is a smooth positive density. We
will also fix a hermitian vector bundle $E\to M$.  

Fix a collar neighborhood $U_Y$ for each boundary component $Y$ of $M$, 
so we have a trivial fiber bundle $\pi_Y:U_Y\to Y$ with fiber
$[0,1)$. We can then canonically identify the bundle of
$1$-densities over $U_Y$  with  $|\Wedge|[0,1)\otimes |\Wedge |Y$ (the
tensor product of the pullback to $[0,1) \times Y$ of the respective
density bundles).  

 Let $\mathfrak m$ be a $b$-density on $|\Wedge
|[0,1)\otimes |\Wedge|Y$. Then there is a smooth defining function
$x:Y\times [0,1)\to
\R$ vanishing at $Y\times {0}$, and a smooth density
$\mathfrak m_Y$ on $Y$, such that 
 \begin{equation*}
 \mathfrak m = \frac {dx} x \otimes \mathfrak m_Y
 \end{equation*}
  Indeed, let $\xi$ be the variable in $[0,1)$. Over the boundary
of $[0,1) \times Y$ we then get, canonically, $\xi \mathfrak m=
d\xi
 \otimes \mathfrak m_Y$. Then, on $[0,1) \times Y$, $\mathfrak m= h
\frac{d\xi}{\xi} \otimes \mathfrak m_Y$ with $h$ smooth, positive,
$h(0,y)=1$. Let $x=g\xi$ where $g$ is determined modulo constant
factor by the requirement that  $\frac {d_\xi x}{x} = h\frac {d_\xi
\xi}{\xi}$ ($d_\xi$ means differential in the variable $\xi$; this
makes sense since we are dealing with a product manifold). Thus
$g(\xi,y)$ should satisfy the equation
 \begin{equation*}
 \frac {\partial g}{\partial \xi} + \frac {1-h}{\xi} g = 0
 \end{equation*}
 Since $h=1$ when $\xi=0$, the solutions $g$ are smooth across $\xi=0$.
Pick the one which is $1$ when $\xi=0$.

We fix the choice of $x$ for each
boundary component. When working near the boundary we will always
assume that the defining function was chosen above, and
that the coordinates, if at all necessary, are consistent with a choice
of product structure as above. By $\partial_x$ we mean the vector field
tangent to the fibers of $U_Y\to Y$ such that $\partial_x x=1$. 

If $E$, $F\to M$ are (smooth) vector bundles and $P\in \diff^m_b(M;E,F)$
is a $b$-elliptic differential operator, then $E$ and $F$ are
isomorphic. This follows from the fact that the principal symbol of
$P$ is an isomorphism $\pi^*E \to \pi^*F$ where $\pi:
{}^b T^*M\minus 0\to M$ is the projection, and the fact that the
compressed cotangent bundle ${}^b T^*M$ admits a global nonvanishing
section (since it is isomorphic to $T^*M$ and $M$ is a manifold with
boundary). Thus when analyzing $b$-elliptic operators
in $\diff^m_b(M;E,F)$ we may assume $F=E$ (for more on this see
\cite{JM}). 

The Hilbert space structure of the space of sections of $E\to M$ is the
usual one, namely integration with respect to $\mathfrak m$ of the
pointwise inner product in $E$:
 \begin{equation*}
 (u,v)_{L^2_b(M;E)}=\int (u,v)_E\, \mathfrak m \quad\text{ if }u,\ v \in
L^2_b(M;E).
 \end{equation*}
Fix a hermitian connection $\nabla$ on $E$. If $P\in \diff^m_b(M;E)$,
then near a boundary component one can write 
 \begin{equation*}
 P = \sum_{\ell=0}^m P_\ell \circ (\nabla_{xD_x})^\ell
 \end{equation*}
 where the $P_\ell$ are differential operators of order $m-\ell$
(defined on $U_Y$) such that for any smooth function $\phi(x)$ and
section
$u$ of $E$ over $U_Y$, $P_\ell(\phi(x) u) =
\phi(x)P_\ell(u)$, in other words, of order zero in $\nabla_{xD_x}$. $P$
is said to have coefficients independent of $x$ near $Y$ if
$\nabla_{\partial_x}P_k(u) = P_k(\nabla_{\partial_x}u)$ for any smooth
section $u$ of $E$ supported in $U_Y$. By means of parallel transport
along the fibers of
$U_Y\to Y$ one can show that if $P\in \diff^m_b(M;E)$, 
then for any $N$ there are operators $P_k$, $\tilde P_N\in
\diff^m_b(M;E)$ such that
 \begin{equation*}
 P = \sum_{k=0}^N P_k x^k + \tilde P_N x^N
 \end{equation*}
 where $P_k$ has coefficients independent of $x$ near $Y$. If $P_k$ 
has coefficients independent of $x$ near $Y$ then so does its formal
adjoint $P_k^\star$. To see this recall that since the connection is
hermitian, $\partial_x (u,v)_E =
(\nabla_{\partial_x}u,v) + (u,\nabla_{\partial_x}v)$ if $u$ and
$v$ are supported near $Y$, so if they vanish on $Y$ then
 \begin{equation*}
 (\nabla_{\partial_x}u,v)_{L^2_b(M;E)} = -
(u,\nabla_{\partial_x}v)_{L^2_b(M;E)}.
 \end{equation*} 
 One derives the assertion easily from this.
 
Fix $\omega\in C_c^\infty(-1,1)$ real valued, nonnegative and such that
$\omega=1$ in a neighborhood of $0$. The Mellin transform of a section
of $C_c^\infty(\Dot M;E)$ is defined to be the entire function $\hat
u:\C\to C^\infty(Y;E)$ such that for any $v\in C^\infty(Y;E|_Y)$
 \begin{equation}\label{MellinDefinition}
 (x^{-\im\sigma} \omega(x) u,\pi_Y^*v)_{L^2_b(M;E)} = \frac 1 {2\pi}
\int_{\Im\sigma=0}(\hat u(\sigma,y),v(y))_{L^2(Y;E|_Y)}\, d\sigma
 \end{equation}
 By $\pi_Y^*v$ we mean the section of $E$ over $U_Y$ obtained by
parallel transport of $v$ along the fibers of $\pi_Y$. Thus if
$u\in C^\infty(U_Y;E)$ is such that
$\nabla_{\partial_x} u=0$ and $\phi \in C_c^\infty(0,1)$, then 
 \begin{equation*}
 \widehat {\phi u}(\sigma)=  \widehat {\omega\phi}(\sigma)\, u
 \end{equation*}
 where $\widehat {\omega\phi}(\sigma)$ is the ``usual'' Mellin transform
of
$\omega\phi$. The only point here is that we incorporate the cut-off
function into the definition. As is well known, the Mellin transform
extends to the spaces $x^\mu L^2_b(M,E)$ in such a way that if $u \in
x^\mu L^2_b(M,E)$ then $\hat u(\sigma)$ is holomorphic in $\set {\Im
\sigma > -\mu}$ and in $L^2(\set{\Im\sigma=-\mu}\times Y)$ with respect
to
$d\sigma\otimes \mathfrak m_Y$. 
 
 The conormal symbol $\hat P_0$ of $P\in \diff^m_b(M)$ is the operator
valued polynomial defined by
 \begin{equation*}
 \hat P_0(\sigma)(u) = x^{-\im\sigma}P(x^{\im\sigma}\pi_Y^*u)|_Y,\quad
u\in C^{\infty}(Y;E_Y), \sigma\in \C.
 \end{equation*}
 It is easy to prove that $\widehat {P_0^\star}(\sigma) = (\hat
P_0(\overline \sigma))^*$. 

\section{Closed extensions}\label{SectGround}

Recall first the abstract situation (cf.~\cite{RiNa}), where
$A:\Dom_{\max}\subset H \to H$ is a densely defined closed operator in a
Hilbert space $H$. Thus $\Dom_{\max}$ is complete with the graph norm
$\|\cdot\|_{A}$ induced by the inner product $(u,v)+(Au,Av)$, and if
$\Dom\subset\Dom_{\max}$ is a subspace, then $A:\Dom\subset H \to H$ is a
closed operator if and only if $\Dom$ is closed with respect to
$\|\cdot\|_{A}$. Fix $\Dom_{\min}\subset \Dom_{\max}$, suppose $\Dom_{\min}$
is dense in $H$ and closed with respect to $\|\cdot\|_{A}$. Let
\begin{equation}\label{domainfamily}
\mathfrak D = \set{\Dom\subset \Dom_{\max} \st \Dom_{\min}\subset \Dom
\text{ and } \Dom \text { is closed w.r.t. }\|\cdot\|_{A}}
\end{equation}
Thus $\mathfrak D$ is in one to one correspondence with the set of
closed operators $A:\Dom\subset H \to H$ such that
$\Dom_{\min}\subset\Dom\subset \Dom_{\max}$.
For our purposes the following restatement is more appropriate.

\begin{proposition}\label{AbstractCharacterization}
The set $\mathfrak D$ is in one to one correspondence with the set of closed
subspaces of the quotient
\begin{equation}\label{E}
   \mathcal E =\Dom_{\max}/\Dom_{\min}
\end{equation}
\end{proposition}

In our concrete case $A\in x^{-\nu}\diff^m_b(M;E)$, $\nu>0$, is a $b$-elliptic
cone operator, considered initially as a densely defined unbounded operator
\begin{equation}\label{SmallDomain}
 A: C_c^\infty (M;E) \subset x^{-\nu/2}L^{2}_b(M;E)\to
    x^{-\nu/2}L^{2}_b(M;E).
\end{equation}
We take $\Dom_{\min}(A)$ as the closure of \eqref{SmallDomain} with
respect to the graph norm, and
\begin{equation*}
 \Dom_{\max}(A)=\set {u\in x^{-\nu/2}L^{2}_b(M;E) \st 
 A u\in x^{-\nu/2}L^{2}_b(M;E) },
\end{equation*}
which is also the domain of the Hilbert space adjoint of
\begin{equation*}
 A^\star:\Dom_{\min}(A^\star)\subset x^{-\nu/2}L^2_b(M;E) \to x^{-\nu/2} L^2_b(M;E).
\end{equation*}
 Thus $\Dom_{\max}$ is the largest subspace of $x^{-\nu/2}L^2_b(M;E)$ on which $A$
acts in the distributional sense and produces an element of $x^{-\nu/2}L^2_b(M;E)$; 
one can define $A$ on any subspace of
$\Dom_{\max}$ by restriction. These definitions have nothing to do with
ellipticity.  The following almost tautological lemma is based on the continuity of
$A:x^{\nu/2}H^m_b(M;E)\to x^{-\nu/2}L^2_b(M;E)$ and the fact that
$C_c^\infty(\Dot M;E)$ is dense in $x^{\nu/2}H^m_b(M;E)$.
\begin{lemma}\label{DMin} Suppose $A \in x^{-\nu/2}L^2_b(M;E)$, let $\Dom\subset
\Dom_{\max}(A)$ be such that $A:\Dom \subset x^{-\nu/2}L^{2}_b(M;E)\to
x^{-\nu/2}L^{2}_b(M;E)$ is closed. If $\Dom$ contains $x^{\nu/2}H^m_b(M;E)$ then
$\Dom$ contains $\Dom_{\min}(A)$. In particular, $x^{\nu/2}H^m_b(M;E)\subset
D_{\min}(A)$.
\end{lemma}

 Adding the $b$-ellipticity of $A$ as a hypothesis provides the following precise
characterization of $\Dom_{\min}(A)$:
\begin{proposition}\label{DomMin} If $A\in x^{-\nu/2}L^2_b(M;E)$ is
$b$-elliptic, then
\begin{enumerate}
\item \label{DomMinA} $\Dom_{\min}(A)=\Dom_{\max}(A)\cap 
 \big(\bigcap_{\eps>0}\,x^{\nu/2-\eps}H^{m}_b(M;E)\big)$
\item \label{DomMinB}
$\Dom_{\min}(A)=x^{\nu/2}H^m_b(M;E)$ if and only if \
$\spec_b(A)\cap \set{\Im\sigma=-\nu/2}=\varnothing$.
\end{enumerate}
\end{proposition}

The proof requires a number of ingredients, beginning with the
following fundamental result \cite{MM}, \cite{Sz91}:
\begin{theorem}\label{cone1}
Let $A\in x^{-\nu}\diff^m_b(M;E)$ be $b$-elliptic.
For every real $s$ and $\gamma$,
\[ A:x^{\gamma}H^s_b(M;E)\to x^{\gamma-\nu}H^{s-m}_b(M;E) \]
is Fredholm if and only if
$\,\spec_b(A)\cap \set{\Im\sigma=-\gamma}=\varnothing$.
In this case, one can find a bounded pseudodifferential parametrix
\[ Q:x^{\gamma-\nu}H^{s-m}_b(M;E)\to x^{\gamma}H^{s}_b(M;E) \]
such that $R=QA-1$ and $\tilde R=AQ-1$ are smoothing cone operators.
\end{theorem}

Note that if $A$ is $b$-elliptic we always can find an operator $Q$
such that
\begin{align*}
 QA-1 &:x^{\gamma}H^{s}_b(M;E)\to x^{\gamma}H^{\infty}_b(M;E)\; \text{ and}\\
 AQ-1 &:x^{\gamma-\nu}H^{s}_b(M;E)\to x^{\gamma-\nu}H^{\infty}_b(M;E)
\end{align*}
are bounded for every $s\in\R$, even if the boundary spectrum intersects the line
$\set{\Im\sigma=-\gamma}$. Also in this case $\ker A$ and $\ker A^\star$
are finite dimensional spaces. Moreover, for every $u\in x^{\gamma}H^{s}_b(M;E)$,
\begin{equation*}
  \|u\|_{x^{\gamma}H^{s}_b}
  \le \|(QA-1)u\|_{x^{\gamma}H^{s}_b}+\|QAu\|_{x^{\gamma}H^{s}_b},
\end{equation*}
which implies
\begin{equation}\label{normest1}
  \|u\|_{x^{\gamma}H^{s}_b}
  \le C_{s,\gamma}\Big(\|u\|_{x^{\gamma}H^{s}_b}
  +\|Au\|_{x^{\gamma-\nu}H^{s-m}_b}\Big)
\end{equation}
for some constant $C_{s,\gamma}>0$.

As noted above, there is a bounded operator
$Q:x^{\gamma-\nu}H^{s-m}_b\to x^{\gamma}H^{s}_b$
such that $R=QA-1:x^{\gamma}H^{s-m}_b\to x^{\gamma}H^{\infty}_b$
is bounded. Hence for $u\in x^{\gamma}H^{s}_b$

\begin{equation*}
  \|u\|_{x^{\gamma}H^{s}_b}
  \le \|Ru\|_{x^{\gamma}H^{s}_b}+\|QAu\|_{x^{\gamma}H^{s}_b},
\end{equation*}
which implies the estimate \eqref{normest1}.

\begin{lemma}\label{cext3}
There exists $\eps>0$ such that
\begin{equation*} \Dom_{\max}(A)\embed x^{-\nu/2+\eps}H^m_b(M;E).
\end{equation*}
\end{lemma}
\begin{proof}
The inclusion follows from \eqref{normest1} and the fact that, if
$u\in\Dom_{\max}(A)$ then $\hat u$ has no poles on
$\set{\Im\sigma=\nu/2}$. Choose $\eps>0$ smaller than the distance
between $\spec_b(A)\cap\set{\Im\sigma<\nu/2}$ and the line
$\set{\Im\sigma=\nu/2}$.
The continuity of the embedding
is a consequence of the closed graph theorem since $\Dom_{\max}(A)$ and
$x^{-\nu/2+\eps}H^m_b$ are both continuously embedded in $x^{-\nu/2}L^2_b$.
\end{proof}

Recall that $x^{-\nu/2+\eps}H^m_b(M;E)$ is compactly embedded in
$x^{-\nu/2}L^2_b(M;E)$. Hence if $A$ with domain
$\Dom$ is closed, then
\begin{equation}\label{CompEmb}
 (\Dom,\|\cdot\|_A)\embed x^{-\nu/2}L^2_b(M;E) \;\text{ compactly.}
\end{equation}
That this embedding is compact is a fundamental difference between the situation at
hand and $b$-elliptic totally characteristic operators and is due to the presence of
the factor $x^{-\nu}$ in $A$. 

\begin{lemma}\label{cext2}
Let $\gamma\in[-\nu/2,\nu/2]$ be such that
$\;\spec_b(A)\cap\set{\Im\sigma=-\gamma}=\varnothing$.
Then $A$ with domain $\Dom_{\max}(A)\cap x^{\gamma}H^{m}_b(M;E)$
is a closed operator on $x^{-\nu/2}L^{2}_b(M;E)$.
\end{lemma}
\begin{proof}
Let $\seq{u_n}\subset \Dom_{\max}(A)\cap x^{\gamma}H^{m}_b$ 
with $u_n\to u$ and $Au_n\to f$ in $x^{-\nu/2}L^{2}_b$. Further, let $Q$ be a
parametrix of $A$ as in Theorem~\ref{cone1}. In particular,
\begin{equation*}
 QA=1+R:x^{\gamma}H^m_b(M;E)\to x^{\gamma}H^m_b(M;E)\text{ is Fredholm}.
\end{equation*}
So, $Au_n\to f$ implies $(1+R)u_n\to Qf$ in $x^{\gamma}H^m_b$, and
$Qf=(1+R)\tilde u$ for some $\tilde u\in x^{\gamma}H^m_b$. Now, since
$\dim\ker(1+R)<\infty$, there is a closed subspace $H\subset x^{-\nu/2}L^2_b$
such that $x^{-\nu/2}L^2_b(M;E)=H\oplus\ker(1+R)$.
If $\pi_H$ denotes the orthogonal projection onto $H$, then
$\pi_H u_n\to \pi_H u$ in $x^{-\nu/2}L^2_b$ and, as above,
\begin{equation*}
 (1+R)\pi_H u_n\to (1+R)\pi_H\tilde u \;\text{ in }
x^{\gamma}H^m_b.
\end{equation*}
 Thus $\pi_H u_n\to\pi_H\tilde u$ in $x^{\gamma}H^m_b\embed
x^{-\nu/2}L^2_b$ which implies that  $\pi_H u=\pi_H\tilde u\in x^{\gamma}H^m_b$.

On the other hand, $(1-\pi_H)u_n\to(1-\pi_H)u\in
\ker(1+R)\subset x^{\gamma}H^m_b$. Therefore, $u\in x^{\gamma}H^m_b$ and $Au=f$. 
In other words, $A$ with the given domain is closed.
\end{proof}


Momentarily returning to the abstract situation of a densely defined closed
operator $A:\Dom_{\max}\subset H \to H$, let $\Dom_{\min}^\star$ be the domain
of its adjoint. Further, let $\Dom_{\max}^\star $ be the domain
of the adjoint of $A:\Dom_{\min}\subset H \to H$, and denote this adjoint by
$A^\star$. Then we have $\Dom_{\min}^\star\subset \Dom_{\max}^\star$.

Let $\mathfrak D^\star$ and $\mathcal E^\star$ be the analogues of
\eqref{domainfamily} and \eqref{E} for $A^\star$. Define
\begin{equation*}
 [\cdot,\cdot]_A:\Dom_{\max}\times \Dom_{\max}^\star\to\C
\end{equation*}
by
\begin{equation}\label{DmaxPairing}
 [u,v]_A = (Au,v) - (u,A^\star v).
\end{equation}
Then $[u,v]_A = 0$ if either $u\in \Dom_{\min}$ or $v \in
\Dom_{\min}^\star$, and $[\cdot,\cdot]_A$ induces a nondegenerate pairing
\begin{equation*}
 [\cdot,\cdot]_A^\flat :\mathcal E\times \mathcal E^\star\to \C.
\end{equation*}
Indeed, suppose $u \in \Dom_{\max}$ is such that $[u,v]_A = 0$ for all $v\in
\Dom_{\max}^\star$. Then $(Au,v) = (u,A^\star v)$ for all $v\in
\Dom_{\max}^\star$, which implies that $u$ belongs to the domain of the
adjoint of $A^\star:\Dom_{\max}^*\subset H \to H$, that is, $u\in
\Dom_{\min}$. Thus the class of $u$ in $\mathcal E$ is zero.
Likewise, if $v\in \Dom_{\max}^\star$ and $[u,v]_A=0$ for all $u \in
\Dom_{\max}$ then $v \in \Dom_{\min}^\star$.

Given $\Dom \in \mathfrak D$, let $\Dom^\perp\subset \Dom_{\max}^\star$
be the orthogonal of $\Dom$ with respect to $[\cdot,\cdot]_A$.

\begin{proposition}\label{orthoadjoint}
Let $\Dom\in \mathfrak D$. The adjoint $A^*$ of
$A:\Dom\subset H\to H$ is precisely the operator $A^\star$
restricted to $\Dom^{\perp}$.
Consequently, $A$ is selfadjoint if and only if $\Dom=\Dom^\perp$.
\end{proposition}
\begin{proof}
Let $A^*:\Dom^*\subset H \to H$ be the adjoint of $A|_{\Dom}$.
We have $(Au,v) = (u,A^\star v)$ for all $u \in \Dom$,
$v\in \Dom^\perp$. Thus $\Dom^\perp \subset \Dom^*$ and
$A^*v=A^\star v$ if $v \in \Dom^\perp$. On the other hand, since
$\Dom^* \subset \Dom_{\max}^\star$ (because $\Dom \supset \Dom_{\min}$),
it makes sense to compute $[u,v]_A$ for $u \in \Dom$ and $v \in \Dom^*$.
For such $u$, $v$ we then have $[u,v]_A =0$ since $A^*v=A^\star v$ for
every $v \in \Dom^*$. Thus $\Dom^* \subset \Dom^\perp$ which completes the proof that
$\Dom^* = \Dom^\perp$.
\end{proof}

\begin{proof}[{\it Proof of Proposition \ref{DomMin}}]
To prove part \ref{DomMinA} we will show inclusion in both directions.
If $\eps>0$ is such that
$\spec_b(A)\cap \set{\Im\sigma=-\nu/2+\eps}=\varnothing$, then
\begin{equation*}
x^{\nu/2}H^m_b(M;E)\subset \Dom_{\max}\cap x^{\nu/2-\eps}H^m_b(M;E),
\end{equation*}
so from Lemmas \ref{DMin} and \ref{cext2} we deduce 
$\Dom_{\min}(A)\subset \Dom_{\max}(A)\cap x^{\nu/2-\eps}H^{m}_b(M;E)$. Thus,
\begin{equation*}
\Dom_{\min}(A)\subset \Dom_{\max}(A)\cap
\bigcap_{\eps>0} x^{\nu/2-\eps}H^{m}_b(M;E).
\end{equation*}
To prove the reverse inclusion let 
$u\in \Dom_{\max}\cap \bigcap_{\eps>0} x^{\nu/2-\eps}H^{m}_b$ and
set $u_n=x^{1/n}u$ for $n\in\mathbb N$.
Then $\seq{u_n}$ is a sequence in $x^{\nu/2}H^{m}_b$ and as $n\to\infty$
\begin{equation*}
u_n\to u \text{ in } x^{\nu/2-\eps}H^{m}_b, \text{ and }
Au_n\to Au \text{ in } x^{-\nu/2-\eps}L^2_b
\end{equation*}
for every $\eps>0$. In particular, $x^{\eps}Au_n\to x^{\eps}Au$ in
$x^{-\nu/2}L^2_b$. Choose $\eps$ sufficiently small such that
$\Dom_{\max}(A^\star)\subset x^{-\nu/2+\eps}H^m_b$ (Lemma~\ref{cext3}).
Then for $v\in\Dom_{\max}(A^\star)$
\begin{equation*}
(Au_n,v)=(x^{\eps}Au_n,x^{-\eps}v)\to(x^{\eps}Au,x^{-\eps}v)=
(Au,v) \text{ as } n\to\infty.
\end{equation*}
On the other hand, $(u_n,A^\star v)\to (u,A^\star v)$ and
$(Au_n,v)=(u_n,A^\star v)$ since $u_n\in\Dom_{\min}(A)$. Hence
$(Au,v)=(u,A^\star v)$ for all $v\in\Dom_{\max}(A^\star)$, that is,
$[u,v]_A=0$ for all $v\in\Dom_{\max}(A^\star)$ which implies
$u\in\Dom_{\min}(A)$ since $\Dom_{\min}(A)=\Dom_{\max}(A^\star)^\perp$.

To prove part \ref{DomMinB}, suppose first that $\spec_b(A)\cap
\set{\Im\sigma=-\nu/2}=\varnothing$. Then Lemma~\ref{cext2} gives that $A$
with domain $\Dom_{\max}(A)\cap x^{\nu/2}H^m_b(M;E)$ is closed. Since
$x^{\nu/2}H^m_b(M;E) \subset \Dom_{\max}(A)$, Lemma~\ref{DMin} implies
$\Dom_{\min}(A)=x^{\nu/2}H^m_b(M;E)$. On the other hand, if
$\Dom_{\min}(A)=x^{\nu/2}H^m_b(M;E)$, then $A$ with domain $x^{\nu/2}H^m_b(M;E)$ 
is Fredholm, so by Theorem~\ref{cone1} $\,\spec_b(A)\cap
\set{\Im\sigma=-\nu/2}=\varnothing$
\end{proof}

We will now prove that $\Dom_{\max}/\Dom_{\min}$ is finite dimensional.
This is a consequence of the following proposition, which is interesting on
its own, see Lesch~\cite[Lemma 1.3.15, Prop. 1.3.16]{Le97}.

\begin{proposition}\label{FredholmProperty}
If $A\in x^{-\nu}\diff^m_b(M;E)$ is $b$-elliptic, every closed extension
\begin{equation*}
  A:\Dom\subset x^{-\nu/2}L^{2}_b(M;E)\to x^{-\nu/2}L^{2}_b(M;E)
\end{equation*}
is a Fredholm operator. Moreover, $\dim\Dom(A)/\Dom_{\min}(A)$ is finite, and
\begin{equation*}
 \Ind A|_{\Dom}=\Ind A|_{\Dom_{\min}}+\dim\Dom(A)/\Dom_{\min}(A).
\end{equation*}
In particular,
\begin{equation*}
 \mathcal{E}(A)=\Dom_{\max}(A)/\Dom_{\min}(A) \text{ is finite dimensional.}
\end{equation*}
\end{proposition}

Below we will give a proof different from that of Lesch. That the dimension of
$\Dom_{\max}(A)/\Dom_{\min}(A)$ is finite can also be proved by
observing that if
$A=x^{-\nu}P$ with
$P\in\diff^m_b(M;E)$ and $Au\in x^{-\nu/2}L^2_b(M;E)$, then $Pu \in
x^{\nu/2}L^2_b(M;E)$, so the Mellin transform of $Pu$ is holomorphic in
$\Im\sigma>-\nu/2$, from which it follows that $\hat u(\sigma)$ is
meromorphic in $\Im\sigma>-\nu/2$, and holomorphic in $\Im \sigma>\nu/2$
since $u\in x^{-\nu/2}L^2_b(M;E)$ (see \cite{K1}, also \cite{MM}).
On the other hand, if $u\in \Dom_{\min}$, then $\hat u(\sigma)$ is
holomorphic in $\Im\sigma>-\nu/2$. This type of argument leads to

\begin{corollary}\label{DMin=Dmax}
If $A\in x^{-\nu}\diff^m_b(M;E)$ is $b$-elliptic then
\[ \Dom_{\min}(A)=\Dom_{\max}(A) \text{ if and only if }\;
\spec_b(A)\cap \set{\Im\sigma\in(-\nu/2,\nu/2)}=\varnothing. \]
\end{corollary}

We will analyze $\mathcal{E}(A)$ more carefully in the next
sections. This will entail some repetition of work done by
Gohberg and Sigal \cite{GohSi}.

Our proof of Proposition~\ref{FredholmProperty} requires the following
classical result \cite{Wlo}:

\begin{lemma}\label{Apriori}
Let $X$, $Y$ and $Z$ be Banach spaces such that $X\embed Y$ is compact.
Further let $T\in\mathcal{L}(X,Z)$. Then the following conditions
are equivalent:
\begin{enumerate}
\item[1)] $\dim\ker T<\infty$ and $\;\rg T$ is closed,
\item[2)] there exists $C>0$ such that for every $u\in X$
\begin{equation*}
 \|u\|_{X}\le C(\|u\|_{Y}+\|Tu\|_{Z}) \quad(\text{a-priori estimate}).
\end{equation*}
\end{enumerate}
\end{lemma}

\begin{proof}[Proof of Proposition~\ref{FredholmProperty}]
Using \eqref{CompEmb} and Lemma~\ref{Apriori} with
$X=(\Dom,\|\cdot\|_A)$ and $Y=Z=x^{-\nu/2}L^{2}_b(M;E)$, we obtain
$\dim\ker A|_\Dom<\infty$ and $\,\rg A|_\Dom$ closed.
On the other hand, the adjoint $A^*$ of a cone operator $A$ is just a closed
extension of its $b$-elliptic formal adjoint $A^\star$,
cf. Proposition~\ref{orthoadjoint}. Applying again Lemma~\ref{Apriori} we get
$\dim\ker(A|_\Dom)^*<\infty$.

To verify the index formula consider the inclusion $\iota:\Dom_{\min}\to
\Dom$ which is clearly Fredholm (use the same argument but now with
$X=(\Dom_{\min},\|\cdot\|_A)$, $Y=x^{-\nu/2}L^{2}_b(M;E)$ and
$Z=(\Dom,\|\cdot\|_A)$). Then $\Ind\iota=-\dim\Dom/\Dom_{\min}$, hence
\[ \Ind A|_{\Dom_{\min}}= \Ind (A|_\Dom \circ \iota)
   =\Ind A|_\Dom + \Ind\iota=\Ind A|_\Dom - \dim\Dom/\Dom_{\min}. \]
\end{proof}

As a consequence of Propositions \ref{AbstractCharacterization} and
\ref{FredholmProperty}, for any closed extension of
\begin{equation*}
 A:C_c^\infty(\Dot M;E) \subset
 x^{-\nu/2}L^2_b(M;E)\to x^{-\nu/2}L^2_b(M;E)
\end{equation*}
with domain $\Dom(A)$ there is a finite dimensional space $\mathcal E'\subset
\Dom_{\max}(A)$ such that 
 \begin{equation*}
 \Dom(A)=\Dom_{\min}(A)\oplus \mathcal E' \;\text{ (algebraic direct sum).}
 \end{equation*}

 From Proposition \ref{FredholmProperty} we also get the following two corollaries
\begin{corollary}
Let $A:\Dom_1\to x^{-\nu/2}L^{2}_b(M;E)$ and 
$A:\Dom_2\to x^{-\nu/2}L^{2}_b(M;E)$ be closed extensions of $A$ such that
$\Dom_1\subset\Dom_2$, $\Ind A|_{\Dom_1}<0$ and $\Ind A|_{\Dom_2}>0$. 
Then there exists a domain $\Dom$ with $\Dom_1\subset\Dom\subset\Dom_2$
such that $\Ind A|_{\Dom}=0$.
\end{corollary}

The interest of this corollary lies in the fact that the vanishing of the
index is necessary for the existence of the resolvent and of
selfadjoint extensions.

\begin{corollary}
Let $A:\Dom_1\to x^{-\nu/2}L^{2}_b(M;E)$ and 
$A:\Dom_2\to x^{-\nu/2}L^{2}_b(M;E)$ be closed extensions of $A$. Then
\begin{equation*}
 \Ind A|_{\Dom_2}-\Ind A|_{\Dom_1}=
   \dim\Dom_{2}/\Dom_{\min}-\dim\Dom_{1}/\Dom_{\min}.
\end{equation*}
In particular, if $\Dom_1\subset\Dom_2$ then
\begin{equation*}
 \Ind A|_{\Dom_2}-\Ind A|_{\Dom_1}=\dim\Dom_{2}/\Dom_{1}.
\end{equation*}
\end{corollary}

This corollary implies the well-known relative index theorems for 
operators acting on the weighted Sobolev spaces $x^\gamma H^m_b(M;E)$, 
cf. \cite{MM}, \cite{Sz91}.

\section{Equality of domains}\label{EqOfDom}
In this section we give sufficient conditions for the domains of different operators
to be equal. Of these results, only the one concerning Friedrichs extensions will be
used later on. 

An operator $A \in x^{-\nu}\diff^m_b(M;E)$ is said to vanish
on $\partial M$ to order $k$ ($k\in \mathbb N$) if for any $u\in
C^\infty(M;E)$, $x^\nu Au$ vanishes to order $k$ on $\partial M$. Let
\begin{equation*}
 \ceil{s}=\min\set{k\in\mathbb N\st s\leq k}.
\end{equation*}

\begin{proposition}\label{EdDom}
Let $A_0$, $A_1 \in x^{-\nu}\diff^m_b(M;E)$ be $b$-elliptic.
\begin{enumerate}
\item\label{EdDom1}
If $A_0-A_1$ vanishes on $\partial M$, then
$\Dom_{\min}(A_0)=\Dom_{\min}(A_1)$.
\item\label{EdDom2}
If $A_0-A_1$ vanishes to order $\ell\leq{\nu-1}$ on $\partial M$, 
$\ell\in\mathbb{N}$, then
\begin{equation*}
 \Dom_{\max}(A_0)\cap x^{\frac{\nu}{2}-\ell-1}H^m_b(M;E)=
 \Dom_{\max}(A_1)\cap x^{\frac{\nu}{2}-\ell-1}H^m_b(M;E).
\end{equation*}
\item \label{EdDom3}
If $A_0-A_1$ vanishes to order $\ceil{\nu-1}$
on $\partial M$, then $\Dom_{\max}(A_0)=\Dom_{\max}(A_1)$.
\item\label{EdDom4} If $A_0$ and $A_1$ are symmetric and bounded from below and
$A_0-A_1$ vanishes to order $\ceil{\nu-1}$ on $\partial M$, then the domains of their
Friedrichs extensions coincide, that is, $\Dom_{F}(A_0)=\Dom_{F}(A_1)$.
\end{enumerate}
\end{proposition}
\begin{proof}
First of all, observe that in all cases it is enough to prove only one inclusion;
the equality of the sets follows then by exchanging the roles of
$A_0$ and $A_1$.

To prove part~\ref{EdDom1}, write
\begin{equation*}
A_1 = A_0 + (A_1-A_0)= A_0 + x^{-\nu}Px
\end{equation*}
with $P\in\diff^m_b(M;E)$ and suppose $u\in \Dom_{\min}(A_0)$. There is then a
sequence $\set{u_n}_{n\in
\mathbb N}\subset C_c^\infty(\Dot M;E)$ such that $u_n\to u$ and
$A_0 u_n \to A_0u$, in $x^{-\nu/2}L^2_b$. Consequently, $x u_n \to x u$ in
$x^{\nu/2} H^m_b$ and $x^{-\nu}Pxu_n \to x^{-\nu}Pxu$ in
$x^{-\nu/2}L^2_b$. Thus $A_1u_n \to A_0 u + x^{-\nu} Pxu$ which implies
$\Dom_{\min}(A_0)\subset \Dom_{\min}(A_1)$.

Now, let $\ell\in\mathbb{N}$, $\ell\leq\nu-1$, and let
$u\in \Dom_{\max}(A_0)\cap x^{\nu/2-\ell-1}H^m_b$.
This means $u\in x^{\nu/2-\ell-1}H^m_b$ and $A_0 u\in x^{-\nu/2}L^2_b$.
To prove that $u\in \Dom_{\max}(A_1)\cap x^{\nu/2-\ell-1}H^m_b$
we only need to show that $A_1 u$ belongs to $x^{-\nu/2}L^2_b$.
Let $P\in\diff^m_b(M;E)$ be such that
$A_1=A_0 + x^{-\nu}Px^{\ell+1}$.
Since $u\in x^{\nu/2-\ell-1}H^m_b$ then $x^{-\nu}Px^{\ell+1}u\in
x^{-\nu/2}L^2_b$. Hence $A_1 u\in x^{-\nu/2}L^2_b$ which proves the second
statement.

To prove the third statement, let $u\in\Dom_{\max}(A_0)$, $\ell=\ceil{\nu-1}$, and $P$
as above. Then $u\in x^{-\nu/2}H^m_b$ and $x^{-\nu}Px^{\ell+1}u\in
x^{-\nu/2}L^2_b$. Thus $A_1u\in x^{-\nu/2}L^2_b$ and
$\Dom_{\max}(A_0)\subset\Dom_{\max}(A_1)$.

To prove part \ref{EdDom4} we first prove the rather useful and well known abstract
characterization of the domain of the Friedrichs extension of a symmetric semibounded
operator given in Lemma~\ref{DFApproximation} below. Suppose
$A:\Dom_{\min}\subset H
\to H$ is a densely defined closed operator which is symmetric and bounded from
below. Let $A^\star:\Dom_{\max}\subset H \to H$ be its adjoint and 
let $A_F:\Dom_F\subset H\to H$ be the Friedrichs extension of $A$. Define
\begin{equation}\label{FthyQprime}
 (u,v)_{A^\star} = c(u, v) + (A^\star u, v)\quad
 \text{for }u,\, v\in\Dom_{\max},
\end{equation}
where $c=1-c_0$ and $c_0\leq 0$ is a lower bound of $A$.

\begin{lemma} \label{DFApproximation}
$u \in \Dom_{\max}$ belongs to $\Dom_F$ if and only if there exists a
sequence $\set {u_n}_{n\in \mathbb N}$ in $\Dom_{\min}$ such that
\begin{equation*}
 (u-u_n,u-u_n)_{A^\star} \to 0 \text{ as } n\to\infty.
\end{equation*}
\end{lemma}
\begin{proof}
Because of the fact that $\Dom_{\min}$ is dense in $\Dom_F$ with respect to
the norm $\|\cdot\|_{A^\star}$ induced by \eqref{FthyQprime}, every
$u\in\Dom_F$ can be approximated as claimed.

Let now $u\in\Dom_{\max}$ and let
$\set{u_n}_{n\in\mathbb N}\subset\Dom_{\min}$ be such that
$(u-u_n,u-u_n)_{A^\star}\to 0$, so $u_n\to u$ in $H$.
Let $\mathcal K\subset H$ be the domain of the positive square root $R$ of $A_F+cI$.
Recall that $\Dom_F=\Dom_{\max}\cap \mathcal K$. Hence $u\in\Dom_{\max}$
belongs to $\Dom_F$ if $u\in\mathcal K$, and the identity
\begin{equation*}
 \|u_n-u_\ell\|_{A^\star} = \|R(u_n-u_\ell)\|
\end{equation*}
implies that $\set{Ru_n}_{n\in \mathbb N}$ also converges in $H$.
Thus $u\in \mathcal K$ since $R$ is closed.
\end{proof}

We now prove part~\ref{EdDom4} of Proposition~\ref{EdDom}. Suppose that $A_0$ and
$A_1$ satisfy the hypotheses there. Then by parts \ref{EdDom1} and \ref{EdDom3},
$\Dom_{\min}=\Dom_{\min}(A_0)=\Dom_{\min}(A_1)$ and
$\Dom_{\max}=\Dom_{\max}(A_0)=\Dom_{\max}(A_1)$, and from
the fact that $A_0+A_1-2A_0$ vanishes to order $\ceil{\nu-1}$ we also get that
$\Dom_{\min}(A_0+A_1)=\Dom_{\min}$ and $\Dom_{\max}(A_0+A_1)=\Dom_{\max}$. Since
$A_0$ and $A_1$ are symmetric and bounded from below, so is $A_0+A_1$. We will
show that these three operators share the same Friedrichs domain by showing
that 
 \begin{equation}\label{EdDom4.1}
 \Dom_F(A_0+A_1) \subset \Dom_F(A_0)\cap \Dom_F(A_1).
 \end{equation}
Suppose
this has been shown. Since $[u,v]_{A_0} =
[u,v]_{A_1}=\frac{1}{2}[u,v]_{A_0+A_1}$, Proposition~\ref{orthoadjoint} implies that 
$A_0$ is selfadjoint with either of the domains $\Dom_F(A_0)$ or
$\Dom_F(A_0+A_1)$, and from the inclusion of the latter in the former one deduces the
equality of these spaces, hence, that $\Dom_F(A_0)=\Dom_F(A_1)$. To prove
\eqref{EdDom4.1}, suppose
$(A_i u,u)_{x^{-\nu/2}L^2_b} \geq c_i(u,u)_{x^{-\nu/2}L^2_b}$, $i=0$, $1$, on
$\Dom_{\min}$, let
$c=1-c_0-c_1$. If $u \in \Dom_F(A_0+A_1)$, then by Lemma~\ref{DFApproximation}
there is a sequence $\set{u_n}_{n\in\mathbb N}\subset \Dom_{\min}$ such that 
\begin{multline*}
 (A_0(u-u_n), u-u_n)_{x^{-\nu/2}L^2_b} + (A_1(u-u_n), u-u_n)_{x^{-\nu/2}L^2_b} \\ +
c(u-u_n,u-u_n)_{x^{-\nu/2}L^2_b} \to 0
\text{ as } n\to\infty. 
\end{multline*}
But then also
\begin{align*}
 (1-c_i)(u-u_n,u-u_n)_{x^{-\nu/2}L^2_b} + (A_i(u-u_n), u-u_n)_{x^{-\nu/2}L^2_b}
\end{align*}
as $n\to\infty$, $i=0$, $1$, so $u\in \Dom_F(A_0)\cap \Dom_F(A_1)$.
\end{proof}

\section{Spaces of Meromorphic Solutions}\label{MeromorphicSolutions}

If $K$ is a finite dimensional complex vector space, we let $\mathfrak
M_{\sigma_0}(K)$ be the space of germs of
$K$-valued meromorphic functions with pole at $\sigma_0$ and
$\mathfrak{Hol}_{\sigma_0}(K)$ be the subspace of holomorphic germs.
These are naturally modules over the ring
$\mathfrak{Hol}_{\sigma_0}(\C)$.  

Let $R^\perp$ be another finite dimensional complex vector space. If
$\mathcal{P}(\sigma): K\to R^\perp$ is a linear map depending holomorphically on
$\sigma$ in a neighborhood of $\sigma_0$, then $\mathcal{P}$ defines a map
$\mathfrak M_{\sigma_0}(K) \to \mathfrak{M}_{\sigma_0}(R^\perp)$, which we also
denote by $\mathcal{P}$.

\begin{lemma}\label{mero1}
Suppose that $\mathcal{P}(\sigma)$ is defined near
$\sigma=\sigma_0$, is invertible for $\sigma\ne \sigma_0$ but
$\mathcal{P}(\sigma_0)=0$. Then there are $\psi_1,\dots,\psi_d\in \mathfrak
M_{\sigma_0}(K)$ be such that
$\beta_j(\sigma)=\mathcal{P}(\sigma)(\psi_j(\sigma))$ is
holomorphic and $\beta_1(\sigma_0),\dots,\beta_d(\sigma_0)$
is a basis of $R^\perp$. For any such $\psi_j$, if $u\in\mathfrak M_{\sigma_0}(K)$
and
$\mathcal{P} u
\in \mathfrak{Hol}_{\sigma_0}(R^\perp)$, then there are $f_j\in \mathfrak
{Hol}_{\sigma_0}(\C)$ such that $u = \sum_{j=1}^d f_j\psi_j$.
\end{lemma}
\begin{proof} Let $\set{b_j}_{j=1}^d$ be a basis of $R^\perp$ and define $\psi_j=
\mathcal{P}^{-1}(b_j)$. Then the $\psi_j$ are meromorphic with pole at
$\sigma_0$, $\mathcal{P}\psi_j=\beta_j=b_j$ is holomorphic, and the
$\beta_j(\sigma_0)$ form a basis of $R^\perp$.

Let now $\psi_1,\dots,\psi_d\in \mathfrak
M_{\sigma_0}(K)$ be such that
$\beta_j(\sigma)=\mathcal{P}(\sigma)(\psi_j(\sigma))$ is
holomorphic and $\beta_1(\sigma_0),\dots,\beta_d(\sigma_0)$
is a basis of $R^\perp$. If $f\in \mathfrak{Hol}_{\sigma_0}(R^\perp)$ then
$f=\sum_j f_j\beta_j$ for some $f_j\in\mathfrak{Hol}_{\sigma_0}(\C)$, because the
$\beta_j(\sigma_0)$ form a basis, and each $\beta_j(\sigma_0)$ can be
written as a linear combination (over $\mathfrak{Hol}_{\sigma_0}(\C)$) of
$\beta_1,\dots,\beta_d$. If $\mathcal{P} u=f$ then
$\mathcal{P}(u-\sum_{j=0}^d f_j\psi_j) = 0$, so $u=\sum_{j=0}^d f_j\psi_j$
for $\sigma \ne \sigma_0$, which is the equality of meromorphic functions.
\end{proof}

The lemma asserts that $\mathcal{P}^{-1}(\mathfrak{Hol}_{\sigma_0}(R^\perp))$
is finitely generated as a submodule of $\mathfrak M_{\sigma_0}(K)$ over
$\mathfrak{Hol}_{\sigma_0}(\C)$. We will be interested in
$\hat{\mathcal E}_{\sigma_0}= \mathcal{P}^{-1}(\mathfrak{Hol}_{\sigma_0}(R^\perp))/
\mathfrak{Hol}_{\sigma_0}(K)$ as a vector space over $\C$.
The following fundamental lemma paves the way to describing a basis
of $\hat{\mathcal E}_{\sigma_0}$.

\begin{lemma}\label{mero2}
Suppose that $\mathcal{P}(\sigma)$ is defined near $\sigma=\sigma_0$,
is invertible for $\sigma\ne \sigma_0$ but $\mathcal{P}(\sigma_0)=0$.
There are $\psi_1,\dots,\psi_d\in \mathfrak M_{\sigma_0}(K)$ such that each
$\beta_j(\sigma)=\mathcal{P}(\sigma)(\psi_j(\sigma))$ is holomorphic,
$\beta_1(\sigma_0),\dots,\beta_d(\sigma_0)$ form a basis of $R^\perp$, and if
 \begin{equation}\label{mero3}
 \psi_j=\sum_{\ell=0}^{\mu_j-1} \frac 1
 {(\sigma-\sigma_0)^{\mu_j-\ell}}\, \psi_{j\ell} + h_j
 \end{equation}
with holomorphic $h_j$ then the $\psi_{j0}$ are linearly independent.
\end{lemma}
\begin{proof}
Without loss of generality assume $\sigma_0=0$.
Pick a basis  $\set{b_j}_{j=1}^d$ of $R^\perp$ and define $\psi_j=
\mathcal{P}^{-1}(b_j)$. Then the $\psi_j$ are meromorphic with pole at
$0$, $\mathcal{P}\psi_j=\beta_j=b_j$ is holomorphic, and the
$\beta_j(0)$ are independent. Each $\psi_j$ can be written as
\eqref{mero3}. Order them so that $\set{\mu_j}_{j=1}^d$
is nonincreasing and let
\begin{equation}
 \begin{gathered}\label{mero4}
\tilde \mu_1=\max\set{\mu_j\st j=1,\dots,d}\\
\tilde \mu_i=\max \set{\mu_j\st \mu_j<\tilde \mu_{i-1}}\quad i=2,\dots,L\\
 s_i=\max \set{j\st \mu_j=\tilde \mu_i},
 \end{gathered}
 \end{equation}
 that is,
 \begin{equation*}
 \tilde \mu_1=\mu_1=\dots=\mu_{s_1}>
 \tilde \mu_2=\mu_{s_1+1}=\dots=\mu_{s_2} >\dots >\tilde \mu_L =
 \mu_{s_{L-1}+1}=\dots=\mu_d
 \end{equation*}
 If the vectors $\psi_{j0}$,
$j=1,\dots,s_1$, are not linearly independent, then order the $\psi_j$ with
$j\leq s_1$ so that $\psi_{10},\dots,\psi_{s'_1,0}$ is a maximal set of
linearly independent vectors among $\set{\psi_{j0}\st 1\leq j\leq s_1}$, write
 \begin{equation*}
 \psi_{k0}=\sum_{j=1}^{s'_1} a_{kj}\psi_{j0} \
 \text{ for } k=s'_1+1,\dots,s_1,
 \end{equation*}
 and replace $\psi_k$ by $\psi_k - \sum_{j=1}^{s'_1} a_{kj}\psi_j$
for $k=s'_1+1,\dots, s_1$. Now
$\mathcal{P}(\psi_k)=\beta_k-\sum_{j=1}^{s'_1} a_{kj}\beta_j$ for these
indices, so it is still true that the $\mathcal{P}(\psi_j)(0)$ form a
basis. With $\mu_j$ denoting the order of the pole of the new $\psi_j$, and
again assuming the orders form a nonincreasing sequence, let $\tilde \mu_i$
and $s_i$ be defined as above. Suppose that already $\psi_{j0}$,
$j=1,\dots, s_i$ is an independent set. If $\psi_{s_i+1,0}$ depends linearly
on $\psi_{1,0},\dots,\psi_{s_i,0}$ then put $s'_{i+1}=s_i$. Otherwise,
reorder $\psi_{s_i+1,0},\dots, \psi_{s_{i+1},0}$ so that
$\psi_{s_i+1,0},\dots,\psi_{s'_{i+1},0}$ together with $\psi_{j0}$,
$j=1,\dots,s_i$ are a maximally independent set in
$\set{\psi_{j0}\st 1\leq j\leq s_{i+1}}$. If $s'_{i+1}<s_{i+1}$ write
 \begin{equation*}
 \psi_{k0}=\sum_{j=1}^{s'_{i+1}} \alpha_{kj}\psi_{j0},
 \quad k=s'_{i+1}+1,\dots, s_{i+1}.
 \end{equation*}
 Replacing $\psi_k$ by $\psi_k-\sum_{j=1}^{s'_{i+1}}
\alpha_{kj}\,\sigma^{\mu_j-\tilde \mu_{i+1}}\psi_j$
$(s'_{i+1}+1\leq k\leq s_{i+1})$, reordering by decreasing order of the pole
(which reorders only $\psi_j$, $j> s'_{i+1}$), now have that the leading
coefficients of the $\psi_j$, $j\leq s_{i+1}$, are independent.
\end{proof}

\begin{lemma}\label{mero8}
With the setup of Lemma~\ref{mero2}, let $\psi_1,\dots,\psi_d\in
\mathfrak M_{\sigma_0}(K)$ be as stated there, and let $\mu_j$ be the order
of the pole of $\psi_j$. Let
\begin{equation*}
 \hat {\mathcal E}_{\sigma_0}
 = \mathcal{P}^{-1}(\mathfrak{Hol}_{\sigma_0}(R^\perp))/
 \mathfrak{Hol}_{\sigma_0}(K),
\end{equation*}
regarded as a vector space over $\C$. Then the images in $\hat {\mathcal
E}_{\sigma_0}$ of the elements 
\begin{equation*}
(\sigma-\sigma_0)^\ell\psi_j,\quad j=1,\dotsc,d,\ \ell=0,\dotsc,\mu_j-1
\end{equation*}
form a basis of this space.
\end{lemma}
\begin{proof}
As before assume $\sigma_0=0$. Because of Lemma~\ref{mero1} the images of the
$\sigma^\ell\psi_j$ span, and we only need to prove linear independence.
Suppose $\sum_{j=1}^d\sum_{k=0}^{\mu_j-1}u_{jk}\sigma^k\psi_j$
is holomorphic. Modulo holomorphic functions,
 \begin{equation*}
 \psi_j=\sum_{\ell=0}^{\mu_j-1}\frac{1}{\sigma^{\mu_j-\ell}}\,\psi_{j\ell}
 \end{equation*}
(where the $\psi_{j0}$ are independent) so
 \begin{equation*}
u(\sigma)=\sum_{j=1}^d\sum_{k=0}^{\mu_j-1}
\sum_{\ell=0}^{\mu_j-k} \frac{u_{jk}}{\sigma^{\mu_j-\ell-k}}\,\psi_{j\ell}
 \end{equation*}
is holomorphic. Thus $\sigma^\nu u(\sigma)$ vanishes at $0$ for $\nu>0$.
Let the $\tilde \mu_i$ be as in \eqref{mero4}. We have that
$\sigma^{\tilde \mu_1} \psi_j(\sigma)$ vanishes at $0$ for $j>s_1$,
and so does $\sigma^{\tilde\mu_1}\sigma^k\psi_{j}(\sigma)$ for $k>0$. Hence
\begin{equation*}
0=\Big(\sigma^{\tilde\mu_1}u(\sigma)\Big)\Big|_{\sigma=0}
 =\sum_{j=1}^{s_1}u_{j0}
  \Big(\sigma^{\tilde\mu_1}\psi_{j}(\sigma)\Big)\Big|_{\sigma=0}
 =\sum_{j=1}^{s_1}u_{j0}\psi_{j0},
\end{equation*}
and so $a_{0j}=0$ for $j=1,\dots,s_1$, since the $\psi_{j0}$ are independent.
If $\tilde \mu_2<\tilde \mu_1-1$ then by the same argument one concludes that
$a_{1j}=0$ for $j=1,\dots, s_1$, and if $\tilde \mu_2=\tilde \mu_1-n$,
($n \geq 1$) then the conclusion is that $a_{kj}=0$ for $j=1,\dots,s_1$ and
$k=0,\dots,n-1$. Having proved this, we conclude
 \begin{equation*}
 u(\sigma)= \sum_{j=1}^{s_1}\sum_{k=n}^{\tilde \mu_1-1}
\sum_{\ell=0}^{\tilde \mu_1-k} \frac {u_{jk}} {\sigma^{\tilde
\mu_1-\ell-k}}\,\psi_{j\ell} +
 \sum_{j=s_1+1}^d\sum_{k=0}^{\mu_j-1}
\sum_{\ell=0}^{\mu_j-k} \frac {u_{jk}}
{\sigma^{\mu_j-\ell-k}}\,\psi_{j\ell}.
 \end{equation*}
Now, since $\sigma^{\tilde \mu_2} u(\sigma)$ also vanishes at $0$ then
 \begin{equation*}
 \sum_{j=1}^{s_1} u_{jn}\psi_{j0}+\sum_{j=s_1+1}^{s_2} u_{j0}\psi_{j0}=0,
 \end{equation*}
 therefore $u_{jn}=0$ for $j=1,\dots,s_1$, and $u_{j0}=0$ for
$j=s_1+1,\dots,s_2$. Continuing in this manner, one obtains $u_{jk}=0$
for all $j$, $k$.
 \end{proof}

\begin{lemma}
With the setup of Lemma \ref{mero2}, let $\psi_1,\dots,\psi_d$ be as stated
there, let $\mu_j$ be the order of the pole of $\psi_j$. Suppose the $\psi_j$
ordered so that $\set{\mu_j}_{j=1}^d$ is nonincreasing. With the notation in
formulas\eqref{mero3} and \eqref{mero4} let
\begin{equation*}
 K_{\tilde \mu_\ell}= \LinSpan_\C\set{\psi_{j0}\st \mu_j\geq\tilde \mu_\ell}.
\end{equation*}
The spaces $K_{\tilde \mu_\ell}$ are independent of the choice of $\psi_j$.
\end{lemma}
\begin{proof}
Let $\mathfrak K_\mu= \set{\psi\in \mathcal P^{-1}(\mathfrak
{Hol}_{\sigma_0}(R^\perp))\st \mathrm{ord}(\psi)\leq \mu}$.
Thus if $\psi\in \mathfrak K_\mu$ then $(\sigma-\sigma_0)^\mu\psi$
is regular; let $m_\mu:\mathfrak K_\mu\to K$ be defined by setting
\begin{equation*}
 m_\mu(\psi)=(\sigma-\sigma_0)^\mu\psi(\sigma)|_{\sigma=\sigma_0}.
\end{equation*}
We will show that $K_{\tilde\mu_\ell}=m_{\tilde\mu_\ell}(\mathfrak
K_{\tilde\mu_\ell})$. To see this, set
\begin{equation*}
 \mathfrak K_{\tilde \mu_\ell}^\psi = \LinSpan_{
 \mathfrak{Hol}_{\sigma_0}(\C)}\set{\psi_j\st \mathrm{ord}(\psi_j)
 =\tilde \mu_\ell}.
\end{equation*}
and note that if $\mu\geq \tilde\mu_1$ then
\begin{equation*}
 \mathfrak K_\mu = \mathfrak K_{\tilde \mu_1}^{\psi} +
 \mathfrak K_{\tilde \mu_2}^{\psi} + \dots +
 \mathfrak K_{\tilde \mu_L}^{\psi}
\end{equation*}
and if $\tilde \mu_{\ell-1}\geq \mu\geq \tilde \mu_{\ell}$ then
\begin{equation*}
 \mathfrak K_\mu
 = (\sigma-\sigma_0)^{\tilde \mu_1-\mu}\mathfrak K_{\tilde\mu_1}^{\psi}
 + \dots + (\sigma-\sigma_0)^{\tilde \mu_{\ell-1}-\mu}
   \mathfrak K_{\tilde\mu_{\ell-1}}^{\psi}
 + \mathfrak K_{\tilde\mu_\ell}^\psi +\dots
 + \mathfrak K_{\tilde\mu_L}^\psi.
\end{equation*}
This is proved using Lemma \ref {mero1}. Thus if
$\,\tilde\mu_{\ell-1}\geq \mu > \tilde \mu_{\ell}$ then
\begin{equation*}
 m_\mu(\mathfrak K_\mu)=\LinSpan_\C\set{\psi_{j0}\st \mu_j
 > \tilde\mu_{\ell}}
\end{equation*}
 and if $\mu=\tilde \mu_\ell$ then
\begin{equation*}
 m_\mu(\mathfrak K_\mu)= \LinSpan_\C\set{\psi_{j0}\st \mu_j
 \geq \tilde \mu_{\ell}}.
\end{equation*}
\end{proof}

Note that
$K_{\tilde \mu_1}\subset \dots \subset K_{\tilde \mu_L}=K$,
$\dim K_{\tilde \mu_\ell} = s_\ell$ and $m_{\tilde \mu_\ell}:\mathfrak
K_{\tilde \mu_\ell}\to K_{\tilde \mu_\ell}$ is surjective. As in see Gohberg and
Sigal \cite{GohSi}, the numbers $\mu_j$ will be called the partial multiplicities of
$\mathcal P$ (at $\sigma_0$) .

\begin{lemma} \label{mero11}
Let $\set{ \psi_i^\star}_{i=1}^d$, $\set{\psi_j}_{j=1}^d \subset
\mathfrak M_{\sigma_0}(K)$ be as in Lemma \ref{mero2}, both sequences ordered so that
the sequences $\set{\mu_i^\star}$, $\set{\mu_j}$ of the orders of the poles is
nonicreasing. Then $\mu_i^\star=\mu_i$ for all $i$ and 
 \begin{equation*}
 \psi_i^\star=\sum_{j=1}^d f_{ij}\psi_j
 \end{equation*}
 where the  $f_{ij}$ are holomorphic, form a nonsingular matrix and
$f_{ij}=(\sigma-\sigma_0)^{\mu_i-\mu_j}\tilde f_{ij}$ for some holomorphic
$\tilde f_{ij}$ if $\mu_i>\mu_j$. Conversely, given holomorphic functions $f_{ij}$
forming a nonsingular matrix and with $f_{ij}/(\sigma-\sigma_0)^{\mu_i-\mu_j}$
holomorphic when $\mu_i>\mu_j$, then the $\psi_j'$ defined by the formula above
satisfy the conclusion of Lemma \ref{mero2}.
\end{lemma}

We leave the proof of this to the reader. It uses the previous lemma and its proof. 

If $K$ is a hermitian vector space, let $\phi_1,\dots,\phi_d$ be an orthonormal basis
of $K$ such that for each $\ell=1,\dots,L$,
 \begin{equation*}
 \phi_{1},\dots,\phi_{s_\ell} \in K_{\tilde \mu_\ell}
 \end{equation*}
 Then for $j=s_{\ell-1}+1,\dots,s_\ell$ we can pick $\psi_j\in \mathfrak
K_{\tilde \mu_\ell}$ such that $m_{\tilde \mu_\ell}(\psi_j) =
\phi_j$, that is, if $K$ is hermitian then the $\psi_j$ can be chosen
to have orthogonal leading coefficients.

\begin{proposition}\label{mero6}
Let $\mathcal{P}(\sigma):K\to R^\perp$ be defined and holomorphic near
$\sigma=\sigma_0$, invertible for $\sigma\ne \sigma_0$
but $\mathcal{P}(\sigma_0)=0$. Then
\begin{enumerate}
\item\label{mero6.1}
there are $\psi_1,\dots,\psi_d\in \mathfrak M_{\sigma_0}(K)$ such that each
$\beta_j=\mathcal{P}\psi_j\in\mathfrak{Hol}_{\sigma_0}(R^\perp)$,
$\beta_1(\sigma_0),\dots,\beta_d(\sigma_0)$ form a basis of $\C^d$, and if
\begin{equation*}
 \psi_j=\sum_{\ell=0}^{\mu_j-1} \frac 1
 {(\sigma-\sigma_0)^{\mu_j-\ell}}\, \psi_{j\ell} + h_j
\end{equation*}
with holomorphic $h_j$ then the $\psi_{j0}$ are linearly independent,
\item\label{mero6.2}
if $K$ is a hermitian vector space, then the $\psi_j$ can even be chosen
such that the $\psi_{j0}$ form an orthonormal basis of $K$ and for $\ell>0$,
$\psi_{j\ell}$ is orthogonal to $\psi_{k0}$ whenever $\mu_k\geq \mu_j-\ell$.
\end{enumerate}
\end{proposition}

\begin{proof}
Because of Lemma~\ref{mero2} there are $\psi_1,\dots,\psi_d$ satisfying
the first statement.

Let now $K$ be a hermitian vector space. We may assume that already the
leading coefficients form an orthonormal basis of $K$. If a coefficient
$\psi_{j\ell}$ with $\ell > 0$ is not already orthogonal to those $\psi_{k0}$
such that $\mu_k \geq \mu_j-\ell$, then write
 \begin{equation*}
 \psi_{j\ell}= \psi_{j\ell}^0 +
 \sum_{\set{k\st \mu_k\geq \mu_j-\ell}} a_k \psi_{k0}
 \end{equation*}
 where $\psi_{j\ell}^0$ is orthogonal to the $\psi_{k0}$ such that
 $\mu_k \geq \mu_j-\ell$. Then
 \begin{equation*}
 \chi(\sigma)=\sum_{\set{k\st \mu_k\geq \mu_j-\ell}} a_k
 (\sigma-\sigma_0)^{\mu_k-\mu_j+\ell}\psi_{k}(\sigma) \in
 \mathfrak K_{\mu_j-\ell}\subset \mathfrak K_{\mu_j}
 \end{equation*}
($\mathfrak K_{\mu}$ being defined using $\sigma_0$) and
$(\sigma-\sigma_0)^{\mu_j}\chi(\sigma)|_{\sigma=\sigma_0}=0$.
So $\psi_j - \chi$ has the same leading term as
$\psi_j$ but now the coefficient of $(\sigma-\sigma_0)^{-\mu_j+\ell}$ is
$\psi_{j\ell}^0$ which is orthogonal to $\psi_{k0}$ for $k$ such that
$\mu_k \geq \mu_j-\ell$. We may then replace $\psi_j$ by $\psi_j-\chi$.
The proof is completed by `reverse' induction on $n=\mu_j-\ell$ beginning
with $n=\tilde \mu_1 -1$, the above being both the first and general steps.
\end{proof}


Let now $Y$ be a compact manifold and $E$
a complex vector bundle over $Y$. We fix a hermitian metric on $E$ and
riemannian metric on $E$ with respect to which we define the Sobolev
spaces $H^s(Y;E)$. Let
$\mathcal{P}(\sigma):H^m(Y;E)\to L^2(Y;E)$ be a holomorphic family of
elliptic operators of order $m$ defined for $\sigma$ near $\sigma_0$ in $\C$.
Suppose $\mathcal{P}(\sigma)$ is invertible for $\sigma \ne \sigma_0$ but
$\mathcal{P}(\sigma_0)$ is not invertible. Let
$K=\ker \mathcal{P}(\sigma_0)$, $R=\rg
\mathcal{P}(\sigma_0)$. Then $K$ and $R^\perp$ are finite dimensional
of the same dimension, say $d$, and consist of smooth sections of $E$.
Regard $\mathcal{P}(\sigma)$ as an operator
\begin{equation*}
 \left[\begin{matrix}
 \mathcal{P}_{11}(\sigma) & \mathcal{P}_{12}(\sigma)\\
 \mathcal{P}_{21}(\sigma) & \mathcal{P}_{22}(\sigma)
 \end{matrix}\right]:
 \begin{matrix} K \\ \oplus \\ K^\perp \end{matrix} \to
 \begin{matrix} R^\perp \\ \oplus \\ R \end{matrix}
\end{equation*}
in the usual way. All the $\mathcal{P}_{ij}$ are holomorphic, and the
operator $\mathcal{P}_{22}(\sigma)$ is invertible for $\sigma$ close
to $\sigma_0$. Thus $\mathcal{P}_{11} - \mathcal{P}_{12}\mathcal{P}_{22}^{-1}
\mathcal{P}_{21}: K\to R^\perp$ depends holomorphically on
$\sigma\in U$ and is invertible for $\sigma \ne \sigma_0$. We can then find
$\tilde\psi_1,\dots,\tilde \psi_d\in \mathfrak M_{\sigma_0}(K)$ such that
$(\mathcal{P}_{11} - \mathcal{P}_{12}\mathcal{P}_{22}^{-1}
\mathcal{P}_{21})\tilde \psi_j\in
\mathfrak{Hol}_{\sigma_0}(R^\perp)$, as in Lemma~\ref{mero2}.
Let $\psi_j$ be the singular part of
$\tilde \psi_j - \mathcal{P}_{22}^{-1}\mathcal{P}_{21}\tilde \psi_j$.
In this last function, $\mathcal{P} _{22}^{-1}\mathcal{P}_{21}
\tilde \psi_j$ has values in $K^\perp$, while $\tilde \psi_j$ has
values in $K$, so the order of the pole of
$\psi_j$ is the same as that of $\tilde \psi_j$ (there is no
cancellation). Note that furthermore the order of the pole of
$\mathcal{P}_{22}^{-1}\mathcal{P}_{21}\tilde
\psi_j$ is lower than that of $\tilde \psi_j$ because
$\mathcal{P}_{22}^{-1}\mathcal{P}_{21}$ vanishes at $\sigma=\sigma_0$.
Let $\mu_j$ be the order of the pole of $\psi_j$.

\begin{proposition}\label{mero20}
Let $\hat u$ be an $H^m(Y;E)$-valued meromorphic function with pole at
$0$. Then $\mathcal{P}(\sigma)(\hat u(\sigma))$ is holomorphic if and
only if there are $\C$-valued polynomials $p_j(\sigma)$ of degree
$\mu_j-1$ such that $\hat u-\sum_{j=1}^d p_j(\sigma)\psi_j(\sigma)$ is
holomorphic. Thus if $\hat f$ is holomorphic and
$\hat u=\mathcal{P}(\sigma)^{-1}(\hat f)$, then $\hat u$ is meromorphic
with singularity of the form $\sum_{j=1}^d p_j(\sigma)\psi_j(\sigma)$.
\end{proposition}
\begin{proof}
Suppose $\hat f = f\oplus g$ is a holomorphic function with values
in $R^\perp \oplus R$ and let
$\hat u = u\oplus v =\mathcal{P}(\sigma)^{-1}(f+g)$,
$\sigma \ne \sigma_0$, decomposed according to $K\oplus K^\perp$, so
 \begin{equation*}\begin{split}
 \mathcal{P}_{11}u + \mathcal{P}_{12}v &=f\\
 \mathcal{P}_{21}u + \mathcal{P}_{22}v &=g
 \end{split}\end{equation*}
 From the second equation, $v=\mathcal{P}_{22}^{-1}(g-\mathcal{P}_{21}u)$,
which replaced in the first gives
\begin{equation*}
 (\mathcal{P}_{11}  -
 \mathcal{P}_{12}\mathcal{P}_{22}^{-1}\mathcal{P}_{21})u
 = f -\mathcal{P}_{12}\mathcal{P}_{22}^{-1}g.
\end{equation*}
Since the $\tilde \beta_j=(\mathcal{P}_{11} -
\mathcal{P}_{12}\mathcal{P}_{22}^{-1}
\mathcal{P}_{21}) \tilde \psi_j$ are holomorphic near $\sigma_0$ and
independent at $\sigma_0$, there are $f_j$, $q_j \in \mathfrak {Hol}_0(\C)$
such that
$f=\sum f_j \tilde\beta_j$,
$\mathcal{P}_{12}\mathcal{P}_{22}^{-1}g =
\sum q_j\tilde\beta_j$ (the $q_j$ vanish at $\sigma_0$ because
$\mathcal{P}_{12}$ does). Then $u=\sum_j (f_j - q_j)\tilde
\psi_j$. Replacing this in the expression for $v$ gives
\begin{equation*}
 v=\mathcal{P}_{22}^{-1}g +\sum_j (f_j - q_j)
\mathcal{P}_{22}^{-1}\mathcal{P}_{21}\tilde \psi_j
\end{equation*}
so
 \begin{equation*}\begin{split}
 u+v&=\sum (f_j - q_j)(\tilde \psi_j -
\mathcal{P}_{22}^{-1}\mathcal{P}_{21}\tilde \psi_j) +
\mathcal{P}_{22}^{-1}g\\
 &=\sum_j p_j(\sigma) \psi_j + h,
 \end{split}\end{equation*}
where the $p_j$ are polynomials and $h$ is holomorphic.
\end{proof}

Note that each $\psi_j$ can be written as
\begin{equation*}
 \psi_j(\sigma,y)=\sum_{\ell=0}^{\mu_j-1}
 \frac{1}{\sigma^{\mu_j-\ell}}\, \psi_{j\ell}(y),
\end{equation*}
 with smooth sections $\psi_{j\ell}$ of $E\to Y$.

\smallskip
\begin{center} {\bf Appendix: Saturated domains} \end{center}
\smallskip

Let $H$ be a Hilbert space, $\Omega\subset\C$ open, and $S\subset\Omega$
a finite set. Let $\mathfrak M_{\Omega,S}(H)$ be the space of meromorphic
$H$-valued functions on $\Omega$ with poles in $S$, and let
$\mathfrak{Hol}_{\Omega}(H)$ be the subspace of holomorphic
elements. Multiplication by a holomorphic function $f(\sigma)$ defines an
operator 
\begin{equation*}
 f(\sigma): \mathfrak M_{\Omega,S}(H)\to\mathfrak M_{\Omega,S}(H) 
 \text{ such that } f(\sigma)\mathfrak{Hol}_{\Omega}(H) \subset
 \mathfrak{Hol}_{\Omega}(H),
\end{equation*}
so it induces an operator on
$\mathfrak M_{\Omega,S}(H)/ \mathfrak{Hol}_{\Omega}(H)$ also
denoted $f(\sigma)$. 
\begin{definition}
A subspace of
$\mathfrak M_{\Omega,S}(H)/\mathfrak{Hol}_{\Omega}(H)$ 
which is invariant under multiplication by $f(\sigma)=\sigma$ will be 
called saturated.
\end{definition}
A saturated subspace $\hat{\mathcal E}$ is thus a module over the ring $\C[\sigma]$. 

\begin{lemma}\label{saturation2}
Let $S=\set{\sigma_1,\dots,\sigma_s}$. If 
$\hat{\mathcal E}\subset \mathfrak M_{\Omega,S}(H)/\mathfrak{Hol}_{\Omega}(H)$
is saturated, then there are saturated spaces $\hat{\mathcal E}_{j}\subset
\mathfrak M_{\Omega,\set{\sigma_j}}(H)/\mathfrak {Hol}_{\Omega}(H)\subset \mathfrak
M_{\Omega,S}(H)/\mathfrak{Hol}_{\Omega}(H)$ such that
\begin{equation}\label{SaturatedSplit}
 \hat{\mathcal E} =\hat{\mathcal E}_{1}\oplus\cdots\oplus\hat{\mathcal E}_{s}.
\end{equation}
\end{lemma}
\begin{proof} 
For every $j$ there is a polynomial $q_j(\sigma)$ such that
$q_j(\sigma_j)=1$ and $q_j(\sigma_k)=0$ for $k\ne j$, with equalities 
satisfied to a sufficiently high order. Take 
$\hat{\mathcal E}_{j}=q_j(\sigma)\hat{\mathcal E}$, which is saturated because $\sigma
q_j(\sigma) = q_j(\sigma)\sigma$. Finally, note that
$q_1+\cdots+q_s=1$  to high order at each $\sigma_j$. For more details see the proof
of the next lemma.
\end{proof}

\begin{lemma}\label{saturation1}
A finite dimensional space 
$\hat{\mathcal E}\subset \mathfrak M_{\Omega,S}(H)/\mathfrak{Hol}_{\Omega}(H)$
is saturated if and only if it is invariant under multiplication by 
$\tau^{\im \sigma}$ for $\tau>0$.
\end{lemma}

\begin{proof} Suppose first that $\hat{\mathcal E}$ is saturated, so by the previous
lemma, $\hat{\mathcal E}=\hat{\mathcal E}_1\oplus\cdots\oplus\hat{\mathcal E}_{s}$
with  saturated spaces $\hat{\mathcal E}_{j}\subset
\mathfrak M_{\Omega,\set{\sigma_j}}(H)/\mathfrak {Hol}_{\Omega}(H)$. It is then
enough to prove that each $\hat{\mathcal E}_{j}$ is invariant under multiplication by
$\tau^{\im \sigma}$. But this is clear, since $\tau^{\im \sigma}$ is a polynomial
plus an entire function vanishing to high order at $\sigma_j$. 

Suppose now that $\hat{\mathcal E}$ is invariant under multiplication by $\tau^{\im
\sigma}$ for any $\tau>0$. We will first reduce the problem to the situation
where $\hat{\mathcal E} \subset \mathfrak M_{\Omega,\set{\sigma_0}}(H)/\mathfrak
{Hol}_{\Omega}(H)$. Let the integers $\mu_j$ be chosen so that for any representative
$\psi$ of en element of $\hat{\mathcal E}$, $(\sigma-\sigma_j)^\mu_j$ is regular at
$\sigma_0$. It is possible to find such numbers because $\hat{\mathcal E}$ is finite
dimensional. For any $(\zeta_1,\dots,\zeta_s)\in \C^s$ with $\zeta_j\ne \zeta_k$ if
$j\ne k$, let $\wp_{jk}(\zeta)=\left[\frac{
\zeta-\zeta_k}{\zeta_j-\zeta_k}\right]^{\mu_k}$, let
$p_j(\zeta;\zeta_1,\dots,\zeta_s)=\prod_{k\ne j}\wp_{jk}$. Then $p_j$ vanishes to
order $\mu_k$ at $\zeta_k$, $k\ne j$, and has value $1$ at $\zeta_j$, so
$p_j=1+(\zeta-\zeta_j)a_j$. The $a_j$ are polynomials in $\zeta$ whose coefficients
depend holomorphically on the $\zeta_j$.  Let
$b_j=\sum_{\ell=0}^{\mu_j-1} (-1)^\ell [(\zeta-\zeta_j)a_j]^\ell$. This is again a
polynomial in $\zeta$, and so is $q_j=b_j p_j$. Thus there are
polynomials in $\zeta$ with coefficients depending on the $\zeta_j$, 
$h_{j,k}(\zeta;\zeta_1,\dots,\zeta_s)$, $j,k=1,\dotsc,s$, such that 
\begin{equation*}
Q_j(\zeta;\zeta_1,\dots,\zeta_s) = \delta_{jk} + (\zeta-\zeta_k)^{\mu_k}
h_{j,k}(\zeta;\zeta_1,\dots,\zeta_s)
\end{equation*}
Let now $q_j(\sigma,\tau) =
Q_j(\tau^{\im\sigma};\tau^{\im\sigma_1},\dots,\tau^{\im\sigma_s})$, defined for those
$\tau$ for which the numbers $\tau^{\im\sigma_j}$ are distinct. Since this is a
polynomial in $\tau^{\im\sigma}$, and since $\hat{\mathcal E}$ is invariant under
multiplication by $\tau^{\im\sigma}$, each $q_j$ defines a linear map
$\pi_j:\hat{\mathcal E} \to \hat{\mathcal E}$. Since $q_j(\sigma,\tau) =
\delta_{jk} + (\sigma-\sigma)^{\mu_k} \tilde h_{j,k}(\sigma,\tau)$, $\pi_j
\psi_j\in \mathfrak M_{\Omega,\set{\sigma_j}}(H)/\mathfrak {Hol}_{\Omega}(H)$ and
$\pi_j\circ\pi_k=\delta_{jk}\pi_j$, and since $\sum_{j=1}^s q_j=1$, $\sum_j\pi_j=I$.
Thus with $\hat{\mathcal E}_j=\pi_j(\hat{\mathcal E})$ we get 
\begin{equation*}
\hat{\mathcal E}=\hat{\mathcal E}_1\oplus\cdots\oplus\hat{\mathcal E}_{s}
\end{equation*}
The spaces $\hat{\mathcal E}_j$ are invariant under multiplication by
$\tau^{\im\sigma}$ because $\tau^{\im\sigma} q_j=q_j\tau^{\im\sigma}$. The $\pi_j$
are independent of $\tau$.

Suppose now that $\hat{\mathcal E}\subset  \mathfrak
M_{\Omega,\set{\sigma_0}}(H)/\mathfrak {Hol}_{\Omega}(H)$ is invariant under
multiplication by $\tau^{\im\sigma}$ for all $\tau$ ($\tau = e$ suffices). Let
$\lambda(w)=\sum_{\ell=1}^N
\frac{(-1)^{\ell+1}}{\ell} w^\ell$, so that  $\lambda(e^\zeta-1)=\zeta+\zeta^{\mu}
h(\zeta)$, with $h(\zeta)$ entire, and some integer $\mu$ depending on $N$ which
may be assumed as large as desired by taking $N$ large enough. Let
$q(\sigma)=\lambda(\tau^{\im\sigma}\tau^{-\im\sigma_0}-1)$. Then 
 \begin{equation*}
 q(\sigma)=\im
(\sigma-\sigma_0)\log\tau + (\sigma-\sigma_0)^{\mu} h(\sigma,\tau)
 \end{equation*}
with $h(\sigma,\tau)$ holomorphic in $\sigma$. If $\psi \in \hat{\mathcal E}$, then
$q\psi\in \hat{\mathcal E}$ because $q$ is a polynomial in $\tau^{\im\sigma}$. But
with $\mu$ large enough, $q\psi = (\sigma-\sigma_0)\log\tau\,\psi$. Thus $\mathcal E$
is saturated. 
\end{proof}

One can also give a proof using that if $\hat{\mathcal E}$ is invariant
under multiplication by $e_t(\sigma) = e^{\im \sigma t}$ for any $t$ then $e_t$
defines a one parameter group on $\mathcal E$. This approach involves the topology of
$\mathfrak M_{\Omega,\set{\sigma_0}}(H)$ (to prove that the group is continuous, hence
differentiable). The proof given is better because it is elementary.


\section{Canonical Pairing}\label{CanonicalPairing}

Suppose that $K$ and $R^\perp$ are hermitian finite dimensional 
vector spaces.
Define a pairing
 \begin{equation*}
 \iota_{\sigma_0,K}: \mathfrak M_{\sigma_0}(K) \times \mathfrak
 M_{\overline \sigma_0}(K) \to \mathfrak M_{\sigma_0}(\C)
 \end{equation*}
by
 \begin{equation*}
 \mathfrak M_{\sigma_0}(K) \times \mathfrak M_{\overline
 \sigma_0}(K)\ni (u,v) \mapsto
 \iota_{\sigma_0,K}(u,v)=(u(\sigma),v(\overline
 \sigma)) \in \mathfrak M_{\sigma_0}(\C).
 \end{equation*}
and likewise a pairing $\iota_{\sigma_0,R^\perp}$ associated with
$R^\perp$. Let $\mathcal{P}(\sigma):K\to R^\perp$ be defined and holomorphic
in a neighborhood of $\sigma_0\in \C$.
Define $P^\star(\sigma) = P(\overline\sigma)^*$, where $*$ denotes the
pointwise adjoint of $\mathcal P:K\to R^\perp$. $P^\star$ is holomorphic in a
neighborhood of $\overline\sigma_0$. Then with the induced map
$\mathcal{P}^\star: \mathfrak M_{\overline\sigma_0}(R^\perp)
\to \mathfrak M_{\overline \sigma_0}(K)$ we have
 \begin{equation*}
 \iota_{\sigma_0,R^\perp}(\mathcal{P}(\sigma) u(\sigma), v(\sigma)) =
 \iota_{\sigma_0,K}(u(\sigma),\mathcal{P}^\star(\sigma)v(\sigma)).
 \end{equation*}
Furthermore, define
$\Theta:\mathfrak M_{\sigma_0}(\C)\to\mathfrak M_{\overline\sigma_0}(\C)$
by $\Theta(f)(\sigma) = \overline{f(\overline \sigma)}$, and likewise
$\Theta:\mathfrak M_{\overline \sigma_0}(\C)\to \mathfrak M_{\sigma_0}(\C)$.

\begin{lemma}
Let $\beta_1,\dots,\beta_d\in \mathfrak{Hol}_{\sigma_0}(R^\perp)$
be such that the $\beta_j(\sigma_0)$ are independent.
Then there are $\tilde \beta_1,\dots,\tilde \beta_d \in
\mathfrak {Hol}_{\overline \sigma_0}(R^\perp)$ such that
 \begin{equation*}
 \iota_{\sigma_0,R^\perp}(\beta_i,\tilde \beta_j) = \delta_{ij}.
 \end{equation*}
 If $\sigma_0$ is real, then 
 \end{lemma}
\begin{proof}
Let $b_i=\beta_i(\sigma_0)$ and write
$\beta_i=b_i+(\sigma-\sigma_0)\sum_k a_{ik}b_k$ with $a_{ik} \in
\mathfrak{Hol}_{\sigma_0}(\C)$. We seek
$\tilde \beta_j=\sum h_{\ell j} b_\ell$ with
$h_{\ell j}\in \mathfrak{Hol}_{\overline\sigma_0}(\C)$. We need
 \begin{equation*}
 \iota_{\sigma_0,R^\perp}(\beta_i,\tilde \beta_j) =
 \sum_{\ell} \Theta(h_{\ell j})(b_i,b_\ell) +
(\sigma-\sigma_0)\sum_{k,\ell} a_{ik}\Theta(h_{\ell j}) (b_k,b_\ell) =
\delta_{ij},
 \end{equation*}
 or, with the matrices $H=[h_{\ell j}]$, $A=[a_{ik}]$,
$B=[(b_k,b_\ell)]$,
 \begin{equation*}
 B\Theta(H) +(\sigma - \sigma_0) A B\Theta(H) = I
 \end{equation*}
 so set
 \begin{equation*}
 \Theta(H) = [I+(\sigma-\sigma_0)B^{-1}A B)]^{-1}B^{-1}.
 \end{equation*}
 \end{proof}

\begin{lemma}\label{mero7}
Let $\mathcal{P}(\sigma):K\to R^\perp$ be defined and
holomorphic near $\sigma_0$, invertible for $\sigma \ne \sigma_0$,
$\mathcal{P}(\sigma_0)=0$. Let $\psi_j\in \mathfrak M_{\sigma_0}(K)$
have independent leading coefficients and be such that
$\beta_j=\mathcal{P}\psi_j\in \mathfrak{Hol}_{\sigma_0}(R^\perp)$,
with the $\beta_j(\sigma_0)$ linearly independent. Let $\tilde \beta_j
\in \mathfrak{Hol}_{\overline\sigma_0}(R^\perp)$ be such that
 \begin{equation*}
 \iota_{\sigma_0,R^\perp}(\beta_i,\tilde \beta_j) = \delta_{ij}.
 \end{equation*}
 Let $\mu_j$ be the order of the pole of $\psi_j$, let $\tilde
\beta^\star_j = (\sigma-\sigma_0)^{\mu_j}\psi_j$, so
$\tilde \beta^\star_j \in\mathfrak{Hol}_{\sigma_0}(K)$ and the
$\tilde\beta^\star_j(\sigma_0)$ are independent. Let $\beta^\star_j \in
\mathfrak {Hol}_{\overline\sigma_0}(K)$ be such that
 \begin{equation*}
 \iota_{\sigma_0,K}(\tilde \beta^\star_i,\beta^\star_j) = \delta_{ij}.
 \end{equation*}
Then $\mathcal{P}^\star \tilde \beta_j = (\sigma-\overline
\sigma_0)^{\mu_j} \beta^\star_j$ so with
 \begin{equation*}
 \psi^\star_j= \frac 1 {(\sigma-\overline\sigma_0)^{\mu_j}}\tilde \beta_j
 \end{equation*}
we have
 \begin{equation*}
 \mathcal{P}^\star (\psi^\star_j) = \beta^\star_j.
 \end{equation*}
Clearly the leading coefficients of the $\psi^\star_j$ are
independent, as are the $\tilde\beta_j(\overline\sigma_0)$, thus the 
multiplicities $\mu_j^\star$ for $\mathcal{P}^\star$ are the same as those for
$\mathcal{P}$, $\mu_j^\star=\mu_j$.
\end{lemma}
\begin{proof}
We have $\mathcal{P} \tilde \beta^\star _j
=(\sigma-\sigma_0)^{\mu_j}\beta_j$, so
 \begin{equation*}\begin{split}
 \iota_{\sigma_0,R^\perp}(\mathcal{P} \tilde \beta^\star_j,\tilde
 \beta_k) &= (\sigma-\sigma_0)^{\mu_k}\delta_{jk}\\
 &=(\sigma-\sigma_0)^{\mu_k}\iota_{\sigma_0,K}(\tilde
 \beta^\star_j,\beta^\star_k)\\
 &=\iota_{\sigma_0,K}(\tilde
 \beta^\star_j,(\sigma-\overline \sigma_0)^{\mu_k}\beta^\star_k)
 \end{split}\end{equation*}
 The last expression must be $\iota_{\sigma_0,K}(\tilde
\beta^\star_j,\mathcal{P}^\star\tilde \beta_k)$, so
$\mathcal{P}^\star\tilde \beta_k
 =(\sigma-\overline \sigma_0)^{\mu_k}\beta^\star_k$.
\end{proof}

With the assumptions on $\mathcal{P}(\sigma)$ as in the previous lemma let
\begin{equation*}
 u\in\mathcal{P}^{-1}(\mathfrak{Hol}_{\sigma_0}(K)) \quad\text{and}\quad
 v \in(\mathcal{P}^\star)^{-1}(\mathfrak{Hol}_{\overline\sigma_0}(R^\perp)).
\end{equation*}
With suitable small $\varepsilon$ (depending on representatives of $u$
and $v$) we get a number
 \begin{equation*}
 [u,v]_{\mathcal P,\sigma_0} = \frac
 {1}{2\pi}\oint_{|\sigma-\sigma_0|=\varepsilon}\!
 \iota_{\sigma_0,R^\perp}(\mathcal{P}u,v)\,d\sigma.
 \end{equation*}
 The circle of integration is oriented counterclockwise. If $v$ is holomorphic, then
$[u,v]_{\mathcal P,\sigma_0}=0$ because
$\mathcal{P}u$ is holomorphic. If
$u$ is holomorphic, then also $[u,v]_{\mathcal P,\sigma_0}=0$, since
 \begin{equation*}
 \frac {1}{2\pi}\oint_{|\sigma-\sigma_0|=\varepsilon}\!
 \iota_{\sigma_0,R^\perp}(\mathcal{P}u,v)\,d\sigma=
 \frac {1}{2\pi}\oint_{|\sigma-\sigma_0|=\varepsilon}\!
 \iota_{\sigma_0,R^\perp}(u,\mathcal{P}^\star v)\,d\sigma
 \end{equation*}
 Thus $[\cdot,\cdot]_{\mathcal P,\sigma_0}$ defines a pairing
 \begin{equation*}
 [\cdot,\cdot]_{\mathcal
 P,\sigma_0}^\flat:\hat{\mathcal E}_{\sigma_0}\times
 \hat{\mathcal E}^\star_{\overline \sigma_0}\to \C,
 \end{equation*}
 where
 \begin{equation}\label{E_sigma0}
 \begin{gathered}
 \hat{\mathcal E}_{\sigma_0}=\mathcal{P}^{-1}
 (\mathfrak{Hol}_{\sigma_0}(K))/\mathfrak {Hol}_{\sigma_0}(K),\\
 \hat{\mathcal E}_{\overline\sigma_0}^\star
 =(\mathcal{P}^\star)^{-1}(\mathfrak{Hol}_{\overline
 \sigma_0}(R^\perp))/\mathfrak {Hol}_{\overline \sigma_0}(R^\perp).
 \end{gathered}
 \end{equation}

\begin{theorem}\label{mero5}
$[\cdot,\cdot]_{\mathcal P,\sigma_0}^\flat$ is a nonsingular paring
of vector spaces.
\end{theorem}
\begin{proof}
Pick $\psi_j \in \mathfrak M_{\sigma_0}(K)$ such that
$\mathcal{P}\psi_j = \beta_j\in \mathfrak {Hol}_{\sigma_0}(R^\perp)$,
with the leading coefficients of the
$\psi_j$ forming a basis of $K$ and with the $\beta_j(\sigma_0)$
forming a basis of $R^\perp$. Let $\tilde\beta_j$, $\tilde\beta_j^\star$,
and $\beta_j^\star$ be as in Lemma~\ref{mero7}.
According to the proof of Lemma~\ref{mero8},
if $u \in \hat{\mathcal E}_{\sigma_0}$ and $v \in \hat{\mathcal E}_{\overline
\sigma_0}^\star$ then $u$ and $v$ are represented by
 \begin{equation}\label{mero10}
 u=\sum_{i=1}^d \sum_{k=0}^{\mu_i-1} (\sigma-\sigma_0)^k u_{ik}\psi_i
 \quad\text{and}\quad
 v=\sum_{j=1}^d \sum_{\ell=0}^{\mu_j-1} (\sigma-\overline \sigma_0)^\ell
 v_{j\ell}\psi^\star_j
 \end{equation}
with constant $u_{ik}$ and $v_{j\ell}$. Now,
 \begin{equation*}
 \mathcal{P}(u) = \sum_{i=1}^d \sum_{k=0}^{\mu_i-1}
 (\sigma-\sigma_0)^k u_{ik}\beta_i
 \end{equation*}
 so
 \begin{equation*}\begin{split}
 \iota_{\sigma_0,R^\perp}(\mathcal{P}u,v) &=
 \sum_{i,j=1}^d \sum_{k,\ell=0}^{\mu_i-1} (\sigma-\sigma_0)^{k+\ell}
 u_{ik}\overline v_{j\ell}\;\iota_{\sigma_0,R^\perp}(\beta_i,\psi^\star_j)\\
 &= \sum_{i,j=1}^d \sum_{k,\ell=0}^{\mu_i-1}
 (\sigma-\sigma_0)^{k+\ell-\mu_j} u_{ik}\overline
 v_{j\ell}\;\iota_{\sigma_0,R^\perp}(\beta_i,\tilde \beta_j)\\
 &=\sum_{j=1}^d \sum_{k,\ell=0}^{\mu_j-1}
 (\sigma-\sigma_0)^{k+\ell-\mu_j} u_{jk}\overline v_{j\ell}
 \end{split}\end{equation*}
 Thus
 \begin{equation}\label{mero9}
 [u,v]_{\mathcal P,\sigma_0}^\flat = \im \sum_{j=1}^d
 \sum_{\;k+\ell-\mu_j=-1}\! u_{jk}\overline v_{j\ell}
 =\im \sum_{j=1}^d \sum_{k=0}^{\mu_j-1} u_{jk}\overline v_{j,\mu_j-k-1}
 \end{equation}
 If $[u,v]_{\mathcal P,\sigma_0}^\flat=0$ for all $v$, pick $v$ as above
so that $v_{j\ell}=u_{j,\mu_j-\ell-1}$, $j=1,\dots,d$,
$\ell=0,\dots,\mu_j-1$. Then
 \begin{equation*}
 [u,v]_{\mathcal P,\sigma_0}^\flat=\im \sum_{j=1}^d \sum_{k=0}^{\mu_j-1}
 u_{jk}\overline u_{jk} = 0
 \end{equation*}
 implies $u=0$.
\end{proof}

\begin{remark}\label{mero24} In the situation of the proposition and with the notation
in the proof, suppose that all $\mu_j$ are even. Let 
 \begin{equation*}
 U=\set{\sum_{j=1}^d\sum_{\,\ell=\mu_j/2}^{\mu_j-1} (\sigma-\sigma_0)^\ell
    u_{j\ell} \psi_j\st u_{j\ell}\in \C}
 \end{equation*}
which can be regarded as a subspace of $\hat{\mathcal E}_{\sigma_0}$. Likewise, since
the $\mu_j^\star$ associated with $\mathcal P$ are equal to the $\mu_j$'s, define
 \begin{equation*}
 V=\set{\sum_{j=1}^d\sum_{\,\ell=\mu_j/2}^{\mu_j-1} (\sigma-\overline \sigma_0)^\ell
    v_{j\ell} \psi_j^\star \st v_{j\ell}\in \C}
  \end{equation*}
 which again can be regarded as a subspace of $\hat{\mathcal
E}_{\overline\sigma_0}^\star$. It follows from \eqref{mero9} that 
 \begin{equation*}
 [u,v]_A=0 \quad\text{ if } u \in U, v\in V
 \end{equation*}
 and therefore, by dimensional considerations and the proposition itself, the
orthogonal of $U$ in $\hat{\mathcal E}_{\overline \sigma_0}^\star$ is $V$. 

The spaces $U$ and $V$ are independent of the $\psi_j$ used to represent them, as long
as these functions are chosen according to Lemma~\ref{mero2}. Indeed, if
$\set{\psi_i'}_{i=1}^d$ is another such choice, then according to Lemma~\ref{mero11},
possibly after reordering, we can write
\begin{equation*}
 \psi_i' = \sum_{\set{j\st \mu_j < \mu_i}}(\sigma-\sigma_0)^{\mu_j-\mu_i}f_{ij}\psi_i
+ \sum_{\set{j\st \mu_j \geq \mu_i}}f_{ij}\psi_i
\end{equation*}
and 
\begin{multline*}
 (\sigma-\sigma_0)^{\mu_j/2} \psi_i' = \sum_{\set{j\st \mu_j <
\mu_i}}(\sigma-\sigma_0)^{(\mu_j-\mu_i)/2}f_{ij}\,(\sigma-\sigma_0)^{\mu_i/2}\psi_i
\\+
\sum_{\set{j\st \mu_j \geq
\mu_i}}(\sigma-\sigma_0)^{(\mu_j-\mu_i)/2} f_{ij}\,(\sigma-\sigma_0)^{\mu_i/2}\psi_i
\end{multline*}
\end{remark}

\begin{lemma}\label{mero25} Suppose $\sigma_0$ is real, let $\psi_1,\dotsc,\psi_d\in
\mathfrak M_{\sigma_0}(K)$ have independent leading coefficients and be such that
their orders $\mu_j$ form a nonincreasing sequence. Then there are holomorphic
functions
$f_{ij}$ with $f_{ij}=0$ if $i>j$ such that 
\begin{equation}\label{mero25.1}
\iota_{\sigma_0,K}(\sum_{k=1}^i
(\sigma-\sigma_0)^{\mu_k-\mu_i}f_{ik}\psi_k,\sum_{\ell=1}^j
(\sigma-\sigma_0)^{\mu_\ell-\mu_j}f_{j\ell}\psi_\ell)=
(\sigma-\sigma_0)^{-\mu_i-\mu_j}\delta_{ij}
\end{equation}
\end{lemma}
This lemma combined with Lemma~\ref{mero11} says that when $\sigma_0$ is real, the
$\tilde \beta_i^\star$ in Lemma~\ref{mero7} can be assumed to form an
``orthonormal'' system: $\iota_{\sigma_0,K}(\tilde \beta_i^\star,\tilde
\beta_j^\star) = \delta_{ij}$.
\begin{proof} We apply the Gram-Schmidt orthonormalization process. There is no loss
of generality if we assume $\sigma_0=0$.  Write $\tilde \beta_i^\star=\sigma^\mu_j
\psi_j$. The holomorphic function $h(\sigma) = \iota_{0,K}(\tilde
\beta_1^\star,\tilde\beta_1^\star)$ is positive when $\sigma$ is real, so there is
$k(\sigma)$ positive defined for $\sigma$ real (close to $0$) such that $k^2 =
h$. Since $k$ is real analytic, it has a holomorphic extension to a
neighborhood of $0$. Since $k(\sigma)\overline{k(\overline \sigma)} = h(\sigma)$
holds when $\sigma$ is real, equality holds also for complex  $\sigma$ near $0$.
Let $f_{11}=1/k$. Then \eqref{mero25.1} holds for $i$, $j=1$, and we replace $\psi_1$
with $k^{-1}\psi_1$ and each $\psi_i$, $i > 1$, by $\psi_i -
(\sigma-\sigma_0)^{\mu_1-\mu_i}\iota_{\sigma_0,K}(\psi_i,\psi_1)h^{-1}\psi_1$. Now
repeat the process with $\psi_2$ and $\psi_i$ with $i>2$. 
\end{proof}

\smallskip
\begin{center} {\bf Appendix: Selfadjoint subspaces} \end{center}
\smallskip

Suppose for the rest of this section that $R^\perp=K$ and $\sigma_0\in \R$. Motivated
by Proposition~\ref{orthoadjoint} a subspace $\hat{\mathcal E}'_{\sigma_0}$ of
$\hat{\mathcal E}_{\sigma_0}$ will be called $\mathcal P$-selfadjoint (or just
selfadjoint if there is no risk of confusion) if
\begin{equation*}
 \hat{\mathcal E}'_{\sigma_0}=\set{u\in \hat{\mathcal E}_{\sigma_0}\st 
 [u,v]_{\mathcal P,\sigma_0}^\flat=0 \text{ for all }
 v\in \hat{\mathcal E}'_{\sigma_0}}.
\end{equation*}
In other words, $\hat{\mathcal E}'_{\sigma_0}$ is selfadjoint if
$\hat{\mathcal E}'_{\sigma_0}=(\hat{\mathcal E}'_{\sigma_0})^\perp$
with respect to $[\cdot,\cdot]_{\mathcal P,\sigma_0}^\flat$.

Let $\mathcal P(\sigma)$ be defined and holomorphic
near $\sigma_0$ (real). We call $\mathcal P$ selfadjoint if
$\mathcal P^\star(\sigma) = \mathcal P(\sigma)$ near $\sigma_0$.
If $\mathcal P$ is selfadjoint we say that $\mathcal P$ is positive if for each
real $\sigma\ne \sigma_0$ close to $\sigma_0$, $\mathcal P(\sigma):K\to K$ is
nonnegative.

\begin{lemma}\label{EvenPowers}
Let $\mathcal P$ be defined near $\sigma_0\in \R$, selfadjoint, positive, with
$\mathcal P(\sigma_0)=0$. Let $\psi_1,\dots,\psi_d \in \mathcal P^{-1}(\mathfrak
{Hol}_{\sigma_0}(K))$ be chosen as in Proposition~\ref {mero6}. Then the numbers
$\mu_j$ are even.
\end{lemma}
\begin{proof}
Without loss of generality we assume that $\sigma_0=0$. We use the notation of
Proposition~\ref{mero6} and assume that the $\psi_j$ are ordered so that the $\mu_j$
are nonincreasing in $j$ and define $\tilde \mu_i$ and $s_i$ as in \eqref{mero4}. We
replace the $\psi_j$ by suitable linear combinations of themselves to arrange that
with $\tilde\beta^\star_j=\sigma^{\mu_j}\psi_j$ we have $\iota_{\sigma_0,K}(\sum_{k}(\tilde \beta_i^\star,\tilde
\beta_j^\star) = \delta_{ij}$. This does not change the multiplicities $\mu_j$. Thus
$\tilde\beta^\star_j$ is holomorphic and the $\tilde \beta^\star_j(0)$ form
a  basis of $K$ (orthonormal), and if $\mathcal P\psi_j=\beta_j$, then
$\beta_j = \sum p_{j\ell} \tilde \beta^\star_\ell$ for some
holomorphic functions $p_{j\ell}$. Thus
\begin{equation*}
 \mathcal P\tilde \beta^\star_j = \sigma^{\mu_j} \sum_{\ell} p_{j\ell}
 \tilde \beta^\star_\ell.
\end{equation*}
Since both the vectors $\beta_j(0)$ and the $\tilde \beta^\star_j$ give bases of $K$,
the matrix $[p_{j\ell}(0)]$ is nonsingular. We have
\begin{equation*}
 \iota_{0,K}(\mathcal P\tilde \beta^\star_j,\tilde \beta^\star_k)
 = \sigma^{\mu_j}p_{jk}(\sigma)
\end{equation*}
and
\begin{equation*}
 \iota_{0,K}(\tilde \beta^\star_j,\mathcal P\tilde \beta^\star_k)
 = \sigma^{\mu_k}\overline {p_{kj}(\overline \sigma)}.
\end{equation*}
Since $\mathcal P$ is selfadjoint, these two functions are equal.
Consequently, if $\mu_k \geq \mu_j$, then
$p_{jk}(\sigma)=\sigma^{\mu_k-\mu_j}\overline{p_{kj}(\overline\sigma)}$. Thus
$[p_{jk}(0)]$ is an upper-triangular block matrix, the $i$-th diagonal block
corresponding to the indices $j$ such that $\mu_j=\mu_i$
($i=1,\dots,L$). Moreover, these diagonal blocks are selfadjoint matrices.
It follows that for each $i$ one can replace the $\psi_j$ ($j$ such that
$\mu_j = \tilde \mu_i$) by linear combinations of themselves with constant
coefficients, and
assume that the diagonal blocks of $[p_{jk}(0)]$ are diagonal themselves.
Let $p_{jj}(0)=\lambda_j$. Since $[p_{jk}(0)]$ is nonsingular, all
$\lambda_j$ are different from $0$. Now,
\begin{align*}
 \iota_{0,K}(\mathcal P\tilde \beta^\star_j,\tilde \beta^\star_j)
 &= \sigma^{\mu_j}(p_{jj} + \sigma\sum_{\ell} p_{j\ell}h_{\ell j})\\
 &= \sigma^{\mu_j}(\lambda_j + \sigma f_j)
\end{align*}
for some holomorphic function $f_j$. Since $\mathcal P$ is positive,
$\sigma^{\mu_j}(\lambda_j + \sigma f_j)$ is nonnegative for real $\sigma$,
and then necessarily $\mu_j$ is even (and $\lambda_j>0$).
\end{proof}

If $\hat{\mathcal E}_{\sigma_0}$ is a saturated subspace of
$\mathfrak M_{\sigma_0}(K)/\mathfrak{Hol}_{\sigma_0}(K)$ (see the appendix of
Section \ref{MeromorphicSolutions} for the definition) then there is a set of elements
$\psi_j\in
\hat{\mathcal E}_{\sigma_0}$, $j=1\dotsc,d$, such that the products
$(\sigma-\sigma_0)^\ell\psi_j$,
$\ell=0,\dots,\mu_j-1$ form a basis of $\hat{\mathcal E}_{\sigma_0}$. The proof of
this is that of Lemma~\ref{mero2} where only the saturation property was used.
We may further assume, as in that lemma, that if the $\psi_j$ are
represented as in \eqref{mero3}, then the $\psi_{j0}$ are independent and
the $\mu_j$ form a nonincreasing sequence.

\begin{proposition}\label{HalfDomains}
Let $\mathcal P$ be defined near $\sigma_0\in \R$, selfadjoint, positive and such that
$\mathcal P(\sigma_0)=0$, and let $\hat{\mathcal E}'_{\sigma_0}$ be a selfadjoint
saturated subspace of $\hat{\mathcal E}_{\sigma_0}$. Then every
$u\in \hat{\mathcal E}'_{\sigma_0}$ can be represented as
\begin{equation*}
 u =\sum_{j=1}^d \sum_{\,\ell=\mu_j/2}^{\mu_j-1} (\sigma-\sigma_0)^\ell
    u_{j\ell} \psi_j
\end{equation*}
with constant $u_{j\ell}$, where the  $\psi_1,\dots,\psi_d \in \mathcal
P^{-1}(\mathfrak {Hol}_{\sigma_0}(K))$ are chosen as in Proposition~\ref {mero6} and
the $\mu_j$ are the respective multiplicities, which by Lemma~\ref{EvenPowers} are
even. The space of such elements will be denoted 
$\hat{\mathcal E}_{\sigma_0,\frac 1{2}}$.
\end{proposition}
\begin{proof} We may assume that $\sigma_0=0$, without loss of generality.
By Remark \ref {mero24}, if we show that for some choice of $\psi_j$ as in the
statement the elements of $\hat{\mathcal E}'_{\sigma_0}$ can be represented as
stated, then for any such choice of $\psi_j$ they are represented as stated. We then
take advantage of Lemma~\ref{mero25} and assume that if
$\tilde\beta^\star_j=\sigma^{\mu_j}\psi_j$ then
$\iota_{\sigma_0,K}(\tilde\beta_j^\star,\tilde\beta_k^\star) = \delta_{ij}$. In the
notation of Lemma~\ref{mero7} we then have $\beta^\star_j=\tilde\beta_j^\star$
($\sigma_0$ is real). As in that lemma let $\beta_j=\mathcal P\psi_j$ and let
$\tilde \beta_j \in \mathfrak{Hol}_{0}(K)$ be such that
$\iota_{0,K}(\beta_j,\tilde \beta_k)=\delta_{ij}$. By Lemma~\ref{mero7}, $\mathcal
P^\star \sigma^{-\mu_j}\tilde \beta_j = \beta_j^\star$, but now the latter is
equal to $\tilde \beta_j^\star$, and $\mathcal P$ is selfadjoint, so
 \begin{equation*}
 \mathcal P(\sigma^{-\mu_j}\tilde \beta_j) = \tilde \beta_j^\star
 \end{equation*}
Since both the $\tilde \beta_j(0)$ and the $\tilde \beta_j^\star(0)$ form bases of
$K$, the functions $\psi_j^\star=\sigma^{-\mu_j}\tilde \beta_j$ satisfy the
conditions in Lemma~\ref{mero11}, so they are related as stated there, and from
Remark~\ref{mero24} we get that if 
 \begin{equation*}
 u =\sum_{j=1}^d \sum_{\,\ell=\mu_j/2}^{\mu_j-1} \sigma^\ell
    u_{j\ell} \psi_j,\qquad v =\sum_{j=1}^d \sum_{\,\ell=\mu_j/2}^{\mu_j-1}
\sigma^\ell
    v_{j\ell} \psi_j
 \end{equation*}
then 
 \begin{equation*}\label{NullPairing}
 [u,v]_{\sigma_0,\mathcal P} = 0.
 \end{equation*}
This said, we now show that if there is an element in $\mathcal E'_{\sigma_0}$
represented by
 \begin{equation*}
 u =\sum_{j=1}^d \sum_{\,\ell=\ell_j}^{\mu_j-1} \sigma^\ell
    u_{j\ell} \psi_j \in \hat{\mathcal E}'_{\sigma_0}
 \end{equation*}
 ($\ell_j\geq 0$) where for some $j$, $\ell_j < \mu_j/2$ with $u_{j,\ell_j}\ne 0$,
then 
$\hat {\mathcal E}_{\sigma_0}$ is not selfadjoint.  Thus let $\delta_j = \mu_j/2 -
\ell_j$, let
$\delta = \max_{j}\delta_j$ and suppose $\delta\geq 1$. Since $\hat {\mathcal
E}_{\sigma_0}$ is saturated both $\sigma^{\delta-1}u$ and $\sigma^{\delta}u$
represent elements in $\hat {\mathcal E}_{\sigma_0}$, and we will show that
$[\sigma^{\delta-1}u,\sigma^\delta u]_{0.\mathcal P}\ne 0$. Let
$J=\set{j:\delta_j=\delta}$. Then 
\begin{equation*}
 \sigma^{\delta-1} u
 \equiv \sum_{j\in J} u_{j\ell_j}  \sigma^{\mu_j/2-1}\psi_j
 + \sum_{j=1}^d \sum_{\ell=\mu_j/2}^{\mu_j-1}
 \sigma^\ell \tilde u_{j,\ell} \psi_j \quad\mod \mathfrak{Hol}_{0}(K)
\end{equation*}
with some $\tilde u_{j\ell}$. Write this as $u_0+ \tilde u_1$. Then 
\begin{equation*}
[u_0+u_1,\sigma
u_0+\sigma u_1]_{\mathcal P, \sigma_0}=[u_0,\sigma
u_0+\sigma u_1]_{\mathcal P, \sigma_0}
\end{equation*}
because of Remark~\ref{mero24}, 
\begin{equation*}
[u_0,\sigma u_0+\sigma u_1]_{\mathcal P, \sigma_0}=[\sigma u_0, u_0+ u_1]_{\mathcal
P, \sigma_0}
\end{equation*}
because multiplication by $\sigma$ is
selfadjoint, and finally, 
\begin{equation*}
[\sigma u_0, u_0+\sigma
u_1]_{\mathcal P, \sigma_0} = [\sigma u_0, u_0]_{\mathcal P, \sigma_0}
\end{equation*}
again by
Remark~\ref{mero24}. Thus $[\sigma^{\delta-1}u,\sigma^\delta]_{\mathcal P, \sigma_0}=
[u_0,\sigma u_0]_{\mathcal P, \sigma_0}$. As in the proof of
Proposition~\ref{HalfDomains}, $\mathcal P \beta_j^\star = \sigma^{\mu_j}\sum_{k=1}^d
p_{jk}\beta_k^\star$ with holomorphic $p_{jk}$
which because of the selfadjointess are such that 
$p_{jk}(\sigma)=\sigma^{\mu_k-\mu_j}\overline{p_{kj}(\overline\sigma)}$ if $\mu_k
\geq \mu_j$. We thus have (since $\psi_j=\sigma^{-\mu_j}\tilde\beta_j^\star$)
\begin{equation*}
\mathcal P u_0 = \sum_{j\in J} \sum_k u_{j,\ell_j} p_{jk}\sigma^{\mu_j/2-1}
\tilde \beta_k^\star
\end{equation*}
and so
\begin{align*}
\iota_{\sigma_0,K}(\mathcal P u_0,\sigma u_0) 
&= \sum_{j,j'\in J} \sum_{k=1}^d p_{j,k}
u_{j,\ell_j}\sigma^{\mu_j/2-1}\overline
u_{j',\ell_{j'}}\iota_{\sigma_0,K}(\beta_k,\sigma^{-\mu_{j'}/2}\psi_{j'})\\
&=\sum_{j,j'\in J}  p_{j,j'}
\sigma^{\mu_j/2-\mu_{j'}/2-1} u_{j,\ell_j}\overline
u_{j',\ell_{j'}}
\end{align*}
It is the residue of this what we will show is nonzero. Terms with 
$\mu_{j'}<\mu_{j}$ clearly do not contribute to the residue. For terms with
$\mu_{j'}>\mu_{j}$ we have
$p_{jj'}(\sigma)=\sigma^{\mu_{j'}-\mu_{j}}\overline{p_{j'j}(\overline\sigma)}$; such
terms again contribute nothing, and we conclude that
\begin{equation*}
 [u_0,\sigma u_0]_{\mathcal P,\sigma_0} = \im \sum_
 {\substack{j,j'\in J\\\mu_j=\mu_{j'}}}
 p_{j,j'}(0)
 u_{j,\ell_j}\overline
u_{j',\ell_{j'}}
\end{equation*}
The positivity of $\mathcal P$ now enters again: as pointed out at the
end of the proof of Lemma~\ref{EvenPowers}, the selfadjoint matrices $[p_{j,j'}]$
with $j,j'$ such that $\mu_j=\mu_{j'}$ are  positive definite. Thus $[u_0,\sigma
u_0]_{\mathcal P,\sigma_0}\ne 0$ since $u_{j,\ell_j}\ne 0$ for at least
for one $j$.
\end{proof}

\begin{remark}
Note that when the $\mu_j$ are even numbers the space
\begin{equation*}
 \hat{\mathcal E}_{\sigma_0,\frac 1{2}}=
\LinSpan_\C\set{(\sigma-\sigma_0)^\ell  \psi_j \st j=0,\dots,d;
\,\ell=\mu_j/2,\dots,\mu_j-1}
\end{equation*}
is a {\em canonical} saturated subspace of $\mathcal E_{\sigma_0}$,
even if the operator $\mathcal P$ is not selfadjoint.
The same holds for $\mathcal E^\star_{\overline\sigma_0,\frac 1{2}}$, 
and we have $(\mathcal E_{\sigma_0,\frac 1{2}})^\perp
= \mathcal E^\star_{\overline\sigma_0,\frac 1{2}}$.
\end{remark}


\section{Structure of the Adjoint Pairing}\label{DomainPairing}

Let now $E\to M$ be a vector bundle, $H^s_b(M;E)$ be the totally characteristic
Sobolev space of order $s$, defined as usual. Suppose
that $A=x^{-\nu}P\in x^{-\nu}\diff^m_b(M;E)$, $\nu>0$, is a $b$-elliptic cone
operator considered initially as
a densely defined unbounded operator
\begin{equation*}
 A: C_c^\infty (\Dot M;E)\subset x^{-\nu/2}L^{2}_b(M;E)\to x^{-\nu/2}L^{2}_b(M;E),
\end{equation*}
and define $\Dom_{\min}=\Dom_{\min}(A)$, $\Dom_{\max}=\Dom_{\max}(A)$ as in
Section~\ref{SectGround}. It is well known from the proof of the existence of
asymptotic expansions of solutions of $Pu=0$ (cf.~\cite{K1},
\cite{MM}, \cite{Sz91}) that if $u\in x^{-\nu/2}L^{2}_b(M;E)$
and $A u \in x^{-\nu/2}L^{2}_b(M;E)$, that is, $u\in \Dom_{\max}$, then $\hat u$
is meromorphic in $\Im\sigma > -\nu/2$ with values in $H^m(\partial M;E|_{\partial
M})$ and poles contained in 
\begin{equation}\label{SpectrumRainStrip}
 \Sigma'(A) = (\spec_b(A)-\im\mathbb N_0)\cap \set{\sigma\st-\nu/2<\Im \sigma<\nu/2}.
\end{equation}
 Moreover, if for $u$ as above, $\hat u$ is holomorphic in $\Im\sigma >-\nu/2$, then
in fact $u \in \Dom_{\max}\cap x^{\nu/2-\eps}H^m_b(M;E)$ for any $\eps>0$, so by
Proposition~\ref{DomMin}, $u\in \Dom_{\min}(A)$. Thus $\mathcal
E(A)=\Dom_{\max}/\Dom_{\min}$ is isomorphic to a certain subspace $\hat{\mathcal
E}(A)$ of 
 \begin{equation}\label{BigQuotient}
 \mathfrak M_{\Omega,\Sigma'(A)}(H^m(\partial M;E|_{\partial
M}))/\mathfrak{Hol}_{\Omega}(H^m(\partial M;E|_{\partial M}))
 \end{equation}
where $\mathfrak M_{\Omega,\Sigma'(A)}(H^m(\partial M;E|_{\partial M}))$ is the space
of meromorphic $H^m(\partial M;E|_{\partial M})$-valued functions on $\Omega =
\set{\sigma \st -\nu/2 < \Im\sigma}$ and
$\mathfrak{Hol}_{\Omega}(H^m(\partial M;E|_{\partial M}))$ is the subspace of
holomorphic elements. It is clear that the space in
\eqref{BigQuotient} is localizable on $\spec_b(A)$ in the sense that it is isomorphic
to the direct sum
 \begin{equation*}\label{SumSmallQuotient}
 \bigoplus_{\Sigma(A)}
 \mathfrak M_{\Omega,\Sigma'_{\sigma}}(H^m(\partial M;E|_{\partial
M}))/\mathfrak{Hol}_{\Omega}(H^m(\partial M;E|_{\partial M}))
 \end{equation*}
 where 
 \begin{equation}
 \Sigma(A)=\set{\sigma\in \spec_b(A)\st -\nu/2<\Im\sigma<\nu/2}
 \end{equation}
 and 
 \begin{equation}
 \Sigma'_\sigma=\set{\sigma -\im \ell\st \ell \in \mathbb N_0,\ \ell <\Im\sigma +
\nu/2}.
 \end{equation}
 It is also the case that $\hat{\mathcal E}(A)$ is localizable: we will
show that 
$\hat{\mathcal E}(A)$ is the direct sum of the spaces 
 \begin{equation*}
 \hat{\mathcal E}_{\sigma}(A) = 
 \hat{\mathcal E}(A)\cap \big[\mathfrak M_{\Omega,\Sigma_{\sigma}(A)}(H^m(\partial
M;E|_{\partial M}))/\mathfrak{Hol}_{\Omega}(H^m(\partial M;E|_{\partial M}))\big].
 \end{equation*}
 To see this, begin by writing 
  \begin{equation}\label{MeroSplit1}
 P=\sum_{k=0}^N P_k x^k + \tilde P_N x^N
 \end{equation}
where the $P_k$ have coefficients independent of $x$ near $\partial M$ and
$N=\min\set{k\in\mathbb N\st \nu\leq k}$. Let $\sigma_0\in \Sigma(A)=\spec_b(A)\cap
\set{\sigma|-\nu/2<\Im \sigma<\nu/2}$. By the
discussion immediately preceding Proposition~\ref{mero20} and the proposition itself
with $\mathcal P=\hat P_0$ at $\sigma_0$, there are elements $\psi_{\sigma_0,j,0,0}\in
\mathfrak M_{\sigma_0}(C^\infty(\partial M;E_{\partial M}))$ and positive integers
$\mu_{\sigma_0,j}$ (forming a nonincreasing sequence), $j=1\dotsc,d_{\sigma_0}$, such
that for any $u\in \Dom_{\max}$, if $\hat P_0(\sigma)\hat u(\sigma)$ is holomorphic
at $\sigma_0$, then  there are polynomials $p_j$ such that $\hat u -
\sum_{j=1}^{d_{\sigma_0,j}} p_j(\sigma)\psi_{\sigma_0,j,0,0}(\sigma)$ is holomorphic
at
$\sigma_0$. If $\Im\sigma_0 -\vartheta >-\nu/2$ define $\psi_{\sigma_0,j,0,\vartheta}$
inductively as the singular part at
$\sigma_0-\im \vartheta$ of
 \begin{equation*}
  - \hat P_0^{-1}(\sigma)\sum_{\zeta=0}^{\vartheta-1}
 \hat P_{\vartheta-\zeta}(\sigma)\psi_{\sigma_0,j,0,\zeta}(\sigma +\im
(\vartheta-\zeta)).
 \end{equation*}
 and for convenience define $\psi_{\sigma_0,j,0,\vartheta} = 0$ if $\Im\sigma_0
-\vartheta
\leq -\nu/2$. The largest index $\theta$ such that $\Im \sigma_0 -\theta>-\nu/2$
will be denoted by $N(\sigma_0)$. Define also
$\psi_{\sigma_0,j,\ell,\vartheta}$ (for
$\ell=0,\dotsc,\mu_j-1)$ to be the principal part of $(\sigma+\im \vartheta)^\ell
\psi_{\sigma_0,j,0,\vartheta}$ (at
$\sigma_0-\im\vartheta$), and finally, let 
 \begin{equation*}
 \Psi_{\sigma_0,j,\ell}=\sum_{\vartheta\geq 0} \psi_{\sigma_0,j,\ell,\vartheta}.
 \end{equation*}
 Then 
 \begin{equation*}
 \sum_{k=0}^N P_k(\sigma) \Psi_{\sigma_0,j,\ell}(\sigma+\im k)
 \end{equation*}
 is holomorphic in $\Im\sigma > -\nu/2-\eps$ for any sufficiently small $\eps>0$. The
claim is now that the images in
$\hat{\mathcal E}(A)$ of the
$\Psi_{\sigma_0,j,\ell}$ form a basis (over $\C$). This is easy to prove
beginning with the fact that the $\psi_{\sigma_0,j,\ell,0}$ form a basis of $\hat
P_0(\sigma)^{-1}/\mathfrak{Hol}_{\sigma_0}$. Once this is proved, the assertion about
$ \hat{\mathcal E}(A)$ being localizable is clear. 
  
 Note that one may multiply each $\psi_{\sigma_0,j,\ell,\vartheta}$ by a suitable
entire function which is equal to $1$ to high order at $\sigma_0-\im \vartheta$ so
that the resulting function is, modulo an entire function, the Mellin transform of an
element in $x^{-\nu/2}H^\infty(M;E)$. This will not change the
fact that the images in $\hat{\mathcal E}(A)$ of the modified
$\Psi_{\sigma_0,j,\ell}$ form a basis. An immediate convenient consequence is 
 
 \begin{lemma} For each $u\in \Dom_{\max}(A)$ there is $u_0\in \Dom_{\min}(A)$ such
that $(u-u_0)\hat{\ }$ is meromorphic on $\C$ with poles only in $\Sigma'(A)$. 
 \end{lemma}
 
 \begin{definition}\label{PoleQuotient} For $\sigma_0 \in \Sigma(A)$,
$\Dom_{\sigma_0}(A)$ is the space of elements $u \in \Dom_{\max}(A)$ such that $\hat
u(\sigma)$ represents an element in
$\hat{\mathcal E}_{\sigma_0}(A)$, that is, $\hat u$ has poles at most at
$\sigma_0-\im\vt$ for $\vt = 0,\dots, N(\sigma_0)$. Further,
 \begin{equation*}
 \mathcal{E}_{\sigma_0}(A)=\Dom_{\sigma_0}(A)/\Dom_{\min}(A).
 \end{equation*}
 For $\sigma_0\in\spec_b(A)\cap \set{\Im\sigma=0}$, and if all the
multiplicities $\mu_{\sigma_0,j}$ associated with $\sigma_0$ are even, we let
$\Dom_{\sigma_0,\frac 12}(A)$ be the space of elements $u\in \Dom_{\sigma_0}(A)$
 such that 
\begin{equation*}
 \hat u \mod \mathfrak{Hol}(\Im\sigma>-\eps) \text{ belongs to }\hat{\mathcal
E}_{\sigma_0,\frac 12}
\end{equation*}
for any small $\eps>0$. The space $\hat{\mathcal E}_{\sigma_0,\frac 12}$
is the one defined in  Proposition~\ref{HalfDomains}, now with $\mathcal P=\hat P_0$.
 \end{definition}
 
 Thus, modulo $\mathfrak {Hol}(\Im\sigma>-\nu)$, the Mellin transform of an element
$u\in \Dom_{\sigma_0}(A)$ can be written as $\hat u(\sigma) =
\sum_{\vartheta=0}^{N(\sigma_0)} \psi_\vartheta(\sigma)$
 where the $\psi_\vartheta(\sigma)$ have poles only at $\sigma_0-\im\vartheta$, and
from the fact that $\widehat {Pu}(\sigma) = \sum_{k,\vartheta} \hat
P_k(\sigma)\psi_\vartheta(\sigma+\im\vartheta)$ is holomorphic in
$\Im
\sigma>-\nu/2$ one deduces that
 \begin{equation}\label{meroHol}
  \sum_{\vartheta=0}^\ell P_{\ell-\vartheta}(\sigma) 
\psi_{\vartheta}(\sigma+\im(\ell-\vartheta)) \text{ is holomorphic if } \ell
\leq N(\sigma_0).
 \end{equation}
 Likewise, if $v\in \Dom_{\sigma_0^\star}(A^\star)$ then $\hat v (\sigma)=
\sum_{\vartheta=0}^{N(\sigma_0^\star)} \psi_\vartheta^\star(\sigma)$ where now
 \begin{equation}\label{meroHolDual}
  \sum_{\vartheta=0}^\ell P_{\ell-\vartheta}(\sigma+\im(\ell-\vartheta)) 
\psi_{\vartheta}^\star(\sigma+\im(\ell-\vartheta)) \text{ is holomorphic if } \ell
\leq N(\sigma_0^\star)
 \end{equation}
 since from \eqref{MeroSplit1}, 
 \begin{equation}\label{MeroSplit1Adjoint}
P^\star=\sum_{\ell=0}^{N-1} x^\ell P_\ell^\star + x^N\tilde P_N^\star.
 \end{equation}

The closed extensions of $A:\Dom_{\min}(A)\subset x^{-\nu/2}L^{2}_b(M;E)\to
x^{-\nu/2}L^2_b(M;E)$ are in one to one correspondence with the subspaces of $\mathcal
E(A)$, therefore with the subspaces of $\hat{\mathcal E}(A)$. Since we are interested
in duality and selfadjoint extensions, we will now turn our attention towards
understanding the pairing $[u,v]_A$ for $u\in\Dom_{\max}(A)$ and
$v\in\Dom_{\max}(A^\star)$ as a pairing of elements of $\hat{\mathcal E}(A)$. In the
following theorem and its proof, the pairing in the integrands is that of
$L^2(\partial M;E|_{\partial M})$. 

\begin{theorem}\label{BasicPairing}
Let $A=x^{-\nu}P\in x^{-\nu}\diff^m_b(M;E)$ be $b$-elliptic, let
$\sigma_0\in\Sigma(A)$ and $\sigma_0^\star\in\Sigma(A^\star)$, suppose
 $u\in \Dom_{\max}(A)$ and $v\in \Dom_{\max}(A^\star)$
are such that
\begin{equation*}
  \hat u=\sum_{\vt=0}^{N(\sigma_0)} \psi_\vt \quad\text{and}\quad
  \hat v=\sum_{\vt=0}^{N(\sigma_0^\star)} \psi_\vt^\star \quad \mod
\mathfrak{Hol}(\Im\sigma>-\nu/2),
\end{equation*}
where $\psi_\vt$ has a pole only at $\sigma_0-\im \vt$, and
$\psi_\vt^\star$ only at $\sigma_0^\star-\im \vt$, in other words, $u\in
\Dom_{\sigma_0}(A)$ and $v\in \Dom_{\sigma_0^\star}(A^\star)$.  If $\sigma_0$ is not
of the form
$\overline{\sigma_0^\star}+\im \tau$ with $\tau \in \mathbb N_0$, then
$[u,v]_A=0$. Otherwise, if $\sigma_0=\overline{\sigma_0^\star} + \im\tau$ for
$\tau\in\mathbb{N}_0$, then
 \begin{equation*}
 [u,v]_A = \frac{1}{2\pi}
 \sum_{\vt=0}^{\tau}
 \oint_{\gamma_{\vt}}(
  \psi_{\tau-\vt}(\sigma) ,{\sum_{\vartheta'=0}^{\vt} \hat P_{\vt-\vartheta'}^\star
 (\overline{\sigma-\im (\vt-\vartheta')})\,
 \psi_{\vartheta'}^\star(\overline{\sigma-\im(\vt-\vartheta')})})\,d\sigma,
 \end{equation*}
where $\gamma_\vt=\gamma_0+\im\vt$ and $\gamma_0$ is a positively oriented
simple closed curve surrounding $\overline{\sigma_0^\star}$.
In particular, if $\tau=0$, i.e., $\sigma_0=\overline{\sigma_0^\star}$, then
 \begin{equation}\label{mero23}
 [u,v]_A = \frac{1}{2\pi}
 \oint_{\gamma_{0}}(
 \psi_0(\sigma) ,{\hat P_0^\star (\overline{\sigma})
 \psi_0^\star(\overline {\sigma})})\,d\sigma.
 \end{equation}
\end{theorem}

\begin{proof} 
For general $u\in \Dom_{\max}(A)$ and $v\in \Dom_{\max}(A^\star)$ and
$\omega\in C^\infty(M)$ supported near the boundary and equal to $1$ near the
boundary one has
 \begin{equation*}
 [u,v]_A=[\omega u,\omega v]_A
 \end{equation*}
because $(1-\omega)u\in \Dom_{\min}(A)$ if $u\in \Dom_{\max}(A)$, and analogously
for $(1-\omega)v$. Recall that the Mellin transform was defined using a cut-off
function like $\omega$. Suppose $P$ is written as in \eqref{MeroSplit1} where the
$P_\ell$ have coefficients independent of $x$ near the boundary, say, in a
neighborhood of the closure of the support of $\omega$. Then, using the expression
for $A^\star$ obtained from \eqref{MeroSplit1Adjoint}, we have
 \begin{multline*}
 [\omega u,\omega v]_A = (x^{-\nu}\sum_{\ell=0}^{N-1} P_\ell\,x^\ell
\omega u,\omega v)_{x^{-\nu/2}L^2_b} - (\omega u,x^{-\nu}\sum_{\ell=0}^{N-1} x^\ell
P_\ell^\star \omega v)_{x^{-\nu/2}L^2_b}\\
 + (x^{-\nu} \tilde P_N\, x^N \omega u,\omega v)_{x^{-\nu/2}L^2_b} - (\omega
u,x^N\tilde P_N^\star x^{-\nu}\omega v)_{x^{-\nu/2}L^2_b}.
 \end{multline*}
The last two terms cancel out since $x^Nu\in x^{\nu/2}H^m_b$ can be
approximated from $C_c^\infty(\Dot M;E)$ in $x^{\nu/2}H^m_b$ norm. Thus
 \begin{equation*}
 [u,v]_A=(x^{-\nu}\sum_{\ell=0}^{N-1} P_\ell\,x^\ell
\omega u,\omega v)_{x^{-\nu/2}L^2_b} - (\omega u,x^{-\nu}\sum_{\ell=0}^{N-1} x^\ell
P_\ell^\star \omega v)_{x^{-\nu/2}L^2_b}.
 \end{equation*}
  It is always the case that $\hat u$ has no poles on $\Im \sigma = \nu/2$, and
adding a suitable element of $\Dom_{\min}(A)$ to $u$ we may assume that
$\hat u$ has no poles on $\Im \sigma = -\nu/2$. A similar remark applies to $\hat v$,
and we may and will assume that neither $\hat u$ nor $\hat v$ has poles on $\Im
\sigma=-\nu/2$. For $\eps>0$ let 
\begin{equation}\label{Betas}
 \beta_0 =-\frac{\nu}{2} +\eps \quad\text{and}\quad
\beta_k =\frac{\nu}{2}-\eps - N + k \;\text{ for } k=1,\dots,N
\end{equation}
We take $\eps>0$ sufficiently small so that $\beta_0<\beta_1$. There is $\eps_0>0$
such that for any $\eps<\eps_0$, if $\sigma_0\in\spec_b(A)\cup\spec_b(A^\star)$ and
$-\nu/2<\Im\sigma<\nu/2$, then for any $\ell\in \mathbb N_0$, the point
$\sigma - \im\ell$ does not lie on a line $\Im\sigma=\beta_k$. That is, no $u\in
\Dom_{\max}(A)$ or $v\in \Dom_{\max}(A^\star)$ has poles on a line $\Im \sigma =
\beta_k$. Fix one such $\eps$ and let
\begin{equation*}
  S_k=\set{\sigma\in\C\st \beta_k\leq \Im \sigma \leq \beta_{k+1}}.
\end{equation*}
These strips partition $\set{ \sigma\in\C\st -\nu+\eps\leq \Im \sigma < \nu/2-\eps}$.
We now show that
\begin{equation}\label{MeroPair1}
 [u,v]_A =
 \sum_{k=0}^{N-1} \oint_{\partial S_k} (\sum_{\ell=0}^{N-k-1} \hat
P_\ell(\sigma)\hat u(\sigma+\im\ell) ,
  {\hat v(\overline{\sigma})})\,d\sigma
\end{equation}
For any $\ell<N$ and $\eps$ small, we have, on the one hand,
\begin{align*}
\quad(\frac 1 {x^\nu}P_\ell x^\ell u,v)
 &= (x^{\eps-\nu}P_\ell x^\ell u,x^{-\eps}v) \\
 &= \frac 1{2\pi}
  \int_{\Im\sigma=\frac\nu{2}} (\hat P_\ell(\sigma+\im(\eps-\nu))
  \hat u(\sigma+\im(\eps-\nu+\ell)) ,{\hat v(\sigma-\im\eps)})\,d\sigma\\
 &= \frac 1{2\pi}
  \int_{\Im\sigma=-\frac\nu{2}} (\hat P_\ell(\sigma+\im\eps)
  \hat u(\sigma+\im(\eps+ \ell)) ,
  {\hat v(\overline{\sigma+\im\eps})})\, d\sigma\\
 &= \frac 1{2\pi}
  \int_{\Im\sigma=\beta_0} (\hat P_\ell(\sigma)\hat u(\sigma+\im\ell) ,
  {\hat v(\overline{\sigma})})\, d\sigma\\
\intertext{and}
(u,\frac 1 {x^\nu} x^\ell P_\ell^\star v)
 &= (x^{-\eps}u, x^{\eps-\nu}x^\ell P_\ell^\star v) \\
 &= \frac 1{2\pi} \int_{\Im\sigma=\frac\nu{2}}
  ( \hat u(\sigma-\im\eps), {\hat P_\ell^\star(\sigma+ \im(\eps-\nu+ \ell))
  \hat v (\sigma+\im(\eps-\nu+\ell))})\, d\sigma \\
 &= \frac 1{2\pi} \int_{\Im\sigma=\frac\nu{2}- \ell}
  (\hat u(\sigma-\im\eps+\im\ell),
  {\hat P_\ell^\star(\overline{\sigma-\im\eps})
  \hat v (\overline{\sigma-\im\eps})})\, d\sigma\\
 &= \frac 1{2\pi} \int_{\Im\sigma=\beta_{N-\ell}}
  (\hat P_\ell(\sigma)\hat u(\sigma+\im\ell) ,
  {\hat v(\overline{\sigma})})\, d\sigma
\end{align*}
so
\begin{align*}
 [u,v]_A&=\sum_{\ell=0}^N \frac{1}{2\pi}
  \oint_{\partial \set{\sigma\st \beta_0\leq \Im\sigma\leq \beta_{N-\ell} }} (\hat
P_\ell(\sigma)\hat u(\sigma+\im\ell) ,
  {\hat v(\overline{\sigma})})\,d\sigma \\
  &=\sum_{\ell=0}^N \sum_{k=0}^{N-\ell-1}
   \oint_{\partial S_k} (\hat P_\ell(\sigma)\hat u(\sigma+\im\ell) ,
  {\hat v(\overline{\sigma})})\,d\sigma\\
  &=\sum_{k=0}^{N-1} \oint_{\partial S_k} (\sum_{\ell=0}^{N-k-1} \hat
P_\ell(\sigma)\hat u(\sigma+\im\ell) ,
  {\hat v(\overline{\sigma})})\,d\sigma
\end{align*}

Let $\sigma_0\in\Sigma$ and suppose $\hat u=\sum_{\vt=0}^\infty \psi_\vt$,
where $\psi_\vt$ has a pole only at $\sigma_0-\im \vt$ and $\psi_\vt = 0$ if
$\vt>N(\sigma_0)$, where $N(\sigma)$ be the number $k$ such that
$\sigma\in S_{k}$. Thus $\hat P_\ell(\sigma)
\psi_\vt(\sigma+\im \ell)$ has a pole at $\sigma_0-\im(\vt+\ell)$ if at
all, and the poles of
\begin{equation}\label{mero21}
 \sum_{\ell=0}^{N-k-1}
 \hat P_\ell(\sigma) \hat u(\sigma+\im \ell) =
 \sum_{\ell=0}^{N-k-1}\sum_{\vt=0}^\infty \hat P_\ell(\sigma)
 \psi_\vt(\sigma+\im \ell),
\end{equation}
if any, that lie in $\beta_k<\Im \sigma < \beta_{k+1}$, come from indices
$\vt,\ell$ with
\begin{equation*}
 \beta_k < \Im \sigma_0 -(\vt+\ell)<\beta_{k+1},
\end{equation*}
that is, $\vt+\ell = N(\sigma_0)-k$. So in \eqref{mero21} only the terms
\begin{equation*}
 \sum_{\substack{ \vt+\ell=N(\sigma_0)-k\\
 0\leq \ell\leq N-k-1\\ 0\leq \vt \leq N(\sigma_0) }}
 \hspace{-.8em} \hat P_\ell(\sigma) \psi_\vt(\sigma+\im \ell)
\end{equation*}
may produce poles in $S_k$. If $k>N(\sigma_0)$ there are no poles.
If $k\leq N(\sigma_0)$, this is
\begin{equation*}
 \sum_{\vartheta=0}^{N(\sigma_0)-k} \hat P_{N(\sigma_0)-k - \vartheta}(\sigma)\,
 \psi_\vartheta(\sigma + \im (N(\sigma_0)-k - \vartheta))
\end{equation*}
 (since $N(\sigma_0)<N$). This is in fact holomorphic, as stated in \eqref{meroHol}.
Thus in
\eqref{MeroPair1}, the integrals
 \begin{equation*}
 \oint_{\partial S_k} (\sum_{\ell=0}^{N-k-1} \hat
P_\ell(\sigma)\hat u(\sigma+\im\ell) ,
  {\hat v(\overline{\sigma})})\,d\sigma
 \end{equation*}
are evaluated as residues on the
conjugates
 of the poles of $\hat v$,
\begin{equation}\label{mero22}
 [u,v]_A
 = \frac 1{2\pi}\sum_{k=0}^{N-1}\sum_s
 \oint_{\gamma_{s,k}}  ( \sum_{\ell=0}^{N-k-1}\hat P_\ell(\sigma) \hat
 u(\sigma+\im \ell) ,{\hat v(\overline \sigma)})\,d\sigma ,
\end{equation}
where the $\gamma_{s,k}$ are simple closed positively oriented curves
surrounding (and separating) the conjugates of the poles of
$\hat v(\sigma)$ in the strip $S_k$.

Suppose now that also $\hat v = \sum_{\vt=0}^\infty \psi_\vt^\star$,
where $\psi_\vt^\star$ has a pole only at $\sigma_0^\star - \im \vt$.
Here $\sigma_0^\star \in \Sigma(A^\star) = \overline{\Sigma(A)}$
and as before, $\psi_\vt^\star=0$ if $\vt > N(\sigma_0^\star)$.
Thus 
 \begin{equation*}
 \sigma\mapsto ( \sum_{\ell=0}^{N-k-1}\hat P_\ell(\sigma) \hat
 u(\sigma+\im \ell) ,{\hat v(\overline \sigma)})
 \end{equation*}
has poles at $\overline{\sigma_0^\star} +\im \vt$, $\vt=0,\dots,N(\sigma_0^\star)$,
and the pole in $S_k$ satisfies
$\beta_k < \Im\overline{\sigma_0^\star} +\vt < \beta_{k+1}$,
that is, $k = N(\overline{\sigma_0^\star}) +\vt$.
In particular, there are poles only in the strips with $k$ satisfying
\begin{equation*}
  N(\overline{\sigma_0^\star})\leq k\leq N-1.
\end{equation*}
Pick a positively oriented simple
closed curve $\gamma_0$ surrounding $\overline{\sigma_0^\star}$, let
$\gamma_\vt=\gamma_0+\im \vt$. Using this, the right hand side of
\eqref{mero22} becomes
\begin{equation*}
 \frac 1{2\pi}\sum_{k=N(\overline{\sigma_0^\star})}^{N-1}
 \oint_{\gamma_{k-N(\overline{\sigma_0^\star})}}( \sum_{\ell=0}^{N-k-1}
 \hat P_\ell(\sigma) \hat u(\sigma+\im \ell) ,{\hat v(\overline
 \sigma)})\,d\sigma.
\end{equation*}
Now replace $\sigma +\im\ell$ by $\tilde\sigma$. The resulting expression
is (dropping the tilde)
\begin{multline*}
 \frac 1{2\pi}\sum_{k=N(\overline
 {\sigma_0^\star})}^{N-1}\sum_{\ell=0}^{N-k-1}
 \oint_{\gamma_{k+\ell-N(\overline{\sigma_0^\star})}}( \hat
 P_\ell(\sigma-\im \ell) \hat u(\sigma) ,{\hat v(\overline
 {\sigma-\im\ell})})\,d\sigma \\
 = \frac 1{2\pi} \sum_{\vt=0}^{N-N(\overline
 {\sigma_0^\star})-1} \oint_{\gamma_{\vt}}( \hat u(\sigma),
 {\sum_{k=0}^{\vt} \hat P_{\vt-k}^\star (\overline{\sigma-\im (\vt-k)})
 \,\hat v(\overline {\sigma-\im(\vt-k)}\,)}) \,d\sigma
\end{multline*}
after reorganizing (notice that
$N-N(\overline{\sigma_0^\star})-1 \le N({\sigma_0^\star})$).
Now, if $\hat v(\sigma)=\sum_{\vartheta'\geq 0}\psi_{\vartheta'}^\star$ is as above,
then for any given $\vartheta$ the only terms in
\begin{equation*}
 {\sum_{k=0}^{\vartheta} \hat P_{\vartheta-k}^\star ({\sigma+\im
 (\vartheta-k)})\sum_{\vartheta'=0}^{N(\sigma_0^\star)}\psi_{\vartheta'}^\star(
 {\sigma+\im(\vt-k)})}
\end{equation*}
which may contribute to the integral along $\gamma_\vartheta$ are those which in
principle have poles at $\sigma_0^\star-\im\vartheta$, namely those in the sum
\begin{equation*}
 {\sum_{\vartheta'=0}^{\vartheta} \hat P_{\vartheta-\vartheta'}^\star ({\sigma+\im
(\vartheta-\vartheta')})
 \psi_{\vartheta'}^\star({\sigma+\im(\vartheta-\vartheta')})}.
\end{equation*}
But according to \eqref{meroHolDual}, this has no poles at
${\sigma_0^\star}-\im\vartheta$ (or anywhere else, for that matter). So, if $\hat u$
has no poles in
$\set{\overline{\sigma_0^\star}+\im\vt\st \vt=0,\dots,N(\sigma_0^\star)}$
then $[u,v]_A=0$. The only case where $[u,v]_A$ may be different from $0$
occurs when there are integers $\vt,\vt'\geq 0$ such that
$\sigma_0 -\im \vt = \overline{\sigma_0^\star} + \im\vt'$, that is,
if $\sigma_0=\overline{\sigma_0^\star}+\im\tau$ for some nonnegative
integer $\tau$, in which case
\begin{equation*}
 [u,v]_A = \frac 1{2\pi}\sum_{\vt=0}^{\tau}
 \oint_{\gamma_{\vt}}(\psi_{\tau-\vt}(\sigma) ,
 {\sum_{\vartheta'=0}^{\vt} \hat P_{\vt-\vartheta'}^\star(\overline{\sigma-\im
(\vt-\vartheta')})\,
 \psi_{\vartheta'}^\star(\overline{\sigma-\im(\vt-\vartheta')})})\,d\sigma
\end{equation*}
as claimed in the theorem.
\end{proof}

If $\sigma_0$, $\sigma_1 = \overline{\sigma_0^\star}\in \Sigma(A)$ are such that
$\sigma_0 = \sigma_1+\im\tau$ with some integer $\tau>0$, then if $u\in
\Dom_{\sigma_0}(A)$ and $v\in \Dom_{\sigma_0^\star}(A^\star)$ it may happen that
$[u,v]_A \ne 0$. 

Since $[u,v]_A=0$ if $u\in \Dom_{\min}(A)$ and $v\in \Dom_{\max}(A^\star)$,
or if $u\in \Dom_{\max}(A)$ and $v\in \Dom_{\min}(A^\star)$, there is a well defied
pairing of elements of $\mathcal E(A)$ and $\mathcal E(A^\star)$.

\begin{theorem}\label{FlatPairing}
$[\cdot,\cdot]_A^\flat$ is a nonsingular paring of
$\mathcal{E}_{\sigma_0}(A)$ and $\mathcal{E}_{\overline
\sigma_0}(A^\star)$.
\end{theorem}
\begin{proof}
We work with \eqref{mero23}. $\hat P_0(\sigma)$ is a
closed operator $L^2(Y) \to L^2(Y)$ with domain $H^m=H^m(Y;E)$.
Let $K=\ker \hat P_0(\sigma_0)$, $R=\hat P_0(\sigma_0)(H^m(Y))$.
Decompose $P(\sigma)$ as
 \begin{equation*}
 \left[\begin{matrix} \hat P_{11}(\sigma) & \hat P_{12}(\sigma)\\
  \hat P_{21}(\sigma)&\hat P_{22}(\sigma)
  \end{matrix}\right]:
  \begin{matrix} K \\ \oplus \\ K^\perp\cap H^m \end{matrix}
  \to \begin{matrix} \,\,\,R^\perp \\ \oplus \\ R \end{matrix}
 \end{equation*}
for $\sigma$ near $\sigma_0$. Here $K^\perp$, $R^\perp$ are computed in
$L^2(Y)$. Since $\hat P_0(\sigma_0)$ is Fredholm,
$\hat P_0^\star(\overline \sigma_0)(H^m)=K^\perp$ and the analogous
decomposition for $\hat P(\overline \sigma)^*=\hat P^\star(\sigma)$ is
 \begin{equation*}
  \left[\begin{matrix} \hat P_{11}^\star(\sigma) &
  \hat P_{21}^\star(\sigma)\\
  \hat P_{12}^\star(\sigma)&\hat P_{22}^\star(\sigma)
  \end{matrix}\right]:
  \begin{matrix} R \\ \oplus \\ R^\perp\cap H^m\end{matrix}
  \to \begin{matrix} K \\ \oplus \\ \,\,K^\perp \end{matrix}
 \end{equation*}
near $\overline \sigma_0$. Since $\Ind \hat P(\sigma)=0$,
$\dim K=\dim R^\perp$.

Let $u \in \mathcal{E}_{\sigma_0}(A)$ represent an element in 
$\Dom_{\sigma_0}(A)$, let $\psi$ be the Mellin transform
of $\phi u$. The principal part of $\psi$ at $\sigma_0$ is the principal
part $\psi_0$ at $\sigma_0$ of a germ of the form
 \begin{equation*}
 \tilde \psi - \hat P_{22}^{-1}\hat P_{21}\tilde \psi
 \end{equation*}
where $\tilde \psi \in \mathfrak M_{\sigma_0}(K)$ is such that
 \begin{equation*}
 (\hat P_{11} - \hat P_{12}\hat P_{22}^{-1}\hat P_{21})\tilde \psi=\beta
 \end{equation*}
is holomorphic near $\sigma_0$. Likewise let $u^\star \in
\mathcal{E}_{\overline \sigma_0}(A^\star)$ represent an element in
$\Dom_{\overline \sigma_0}(A^\star)$, $\psi^\star$ the
Mellin transform of $\phi u^\star$. Again the principal part
$\psi_0^\star$ of $\psi^\star$ at $\overline \sigma_0$ is the principal
part at $\overline \sigma_0$ of germ of the form
 \begin{equation*}
 \tilde \psi^\star -
 (\hat P_{22}^{-1})^\star\hat P_{12}^\star\tilde\psi^\star
 \end{equation*}
where $\tilde\psi \in\mathfrak M_{\overline \sigma_0}(R^\perp)$
is such that
 \begin{equation*}
 (\hat P_{11}^\star -
 \hat P_{21}^\star(\hat P_{22}^{-1})^\star\hat P_{12}^\star)
 \tilde\psi^\star=\beta^\star
 \end{equation*}
is holomorphic, near $\overline\sigma_0$. Let
 \begin{gather*}
 \mathcal{P}=\hat P_{11} - \hat P_{12}\hat P_{22}^{-1}\hat P_{21}\\
 \mathcal{P}^\star = \hat P_{11}^\star -
 \hat P_{21}^\star(\hat P_{22}^{-1})^\star\hat P_{12}^\star
 \end{gather*}
 Then as discussed before the lemma,
 \begin{equation*}
 [u,v]_A =\frac 1{2\pi} \oint_{\gamma_0} (\psi(\sigma),
 \hat P_0^\star(\overline\sigma)\psi^\star(\overline\sigma))d\sigma
 \end{equation*}
with a positively oriented curve $\gamma_0$ surrounding $\sigma_0$
and no other pole. Since
\begin{align*}
\intertext{ \quad
$\hat P_0^\star(\tilde \psi^\star - (\hat P_{22}^{-1})^\star\hat
 P_{12}^\star \tilde \psi^\star)$
}
 &= \hat P_{11}^\star \tilde \psi^\star -
 \hat P_{21}^\star (\hat P_{22}^{-1})^\star\hat P_{12}^\star
 \tilde\psi^\star + \hat P_{12}^\star \tilde \psi^\star
 - \hat P_{22}^\star (\hat P_{22}^{-1})^\star
 \hat P_{12}^\star\tilde\psi^\star \\
 &= \mathcal{P}^\star \tilde \psi^\star
\end{align*}
and
\begin{align*}
 (\psi(\sigma),\mathcal{P}^\star(\overline \sigma)\tilde
 \psi^\star(\overline \sigma)) &=
 (\tilde \psi(\sigma)-\hat P_{22}(\sigma^{-1}\hat P_{21}(\sigma)\tilde
 \psi(\sigma),\mathcal{P}^\star(\overline \sigma)
 \tilde \psi^\star(\overline \sigma))\\
  &=(\tilde \psi(\sigma),\mathcal{P}^\star(\overline \sigma)
 \tilde\psi^\star(\overline \sigma))
\end{align*}
we have
\begin{equation*}
[u,v]_A =\frac 1{2\pi}\oint_{\gamma_0} (\tilde \psi(\sigma),
\mathcal{P}^\star(\overline \sigma)\tilde \psi^\star(\overline\sigma))d\sigma.
\end{equation*}
Thus the pairing of $\mathcal{E}_{\sigma_0}(A)$ and
$\mathcal{E}_{\overline \sigma_0}(A^\star)$ is the pairing of the
spaces associated to $\mathcal{P}$ at $\sigma_0$ and
$\mathcal{P}^\star$ at $\overline \sigma_0$ which Theorem~\ref{mero5}
asserts is nonsingular.
\end{proof}


\section{Friedrichs Extension}\label{FriedrichsExtension}

Suppose $A\in x^{-\nu}\diff^m_b(M;E)$, $\,A=x^{-\nu}P$, is $b$-elliptic,
symmetric and bounded from below by some $c_0\leq 0$, as an operator
$C_c^\infty(\Dot M;E)\subset x^{-\nu/2}L^2_b(M;E) \to x^{-\nu/2}L^2_b(M;E)$.
The domain of the Friedrichs extension is denoted $\Dom_F(A)$.  Recall (Definition
\ref{PoleQuotient}) that we denote by $\Dom_{\sigma_0}(A)$ the space of functions
$u\in \Dom_{\max}(A)$ such that $\hat u(\sigma)$ has poles at most at
$\sigma_0-\im\vt$ for $\vartheta = 0,1,\dots$, by
$\mathcal E_{\sigma_0}$ the quotient $\Dom_{\sigma_0}(A)/\Dom_{\min}(A)$ and by 
$[u,v]_A = (Au,v) - (u,A^\star v)$, as introduced in \eqref{DmaxPairing}. 

\begin{lemma}\label{Fried1}
$\Dom_{F}(A)$ contains all $u\in \Dom_{\max}(A)$ such
that $\hat u(\sigma)$ has no poles in $\set{\Im\sigma \geq 0}$.
That is, $\Dom_{\max}(A)\cap H^m_b(M;E)\subset\Dom_{F}(A)$.
\end{lemma}
\begin{proof}
We will show that if $u\in\Dom_{\max}(A)\cap H^m_b(M;E)$,
there is a sequence $\seq{u_n}$ in $C_c^\infty(\Dot M;E)$ such that
\[ c\|u-u_n\|^2_{x^{-\nu/2}L^2_b} + (A(u-u_n),u-u_n)_{x^{-\nu/2}L^2_b}
   \to 0\ \text{ as } n\to \infty.\]
This will imply $u\in \Dom_F(A)$ by  Lemma \ref{DFApproximation}.
Consider $P=x^{\nu}A$ as
an unbounded operator on $L^2_b(M;E)$. Since it is $b$-elliptic, we have
$H^m_b(M;E)\subset \Dom_{\min}(P)$. Therefore, if $u\in\Dom_{\max}(A)\cap
H^m_b(M;E)$,  there is a sequence $\seq{u_n}\subset C_c^\infty(\Dot M;E)$ such that
\[ \|u-u_n\|_{L^2_b(M;E)}\to 0\ \text{ and }\
\|P(u-u_n)\|_{L^2_b(M;E)}\to 0\ \text{ as } n\to\infty. \]
With this sequence we have
\[ \|u-u_n\|_{x^{-\nu/2}L^2_b}\to 0\ \text{ as } n\to\infty \]
since $L^2_b(M;E)\embed x^{-\nu/2}L^2_b(M;E)$. Also,
\begin{align*}
(A(u-u_n),u-u_n)_{x^{-\nu/2}L^2_b}
 &= (P(u-u_n),u-u_n)_{L^2_b}\\
 &\leq
  \|P(u-u_n)\|_{L^2_b} \|u-u_n\|_{L^2_b}
\end{align*}
so
\begin{equation*}
(A(u-u_n),u-u_n)_{x^{-\nu/2}L^2_b}\to 0\ \text{ as } n\to \infty
\end{equation*}
and the proof is complete.
\end{proof}

\begin{lemma}\label{Fried2} 
$\Dom_F(A)$ contains no $u\in \Dom_{\max}(A)$ such
that $\hat u(\sigma)$ has a pole in $\set{\Im\sigma> 0}$. Thus
$\Dom_F(A)\subset\Dom_{\max}(A)\cap x^{-\eps}H^m_b(M;E)$ for any $\eps >0$.
\end{lemma}
\begin{proof}
Let $\sigma_0\in\spec_b(A)$ be such that $\Im\sigma_0>0$.
Suppose that $u\in \Dom_{\sigma_0}\cap\Dom_F$. In particular, $[u,v]_A=0$
for all $v\in \Dom_F(A)$ since $A_F$ is selfadjoint. From the previous lemma we know
that  $\Dom_{\overline\sigma_0}(A)\subset\Dom_F(A)$, hence $[u,v]_A=0$ for all
$v\in \Dom_{\overline\sigma_0}(A)$, hence $u=0$ since by
Theorem~\ref{FlatPairing} the induced pairing
$[\cdot,\cdot]^\flat_A$ of $\mathcal{E}_{\overline\sigma_0}$ and
$\mathcal{E}_{\sigma_0}$ is nonsingular.
\end{proof}

As a consequence of these two lemmas we get
\begin{theorem}\label{Fried7} Suppose $A\in x^{-\nu}\diff^m_b(M;E)$ is $b$-elliptic
and semibounded. If 
\begin{equation*}
\spec_b(A)\cap\set{\Im\sigma=0}=\varnothing,
\end{equation*}
then the domain of the Friedrichs extension of $A$ is 
\begin{equation*}
\Dom_F(A)=\sum_{
 \substack{ \sigma\in \spec_b(A)\\  -\nu/2<\Im\sigma<0 }
 }\hspace{-.8em}
 \Dom_{\sigma}(A)
\end{equation*}
That is, $\Dom_F(A)=\Dom_{\max}(A)\cap H^m_b(M;E)$.
\end{theorem}

This finishes the discussion of the Friedrichs  extension of $A$ when $\spec_b(A)\cap
\set{\Im\sigma=0}$. In order to determine the domain of the Friedrichs extension if
$\spec_b(A)$ does contain real elements, we need two more ingredients. The first is
an invariance property, under certain circumstances (small $\nu$), of $\Dom_F(A)$,
which translates into the saturation property on the Mellin transform side. The second
is the positivity of the conormal symbol of $A$ when $A$ is bounded from below.  

Let $\omega\in C_c^\infty(\R)$ be a function with sufficiently small support, and equal
$1$ near the origin. Let $\phi$ be the flow of $X$, which we shall write
multiplicatively: the integral curve of $X$ through $p$ is $t\mapsto \phi_{e^t}(p)$.  
We can write
 \begin{equation*}
 \phi^*_{1/\tau}\mathfrak m = c_\tau^2 \mathfrak m,\quad \phi^*_{1/\tau}x =
\tau^{-1}\xi_\tau^2 x
 \end{equation*}
with smooth positive functions $c_\tau$ and $\xi_\tau$ which are equal to $1$ if
$\tau=1$, or if $x$ is close to $\partial M$ (how close depends on $\tau$), or in
complement of the support of $\omega$. Thus we have
 \begin{equation*}
 \int \phi_\tau^*f\, \mathfrak m = \int \phi_\tau^*f\,
\phi_\tau^*(\phi_{1/\tau}^*\mathfrak m) = \int f\, c_\tau^2 \mathfrak m
 \end{equation*}
 Recall that for sections $u$ of $E$, $\phi_\tau^* u$ is the section whose value at
$p$ is the result of parallel transport of $u(\phi_\tau(p))$ to $p$ along the curve
through $p$. The connection is compatible with the hermitian form on $E$, so for
sections $u$, $v$ of $E$,
 \begin{equation*}
 (\phi_\tau^* u, \phi_\tau^* v)_p = ( u, v)_{\phi_\tau(p)}.
 \end{equation*}
 Let $\gamma_\tau = \phi_\tau^* \frac{1}{\xi_\tau^{\nu} c_\tau}$, define
 \begin{equation*}
 \kappa_\tau u = \tau^{\nu/2} \gamma_\tau \phi_\tau^* u,\quad u\in C^\infty(M;E).
 \end{equation*}
 Then $\kappa_\tau$ defines an isometry
 \begin{equation*}
 x^{-\nu/2}L^2_b(M;E)\to x^{-\nu/2}L^2_b(M;E)
 \end{equation*}
 and $\kappa_\tau^*=\kappa_{1/\tau}$. On functions $f$, $\kappa_\tau$ is defined as
$\kappa_\tau f = \phi_\tau^* f$, so that if $f$ is a function and $u$ a section of
$E$ then $\kappa_\tau(fu)=\kappa_\tau(f)\kappa_\tau(u)$. 

A subspace $\Dom\subset x^{-\nu/2}L^2_b(M;E)$ is
$\kappa$-invariant if
\begin{equation*}
\kappa_\tau u\in\Dom \text{ for every $u\in\Dom$ and $\tau>0$}.
\end{equation*}
For example $C_c^\infty(M;E)$ is $\kappa$-invariant.  Let $A\in
x^{-\nu}\diff^m_b(M;E)$ be arbitrary, write $A= \sum_{k=0}^{N-1} x^k A_k + x^N
\tilde A_N$ where
$x^{\nu}A_k=P_k \in \diff^m_b(M;E)$ has coefficients independent of $x$ near
$\partial M$, and $ \tilde A_N\in x^{-\nu}\diff^m_b(M;E)$. Then
\begin{equation*}
 \kappa_\tau A \kappa^{-1}_\tau= 
 \sum_{k=0}^{N-1} \tau^{\nu-k} x^k A_k + \tau^{\nu-N}x^N \tilde A_{N,\tau} 
\end{equation*}
for some $\tilde A_{N,\tau}\in x^{-\nu}\diff^m_b(M;E)$. In particular, for 
$x$ near $\partial M$ and $\tau$ smaller than some $\tau_0$ (depending on $x$),
\begin{equation*}
 \kappa_\tau A_0 \kappa^{-1}_\tau  = \tau^\nu A_0 .
\end{equation*}
This identity and the $\kappa$-invariance of $C_c^\infty(M;E)$ easily imply that
the canonical domains $\Dom_{\min}(A_0)$, $\Dom_{\max}(A_0)$ and $\Dom_{\max}(A_0)\cap
x^\gamma H^m_b(M;E)$ are $\kappa$-invariant. If $A_0$ is symmetric and bounded from
below, then using Lemma \ref{DFApproximation} one also proves easily that the domain
$\Dom_F(A_0)$ of the Friedrichs extension is also
$\kappa$-invariant.

\begin{lemma}\label{K-invariant}
Let $A\in x^{-\nu}\diff^m_b(M;E)$ be $b$-elliptic. 
Let $\Dom\subset\Dom_{\max}(A)$ be a domain on which $A$ is closed. 
$\Dom$ is $\kappa$-invariant if and only if the finite dimensional space
\begin{equation*}
 \hat{\mathcal E}_{\Dom}=\set{\hat u:
u\in\Dom}/\mathfrak{Hol}(\Im\sigma>-\nu/2)
\end{equation*}
is a saturated space.
\end{lemma}
\begin{proof} If $u$ is a smooth function on $M$, then for 
$\tau>0$
\begin{align*}
\widehat{\kappa_\tau u}(\sigma) 
 &=\int x^{-\im\sigma} \omega(x)u(\tau x,y)\frac{dx}{x}\\
 &=\tau^{\im\sigma}\int x^{-\im\sigma} \omega(\tau^{-1}x)u(x,y)\frac{dx}{x}\\
 &=\tau^{\im\sigma}\hat u(\sigma) + \tau^{\im\sigma}\hat w_\tau(\sigma)
\end{align*}
for $w_\tau= (\omega(\tau^{-1}x)-\omega(x))u$.
Now, since $\omega(\tau^{-1}x)-\omega(x)$ is a smooth function supported in 
the interior of $M$, then $\hat w_\tau$ is 
an entire function, that is, $\widehat{\kappa_\tau u}(\sigma)-\tau^{\im\sigma}
\hat u(\sigma)$ is entire. The same conclusion holds when $u$ is a smooth
section of $E$:
 \begin{equation*}
 \widehat{\kappa_\tau u}(\sigma)  = \tau^{\im\sigma}\hat u(\sigma) \mod
\mathfrak{Hol}_{\C}(C^\infty(\partial M;E|_{\partial M}))
 \end{equation*} 
 This proves that 
$\widehat{\kappa_\tau u}\mod \mathfrak{Hol}(\Im\sigma>-\nu/2)$ is an
element of $\hat{\mathcal E}_{\Dom}$ if and only if
$\hat{\mathcal E}_{\Dom}$ is invariant under multiplication by $\tau^{\im\sigma}$,
i.e., if and only if
$\hat{\mathcal E}_{\Dom}$ is saturated, due to Lemma~\ref{saturation1}.
The assertion thus follows from the isomorphism between
$\Dom/\Dom_{\min}$ and $\hat{\mathcal E}_{\Dom}$ given by the Mellin transform.
\end{proof}

The second ingredient we need to determine the Friedrichs extension is the positivity
of the conormal symbols of operators bounded from below. This is standard but we
provide a proof.
\begin{lemma}\label{Positivity}
Let $A\in x^{-\nu}\diff^m_b(M;E)$ be $b$-elliptic, symmetric and bounded
from below. For every $\sigma\in\R$ the conormal symbol $\hat P_0(\sigma)$
of $A$ is nonnegative.
\end{lemma}
\begin{proof}
Suppose  $v\in C^{\infty}(\partial M;E)$. 
Let $\phi\in C_c^\infty(0,1)$ be such that $\int|\phi(x)|^2\frac{dx}x=1$,
and let
\begin{equation*}
 \phi_n(x)=\frac{1}{n^{1/2}}\phi(x^{1/n}), \quad n\in\mathbb N.
\end{equation*}
Then $\phi_n v\in\Dom_{\min}(A)$ since it is smooth and supported in the
interior of $M$. It is easy to prove that for real $\sigma$ one has
 \begin{equation*}
 (A(x^{\im \sigma}\phi_n v), x^{\im \sigma }\phi_n v)_{x^{-\nu/2}L^2_b(M;E)}
 \to (\hat P_0(\sigma )v,v)_{L^2(\partial M;E|_{\partial M})}\ \text{ as } n\to
\infty.
 \end{equation*}
Pick $c$ real such that $A-cI\ge 0$. The conormal symbol
of $A-cI$ is then the same as that of $A$. Thus
\begin{equation*}
 0\leq ((A-cI)x^{\im\sigma}\phi_n v,
 x^{\im\sigma}\phi_n v)_{x^{-\nu/2}L^2_b(M;E)}
 \to(\hat P_0(\sigma) v,v)_{L^2(\partial M;E|_{\partial M})}.
\end{equation*}
\end{proof}


 From Lemma \ref{EvenPowers}, the multiplicities associated to each of the points of
$\spec(A)$ lying on the real line are even, and the last part of
Definition~\ref{PoleQuotient} makes sense: there are well defined spaces
$\Dom_{\sigma,\frac 12}(A)$ for each $\sigma \in \spec_b(A)\cap \set{\Im\sigma=0}$. 

\begin{theorem}\label{FriedExtensionConst}
Let $A \in x^{-\nu}\diff^m_b(M;E)$ be $b$-elliptic, symmetric, 
bounded from below, and such that $P$ has coefficients independent of $x$ for
$x$ small. Suppose $\sigma\in\spec_b(A) \Longrightarrow \Im\sigma=0$ or
$|\Im\sigma|>\nu/2$. Then the domain 
of the Friedrichs extension of $A$ is given by
\begin{equation*}
 \Dom_F(A)=
 \sum_{\substack{\Im\sigma=0\\ \sigma\in\spec_b(A)}}
 \hspace{-.8em}\Dom_{\sigma,\frac 1{2}}(A).
\end{equation*}
\end{theorem}
In the situation of the theorem, the spaces $\Dom_{\sigma,\frac
1{2}}(A)/\Dom_{\min}(A)$ agree via the Mellin transform, with those defined in
Proposition~\ref{HalfDomains} since $P$ has coefficients independent of $x$ near the
boundary.
\begin{proof} With the hypotheses of the proposition,  
 \begin{equation*} \mathcal E(A)=\bigoplus_{\substack{\Im\sigma=0\\
\sigma\in\spec_b(A)}}
 \hspace{-.8em} \mathcal
 E_\sigma(A),
 \end{equation*}
where $\mathcal E_\sigma(A)$ was defined in \ref{PoleQuotient}. Let
\begin{equation*}
\mathcal E_F(A) =\Dom_F(A)/\Dom_{\min}(A),
\end{equation*}
a subspace of $\mathcal E(A)$. Passing to the Mellin transform side, $\hat{\mathcal
E}_F$ is saturated since $\Dom_F$ and $\Dom_{\min}(A)$ are $\kappa$-invariant, and
selfadjoint in $\hat{\mathcal E}(A)$ in the sense of the appendix of
Section~\ref{CanonicalPairing} since
$A$ with domain $\Dom_F$ is selfadjoint.  Since $\hat{\mathcal E}_F$ is saturated, by
Lemma~\ref{saturation2} there are saturated subspaces $\hat{\mathcal E}_{\sigma_j,F}
\subset \hat{\mathcal E}_{\sigma_j}(A)$ such that 
\begin{equation*}
\hat{\mathcal E}_F=\bigoplus_{\sigma_j\in S} \hat{\mathcal E}_{\sigma_j,F}.
\end{equation*}
Since $\hat{\mathcal E}_F$ is selfadjoint in $\hat{\mathcal E}$,
since $[\cdot,\cdot]^\flat_A$ is
nondegenerate, and since $[\hat u,\hat v]^\flat_A=0$ if $\hat u \in \hat{\mathcal
E}_{\sigma_j}(A)$ and $v \in \hat{\mathcal E}_{\sigma_k}(A)$ with $\sigma_j \ne
\sigma_k$ (Proposition~\ref{BasicPairing}), each $\hat{\mathcal E}_{\sigma_j,F}$ is
selfadjoint in $\hat{\mathcal E}_{\sigma_j}$. Moreover, because $A$ is bounded from
below,
$\mathcal P = \hat P_0$ is nonnegative by Lemma~\ref{Positivity}, and so, by
Lemma~\ref{EvenPowers} the $\mu_j$ are even. Now Proposition~\ref{HalfDomains}
applies and we deduce 
\begin{equation*}
\hat{\mathcal E}_{\sigma_j,F}=\hat{\mathcal E}_{\sigma_j,\frac 12}.
\end{equation*}
But by definition, $\Dom_{\sigma_j,\frac 12}$ is the space of elements on
$\Dom_{\max}$ such that $\hat u$ represents an element in $\hat{\mathcal
E}_{\sigma_j,\frac 12}$. 
\end{proof}

\begin{lemma}\label{Fried5}
Let $A\in x^{-\nu}\diff^m_b(M;E)$ be $b$-elliptic, symmetric and bounded
from below, and let $P=x^{\nu}A$. Then 
\begin{equation*}
 \Dom_{\max}(A)\cap\Dom_{F}(x^{-\eps}P) \subset \Dom_{F}(A)
 \ \text{ for any positive } \eps<\nu.
\end{equation*}
\end{lemma}
\begin{proof}
Note that
\begin{equation}\label{EquivPairings}
 (x^{-\nu}Pv,v)_{x^{-\nu/2}L^2_b}=
 (Pv,v)_{L^2_b}=(x^{-\eps}Pv,v)_{x^{-\eps/2}L^2_b}
\end{equation}
whenever all three expressions exist. It is always true if
e.g.~$v\in C_c^{\infty}(\Dot M;E)$. This shows, in particular, that $A=x^{-\nu}P$
is symmetric in $x^{-\nu/2}L^2_b$ if and only if $x^{-\eps}P$ is
symmetric in $x^{-\eps/2}L^2_b$. Suppose $\eps<\nu$. Then
$x^{-\eps}P$ is also bounded from below because
\begin{equation}\label{embedding}
 x^{-\eps/2}L^2_b(M;E)\embed x^{-\nu/2}L^2_b(M;E).
\end{equation}
Let $u\in\Dom_{\max}(A)\cap\Dom_{F}(x^{-\eps}P)$ and let $\seq{u_n}\in C_c^\infty(\Dot
M;E)$ such that
\begin{equation*}
 c\|u-u_n\|^2_{x^{-\eps/2}L^2_b}+ (x^{-\eps}P(u-u_n),u-u_n)_{x^{-\eps/2}L^2_b}\to 0
 \ \text{ as } n\to\infty.
\end{equation*}
Then $(x^{-\nu}P(u-u_n),u-u_n)_{x^{-\nu/2}L^2_b} =
(x^{-\eps}P(u-u_n),u-u_n)_{x^{-\eps/2}L^2_b} \to 0$ because of \eqref{EquivPairings}.
Since also $\|u-u_n\|_{x^{-\nu/2}L^2_b}\to 0$ because
$\|u-u_n\|_{x^{-\eps/2}L^2_b}\to 0$, the lemma is proved.
\end{proof}

\begin{lemma}\label{Fried9}
Let $A=x^{-\nu}P$ be as in Lemma~\ref{Fried5}. Then $P$ can be written as
\begin{equation*}
 P=P_0 + xP_1
\end{equation*}
with $P_0$, $P_1\in\diff^m_b(M;E)$ such that $x^{-\nu}P_0$ is
$b$-elliptic, symmetric, bounded from below, and has coefficients independent
of $x$ for $x$ small.
\end{lemma}
\begin{proof}
Near the boundary $\partial M$, $P$ can be written as $P=\tilde P_0 +
x\tilde P_1$, where $\tilde P_0$ has coefficients independent of $x$.
Let $\omega\in C_c^\infty(\R)$ be equal to $1$
near $\partial M$, define  
\begin{equation*}
 P_0=\omega\tilde P_0\,\omega + (1-\omega)P(1-\omega).
\end{equation*}
Clearly $(1-\omega)P(1-\omega)$ is symmetric and bounded from below. As the conormal
symbol of $P$, $\mathcal P=\widehat{\tilde P}_0$ is a selfadjoint
holomorphic family in the sense that $\mathcal P(\sigma)^* =\mathcal P(\overline
\sigma)$ on $H^m(\partial M;E|_{\partial M})$, and positive by
Lemma~\ref{Positivity}. From the Mellin transform version of Plancherel's identity it
follows that $\omega\tilde P_0\,\omega$ is also bounded from below if $\omega$ has 
sufficiently small support. Evidently, if the support of
$\omega$ is small enough, then $P_0$ is elliptic in the interior and therefore
$b$-elliptic. 
\end{proof}

\begin{lemma}\label{Fried88}
Let $A=x^{-\nu}P$ with $P=P_0+xP_1$ as in Lemma~\ref{Fried9}. Then for $0<\eps < \nu$,
$x^{-\eps}P$ and $x^{-\eps}P_0$ are both $b$-elliptic, symmetric
and bounded from below as operators on $x^{-\eps/2}L^2_b(M;E)$, and for small
$\eps$, $\Dom_{F}(x^{-\eps}P) = \Dom_{F}(x^{-\eps}P_0)$.
\end{lemma}

The first statement follows from Lemma~\ref{Fried9} and the proof of
Lemma~\ref{Fried5}, and the equality of the Friedrichs domains is a consequence of
part~\ref{EdDom4} of Proposition~\ref{EdDom}.


\begin{theorem}\label{FriedExtension}
Let $A\in x^{-\nu}\diff^m_b(M;E)$ be $b$-elliptic, symmetric and bounded
from below. Then the domain of the Friedrichs extension of $A$ is 
\begin{equation*}
 \Dom_F(A)=\sum_{\substack{-\nu/2<\Im\sigma<0\\ \sigma\in\spec_b(A)}}
 \hspace{-1em}\mathcal D_{\sigma}(A) + \hspace{-1em}
 \sum_{\substack{\Im\sigma=0\\ \sigma\in\spec_b(A)}}
 \hspace{-.8em}\Dom_{\sigma,\frac 1{2}}(A).
\end{equation*}
\end{theorem}
\begin{proof} Let $\Dom_F'(A)$ be the space on the right in the statement. Then $A$
with domain $\Dom_F'(A)$ is selfadjoint, so we only need to prove that
$\Dom_F'(A)\subset \Dom_F(A)$, and we proceed to do so. Write
$x^\nu A=P_0+xP_1$ as in Lemma \ref{Fried9}. From Lemmas \ref{Fried5} and
\ref{Fried88}  we get that if $\eps>0$ is small enough then 
$\Dom_{\max}(A)\cap\Dom_{F}(x^{-\eps}P_0) \subset \Dom_{F}(A)$. We may apply Theorem
\ref{FriedExtensionConst} to $x^{-\eps}P_0$ and deduce that
\begin{equation*}
 \Dom_F(x^{-\eps} P_0)=
 \sum_{\substack{\Im\sigma=0\\ \sigma\in\spec_b(A)}}
 \hspace{-.8em}\Dom_{\sigma,\frac 1{2}}(x^{-\eps} P_0).
\end{equation*}
But the intersection of this space and $\Dom_{\max}(A)$ is $\Dom_F'(A)$. Thus
$\Dom_F'(A)\subset \Dom_F(A)$ and therefore $\Dom_F'(A) =\Dom_F(A)$. 
\end{proof}

Together with Theorem~\ref{Fried7} we in particular obtain
\begin{equation*}
 \Dom_F=\Dom_{\max}(A)\cap H^m_b(M;E)\ \text{ if and only if }
 \spec_b(A)\cap\set{\Im\sigma=0}=\varnothing.
\end{equation*}

The following corollary improving part \ref {EdDom4} of Proposition \ref{EdDom} is an
immediate consequence of the theorem, since the hypothesis implies that for $-\nu\leq
\Im\sigma\leq 0$ the spaces $\Dom_\sigma$ for both operators are equal:

\begin{corollary}
Suppose $A_0$, $A_1 \in x^{-\nu}\diff^m_b(M;E)$ are $b$-elliptic, symmetric and
bounded from below. If $A_0-A_1$ vanishes to order $k$, $k>\nu/2$, then $\Dom_F(A_0)=
\Dom_F(A_1)$.
\end{corollary}

\section{Applications and Examples}\label{Examples}

Let $A\in x^{-\nu}\diff_b^m(M;E)$
be $b$-elliptic and assume that for any two distinct $\sigma_0$ and 
$\sigma_1$ in $\Sigma(A)$, $\Im(\sigma_0-\sigma_1) \not\in \mathbb Z$.
Using that $[\cdot,\cdot]_A$ pairs $\mathcal E_{\sigma_0}(A)$ and $\mathcal
E_{\overline \sigma_0}(A^\star)$ nonsingularly, one can extend the examples 
to the excluded case. Recall that $A^\star$ denotes the formal adjoint of $A$,
and $\Sigma(A)=\spec_b(A)\cap\set{-\nu/2<\Im\sigma<\nu/2}$.

\begin{example}
Regard $A$ with domain $\Dom_{\sigma_0}(A)$ for some $\sigma_0\in \Sigma(A)$.
Then,  
\begin{equation*}
 \Dom(A^*)=\sum_{\substack {\sigma \in \Sigma(A)\\ \sigma\ne \sigma_0}}
 \Dom_{\overline \sigma}(A^\star),
\end{equation*}
where $A^*$ denotes the Hilbert space adjoint of $A$.
If there there are one or more $\sigma_j\in \Sigma(A)$ such that
$\sigma_0-\sigma_j=\im\tau$ with $\tau\in \mathbb N$, then in place of the
$\Dom_{\overline \sigma_j}(A^\star)$ one must use certain subspaces. 
\end{example}

\begin{example} 
If $A$ is given the domain 
\begin{equation*}
 \sum_{\substack{\sigma
 \in \Sigma(A)\\ \Im\sigma< 0}}\Dom_{\sigma}(A),
\end{equation*}
then the domain of the adjoint is
\begin{equation*}
 \Dom(A^*)=\sum_{\substack
 {\sigma \in \Sigma(A)\\ \Im\sigma\geq  0}}\Dom_{\overline\sigma}(A^\star).
\end{equation*} 
Note that the poles of the Mellin transforms of elements in this space are on
or below the real axis.
\end{example}

\begin{example}
Suppose  $\spec_b(A)$ does contain points on $\Im\sigma = 0$, but that all the
partial multiplicities $\mu_{\sigma_s,j}$ of each such point are even. Referring to
Definition~\ref{PoleQuotient} for the notation, let $A$ have domain
\begin{equation*}
\sum_{\substack{\sigma\in \Sigma(A)\\ \Im\sigma<0}}\Dom_{\sigma}(A) +
\sum_{\substack{\sigma\in \Sigma(A)\\ \Im\sigma=0}}\Dom_{\sigma,\frac 12}(A)
\end{equation*}
Then,
\begin{equation*}
\Dom(A^*)= \sum_{\substack {\sigma \in \Sigma(A)\\ \Im\sigma> 0}}
 \Dom_{\overline\sigma}(A^\star) +
 \sum_{\substack {\sigma \in \Sigma(A)\\ \Im\sigma= 0}}
 \Dom_{\sigma,\frac 12}(A^\star).
\end{equation*}
Thus if $A$ is symmetric then $A$ with the given domain is selfadjoint. 
In Section~\ref{FriedrichsExtension} we proved that if $A$ is symmetric and
bounded from below then this is the domain of the Friedrichs extension.
\end{example}

In some cases, in particular geometric problems, one encounters
operators of the form $B^*B$ and $BB^*$. Below we discuss some aspects of
$B^*B$ using some of the results of this paper. The operator $BB^*$ can be
treated in the same manner.

Assume that $B\in x^{-\nu}\diff_b^m(M)$ is $b$-elliptic, $B=x^{-\nu}Q$. Let
$\Dom(B)\subset\Dom_{\max}(B)$ be such that $B:\Dom(B)\to x^{-\nu/2}L^2_b(M)$ 
is closed. Then $B^*B$ is a selfadjoint extension of the symmetric operator
$B^\star B$, considered as an unbounded operator on $x^{-\nu/2}L^2_b(M)$. 
Recall that
\begin{equation*}
 \Dom(B^*B) =\set{u\in\Dom(B)\st Bu\in \Dom(B^*)} 
\end{equation*}
Since $B^*$ is a closed extension of the formal adjoint $B^\star$, $B^*B$ is indeed a
closed extension of $B^\star B$ with
\begin{equation*}
 \Dom_{\min}(B^\star B)\subset \Dom(B^*B)\subset \Dom_{\max}(B^\star B).
\end{equation*}
Note that if $u\in\Dom_{\max}(B^\star B)$, then $\hat u$ is meromorphic in
$\set{\Im\sigma>-3\nu/2}$ with poles on the strip
$\set{\nu/2>\Im\sigma>-3\nu/2}$. 
Write $B^\star B=x^{-2\nu} P$ with $P\in\diff_b^{2m}(M)$.
The conormal symbol of $B^\star B$ is then given by
\begin{equation*}
  \hat P_0(\sigma)= \hat Q_0(\overline{\sigma}-\im\nu)^* 
  \circ \hat Q_0(\sigma),
\end{equation*}
where $\hat Q_0(\sigma)$ is the conormal symbol of $B$. Thus, the boundary
spectrum of $B^\star B$ contains $\spec_b(B)$ and its reflection with respect
to $\set{\Im\sigma=-\nu/2}$ (line of symmetry). In particular, every
$\sigma\in\spec_b(B^\star B)\cap \set{\Im\sigma=-\nu/2}$ has even
multiplicities. For $\sigma_0\in\spec_b(B)$ define 
$\Dom^s_{\sigma_0}(B)$ as the space of elements $u\in\Dom_{\max}(B)$ such that
$Bu\in\Dom_{\max}(B^\star)$, and such that $\hat u$ is meromorphic in $\C$
with poles at most at $\sigma_0-\im\vartheta$ for
$\vartheta=0,\dots,\ceil{\nu}$. Thus, for $\Im\sigma_0=-\nu/2$, we have
$\Dom^s_{\sigma_0}(B)\subset\Dom_{\min}(B)$.

\begin{lemma}
If $u\in\Dom^s_{\sigma_0}(B)$ and $\Im\sigma_0=-\nu/2$, then 
$Bu\in\Dom_{\min}(B^\star)$. 
\end{lemma}
\begin{proof}
Let $\sigma_0\in\spec_b(B)$ be such that $\Im\sigma_0=-\nu/2$. If
$u\in\Dom^s_{\sigma_0}(B)$, then $Bu\in x^{-\nu/2}L^2_b$ and 
$\hat u$ is meromorphic with poles at $\sigma_0-\im\vartheta$ for
$\vartheta=0,\dots,\ceil{\nu}$. It follows that $\widehat{Bu}$ is then
holomorphic in $\Im\sigma>-\nu/2$ which implies
$Bu\in\Dom_{\min}(B^\star)$.
\end{proof}

Since $\Dom_{\min}(B)\subset\Dom(B)$ and
$\Dom_{\min}(B^\star)\subset\Dom(B^*)$, the previous lemma implies 
\begin{lemma}
If $\Im\sigma_0=-\nu/2$ then $\Dom^s_{\sigma_0}(B)\subset \Dom(B^*B)$.
\end{lemma}
Note that $\Dom^s_{\sigma_0}(B)=\Dom_{\sigma_0,\frac 12}$.

\begin{example}
If $\Dom(B)=\Dom_{\min}(B)$, $B^*B$ is the Friedrichs extension of
$B^\star B$ and
\begin{equation}\label{B*BMin}
 \Dom(B^*B)=\set{u\in\Dom_{\min}(B)\st Bu\in\Dom_{\max}(B^\star)}.
\end{equation}
Denote the set on the right by $\Dom_{F}$. Since $\Dom(B^*B)\subset\Dom_{F}$,
then $\Dom_{F}^\perp \subset \Dom(B^*B)^\perp=\Dom(B^*B)$,
where $\perp$ means the orthogonal in $\Dom_{\max}(B^\star B)$ with
respect to $[\cdot,\cdot]_{B^\star B}$.
Now let $u\in \Dom(B^*B)$, so $Bu\in\Dom(B^*)$. Then, for every $v\in\Dom_{F}$
\begin{equation*}
 0= [v,Bu]_B= (Bv,Bu) - (v,B^\star B u)= (B^\star Bv,u) - (v,B^\star B u) 
  = [v,u]_{B^\star B}.
\end{equation*}
This implies $u\in \Dom_{F}^\perp$ and we get \eqref{B*BMin}. In
this case, the Mellin transform $\hat u$ of an element $u\in\Dom(B^*B)$ is
holomorphic in $\set{\Im\sigma>-\nu/2}$ and meromorphic in
$\set{\Im\sigma>-3\nu/2}$ with poles on 
$\set{\overline\sigma-\im\nu \st \sigma\in \spec_b(B)}$.

If $\Dom(B)=\Dom_{\max}(B)$, then
\begin{equation*}\label{B*BMax}
 \Dom(B^*B)=\set{u\in\Dom_{\max}(B)\st Bu\in\Dom_{\min}(B^\star)}
\end{equation*}
and for $u\in\Dom(B^*B)$, $\hat u$ has poles at most on
$\spec_b(B^\star)\cap\set{-\nu/2\le \Im\sigma <\nu/2}$.

If $\Dom(B)=\Dom_{\sigma_0}(B)$ for some
$\sigma_0\in\Sigma(B)=\spec_b(B)\cap\set{-\nu/2< \Im\sigma< \nu/2}$, then
\begin{equation*}
\Dom(B^*B)=\set{u\in\Dom_{\sigma_0}(B)\st Bu\in
\mbox{$\sum_\sigma$}\Dom_{\sigma}(B^\star) \text{ for }
\sigma \in \Sigma(B^\star),\; \sigma\ne\overline\sigma_0}
\end{equation*}
and for $u\in\Dom(B^*B)$, $\hat u$ has poles at $\sigma_0$ and on
$\set{\overline\sigma-\im\nu \st \sigma\in\spec_b(B),\; \sigma\ne \sigma_0}$.
\end{example}

We finish this section with some concrete examples.
Assume now that, near the boundary, $M$ is of the form $[0,1)\times S^n$ with
$\partial M=\set{0}\times S^n$. Let
$D_y^2$ be the Laplacian on $S^n$, let $x$ be the variable in $[0,1)$. For
$0<\nu\le 2$ let $A=x^{-\nu}P$ with
\begin{equation}\label{ExGeneric}
 P = (xD_x)^2 + a^2 D_y^2 +\beta b^2,
\end{equation}
$\beta=\pm 1$ and nonnegative constants $a$ and $b$. They arise, for
instance, as Laplace-Beltrami operators on scalar functions associated to metrics
of the form 
\[ g=dx^2 /{x^{2-\nu}} + x^{\nu}dy^2. \] 
Compare Br\"uning and Seeley \cite{BrSe85}, Lesch \cite{Le97}, and especially Mooers
\cite{Moo99} on $k$-forms.

The conormal symbol of $A$ is defined to be
\begin{equation*}
 \widehat P_0(\sigma) = \sigma^2 + a^2 D_y^2 +\beta b^2,
\end{equation*}
and the boundary spectrum $\spec_b(A)$ is the set of points in $\C$ such that
$\hat P(\sigma)$ is not invertible. Since the set of eigenvalues of $D_y^2$ is
$\set{k(k+n-1) \st k\in\mathbb N_0}$, then 
\begin{equation*}
 \spec_b(A)=\set{\sigma\in\C \st \sigma^2+a^2 k(k+n-1) +\beta b^2=0
  \;\text{ for some } k\in\mathbb N_0} 
\end{equation*}
Consider $A:C_c^\infty(M)\subset x^{-\nu/2} H^2_b(M)\to x^{-\nu/2}H^2_b(M)$
as an unbounded operator. Then 
 $\Sigma=\spec_b(A)\cap \set{-\nu/2\le \Im\sigma\le \nu/2}$ 
is the only set that matters when looking for the closed extensions 
of $A$. We look separately at the cases $b=0$ and $b\not=0$.

\begin{example}\label{CEx1}
Let $b=0$ and $a>0$ in \eqref{ExGeneric}. In this case, $A$ is symmetric
and bounded from below. Note that $\sigma^2+a^2 k(k+n-1)$ has simple roots in
$\Sigma$ when $k\not=0$, and a root of order $2$ when $k=0$. For 
illustration purposes it is enough to look at the case when $\nu=2$ and $n=1$.
Thus $\Sigma=\set{\pm\im ka \st ka\le 1, k\in\mathbb N_0}$. 
Suppose $a > 1$ so that the only point in $\Sigma$ is $0$. Then
$\Dom_{\min}(A)=x H^2_b(M)$ and 
\begin{equation*}
\mathcal E(A)=\Dom_{\max}(A)/\Dom_{\min}(A)=\LinSpan\set{\omega,\omega\log x}
\end{equation*}
for some $\omega\in C_c^\infty(\R)$ such that $\omega(x)=1$ near the origin.

Let $\psi_0=\omega$ and $\psi_1=\omega\log x$. Then,
\begin{gather*}
 \hat \psi_0(\sigma)= \frac {\Phi(\sigma)}{\sigma} \quad\text{and}\quad
 \hat \psi_1(\sigma)= \frac {\Phi(\sigma)}{\sigma^2}-
 \frac{\Phi'(\sigma)}{\sigma}
\intertext{with}
 \Phi(\sigma) = \int_{\R^+} x^{-\im\sigma}D_x\omega(x)\,dx.
\end{gather*}
Let $u=u_0\psi_0+u_1\psi_1$ and $v=v_0\psi_0+v_1\psi_1$ be elements of 
$\mathcal E(A)$. Then, with $\gamma$ a positively oriented simple
closed curve in $\C$ surrounding $0$,
\begin{align*}
[u,v]_A &= \frac 1{2\pi} \oint_\gamma \sigma^2 \hat u(\sigma)
 \overline {\hat v(\overline \sigma)}\,d\sigma\\
 &=\frac 1{2\pi}\oint_\gamma 
   \Big(\frac{u_0\overline v_1 + u_1\overline v_0}{\sigma}\Big)
   \Phi(\sigma)\overline{\Phi(\overline\sigma)} d\sigma\\
 &= \im (u_0\overline v_1 + u_1\overline v_0). 
 \end{align*}
If $\Dom$ is a domain on which $A$ is selfadjoint, then
$\Dom=\Dom_{\min}(A)\oplus\, \LinSpan{u}$ with some $u$ as above such that
$[u,u]_A=0$. Thus $u_0\overline u_1+u_1\overline u_0=0$ which
implies $u_0\overline u_1\in \im\R$, so $\LinSpan{u}$ is one dimensional. 
Moreover, all selfadjoint extensions of $A$ are of the form
$A_\lambda:\Dom^\lambda\to x^{-1} H^2_b(M)$ with
\begin{equation*}
 \Dom^\lambda/\Dom_{\min}(A)= 
 \LinSpan\set{(e^{\im\lambda}+1)\omega +(e^{\im\lambda}-1)\omega\log x}.
\end{equation*}
In particular, $\Dom^0$ is precisely the domain of the Friedrichs
extension of $A$. 
If $\frac 1{\ell+1}< a\le \frac 1{\ell}$, $\ell\in\mathbb N$, then 
$\Sigma=\set{\pm \im ka\st k=0,1,\dots,\ell}$ and we have
\begin{equation*}
\Dom_F(A)=\Dom^0\oplus\sum_{k=1}^\ell\Dom_{\set{-\im ka}}(A). 
\end{equation*}
Note that the partial multiplicities of the poles $\im ka$, $k\ne 0$ are equal to
$1$, but the total multiplicity is $2$. 
\end{example}

\begin{example}
Let $\alpha\in C_c^\infty(\R)$ such that $\alpha(0)>1$. Then
\begin{equation*}
 A = x^{-2}\left[(xD_x)^2 + \alpha(x)^2 D_y^2\right]  
\end{equation*}
is symmetric and bounded from below.  Moreover, $0$ is the only point of the
boundary spectrum in $\set{-1\leq \Im \sigma\leq 1}$. If 
$A_0=x^{-2}[(xD_x)^2 + \alpha(0)^2 D_y^2]$, then
\begin{equation*}
 \Dom_{\min}(A)=\Dom_{\min}(A_0),\;\; \Dom_{\max}(A)=\Dom_{\max}(A_0),
 \;\text{ and }\; \Dom_{F}(A)=\Dom_{F}(A_0).
\end{equation*}
\end{example}

\begin{example}
Let $b\not=0$. If $\beta=1$ in \eqref{ExGeneric},
then $A=x^{-2}P$ behaves similarly but `nicer'
than the operator in Example~\ref{CEx1} since $\sigma^2+a^2 k(k+n-1)+b^2$ 
has only simple roots, and $\spec_b(A)\cap \set{\Im\sigma=0}=\varnothing$.

If $\beta=-1$, the situation is different. 
In this case, the conormal symbol of $A$ 
\begin{equation*}
 \hat P(\sigma) = \sigma^2 + a^2 D_y^2 - b^2
\end{equation*}
fails to be nonnegative for real $\sigma$ and therefore $A$ is not bounded
from below. Moreover, $\sigma^2+a^2 k(k+n-1)-b^2$ has two simple real roots,
$-b$ and $b$. We assume $b\le 1< a$, so $\spec_b(A)\cap 
\set{-1\leq \Im\sigma \leq 1}$ contains only $-b$ and $b$. Thus 
\begin{equation*}
\mathcal E(A)=\Dom_{\max}(A)/\Dom_{\min}(A)=
  \LinSpan\set{\omega x^{\im b},\omega x^{-\im b}}
\end{equation*}
for some $\omega\in C_c^\infty(\R)$ such that $\omega(x)=1$ near the origin.

Let $\psi_+ = \omega x^{\im b}$ and $\psi_{-} = \omega x^{-\im b}$. Then,
with $\Phi$ as above, 
\[ \hat \psi_\pm(\sigma)= \frac{\Phi(\sigma\mp b)}{\sigma\mp b}.\]
Let $u=u_+\psi_+ + u_-\psi_-$, $v=v_+\psi_+ + v_-\psi_-$. Then, with $\gamma$ 
a positively oriented simple curve in $\C$ surrounding $b$ and $-b$,
\begin{align*}
 [u,v]_A &=\frac 1{2\pi} \oint_{\gamma} (\sigma^2 - b^2)
          \hat u(\sigma)\overline {\hat v(\overline \sigma)}\,d\sigma\\       
 &=\frac 1{2\pi}\oint_\gamma \Big\{
   \left(\tfrac{\sigma+b}{\sigma-b}\right) u_+\overline v_+
   \Phi(\sigma-b)\overline{\Phi(\overline \sigma-b)}\\
 &\hspace{8em}
   +\left(\tfrac{\sigma-b}{\sigma+b}\right) u_-\overline v_-
   \Phi(\sigma+b)\overline{\Phi(\overline\sigma+b)}\Big\}d\sigma\\
 &=2b\im (u_+\overline v_+ - u_-\overline v_-)
 \end{align*}
Thus $[u,u]_A = 2b \im (|u_+|^2 -|u_-|^2)$ and if $A$ with
domain $\Dom_{\min}(A)\oplus\,\LinSpan u$ is selfadjoint then $|u_+|=|u_-|$.
Thus $\LinSpan u$ is one dimensional and all selfadjoint extensions of $A$
are of the form $A_\lambda:\Dom^\lambda\to x^{-1}H^2_b(M)$ with
\begin{equation*}
 \Dom^\lambda/\Dom_{\min}(A)= 
 \LinSpan\set{\omega x^{\im b}+e^{\im \lambda}\omega x^{-\im b}}.
\end{equation*}
\end{example}


\begin{thebibliography}{10}

\bibitem{BrSe85} 
J.~Br\"uning and R.~Seeley, \emph{Regular singular asymptotics}, Adv. in Math.
\textbf{58} (1985), 133--148.

\bibitem{BrSe88} 
\bysame, \emph{An index theorem for first order regular singular operators},
Amer. J. Math. \textbf{110} (1988), 659--714.

\bibitem{Ch79} 
J.~Cheeger, \emph{On the spectral geometry of spaces with cone-like singularities},
Proc. Nat. Acad. Sci. USA \textbf{76} (1979), 2103--2106.

\bibitem{GohSi} 
I.~C. Gohberg and E.~I. Sigal, \emph{An operator generalization of the logarithmic
residue theorem and {R}ouch\'e's theorem}, Math. USSR Sbornik \textbf{13}
(1971), no. 4, 603--625.

\bibitem{K1} 
V.A. Kondrat'ev, \emph{Boundary problems for elliptic equations in domains with
conical or angular points}, Trans. Mosc. Math. Soc. \textbf{16} (1967), 227--313.

\bibitem{JM}
H. Jacobowitz and G. Mendoza, \emph{Elliptic equivalence of vector bundles},
preprint, 2000.

\bibitem{Le97} 
M.~Lesch, \emph{Operators of {F}uchs type, conical singularities, and asymptotic
methods}, Teubner-Texte zur Math. vol 136, B.G. Teubner, Stuttgart, Leipzig, 1997.

\bibitem{RBM1} 
R.~Melrose, \emph{Transformation of boundary value problems}, Acta Math. \textbf{147}
(1981), no. 3-4, 149--236.

\bibitem{RBM2}
\bysame, \emph{The Atiyah-Patodi-Singer index theorem}, 
Research Notes in Mathematics, A~K~Peters, Ltd., Wellesley, MA, 1993.

\bibitem{MM} 
R.~Melrose and G.~Mendoza, \emph{Elliptic operators of totally characteristic type},
MSRI Preprint, 1983.

\bibitem{Moo99} 
E.~Mooers, \emph{Heat kernel asymptotics on manifolds with conic singularities},
J. Anal. Math. \textbf{78} (1999), 1--36.

\bibitem{RiNa} 
F.~Riesz and B.~Sz.-Nagy, \emph{Functional analysis}, Dover Publications Inc., New York,
1990. Translated from the second French edition by Leo F. Boron,
reprint of the 1955 original.

\bibitem{CSS01} 
S. Coriasco,  E. Schrohe, and J. Seiler, \emph{Bounded imaginary powers of
differential operators on manifolds with conical
singularities}, preprint 2001 (math AP/0106008).

\bibitem{Sz91} 
B.-W.~Schulze,  \emph{Pseudo-differential operators on manifolds with singularities},
North-Holland, Amsterdam, 1991.

\bibitem{Sz94} 
\bysame, \emph{Pseudo-differential boundary value problems, conical singularities,
and asymptotics}, Akademie Verlag, Berlin, 1994.

\bibitem{Wlo} 
J.~Wloka, \emph{Partial differential equations}, Cambridge University Press, Cambridge,
1987. Translated from the German by C. B. Thomas and M. J. Thomas.

\end{thebibliography}
\end{document}